\documentclass{amsart}
\usepackage{amsmath,amscd}
\usepackage{amsfonts,amssymb}
\usepackage[all]{xy}
\usepackage{enumerate}
\numberwithin{equation}{subsection}
\theoremstyle{plain}
\newtheorem{thm}[subsection]{Theorem}
\newtheorem{prop}[subsection]{Proposition}
\newtheorem{lemma}[subsection]{Lemma}
\newtheorem{cor}[subsection]{Corollary}
\newtheorem*{prin}{Homotopy Section Principle (HSP)}
\theoremstyle{definition}
\newtheorem{defn}[subsection]{Definition}
\newtheorem{example}[subsection]{Example}
\newtheorem{notn}[subsection]{Notation}
\newtheorem{cont}[subsection]{Contents}
\newtheorem{ackn}[subsection]{Acknowledgement}
\theoremstyle{remark}
\newtheorem{rem}[subsection]{Remark}

\input cyracc.def
\font\tencyr=wncyr10
\font\eightcyr=wncyr7
\def\cyr{\tencyr\cyracc}
\def\cyri{\eightcyr\cyracc}
\def\Sh{\textrm{{\cyr Sh}}}
\def\Shi{\textrm{{\cyri Sh}}}
\begin{document}
\title{\'Etale homotopy equivalence of rational points on algebraic varieties}
\author{Ambrus P\'al}
\date{February 10, 2015.}
\address{Department of Mathematics, 180 Queen's Gate, Imperial College, London, SW7 2AZ, United Kingdom}
\email{a.pal@imperial.ac.uk}
\begin{abstract} It is possible to talk about the \'etale homotopy equivalence of rational points on algebraic varieties by using a relative version of the \'etale homotopy type. We show that over $p$-adic fields rational points are homotopy equivalent in this sense if and only if they are \'etale-Brauer equivalent. We also show that over the real field rational points on projective varieties are \'etale homotopy equivalent if and only if they are in the same connected component. We also study this equivalence relation over number fields and prove that in this case it is finer than the other two equivalence relations for certain generalised Ch\^atelet surfaces. 
\end{abstract}
\footnotetext[1]{\it 2000 Mathematics Subject Classification. \rm 14F35, 
14G05.}
\maketitle
\pagestyle{myheadings}
\markright{\'Etale homotopy equivalence of rational points on algebraic varieties}

\section{Introduction}

For every field $K$ let $\overline K$ denote its separable closure. For every variety $X$ defined over a $K$ as above and for every field extension $L|K$ let $X_L$ denote the base change of $X$ to $L$. Let $K$ be a field and let $X$ be a variety defined over $K$. In the paper \cite{HS} Harpaz and Schlank define a relative version $Et_{/K}(X)$ of the \'etale homotopy type of $X$ by looking at the action of the absolute Galois group Gal$(\overline K|K)$ of $K$ on the \'etale hypercoverings of $X_{\overline K}$. With the aid of the action of Gal$(\overline K|K)$ on $Et_{/K}(X)$ they define a pro-object $X(hK)$ in the category of sets, which they call the homotopy fixed points set of $X$, that serves as a certain homotopical approximation of the set $X(K)$ of rational points. By slight abuse of notation we will use the same symbol to denote the projective limit of $X(hK)$ which we will consider as a topological space equipped with its pro-discrete topology. It is possible to define a natural map: 
$$\iota_{X/K}:X(K)\longrightarrow X(hK)$$
which can be thought of as a homotopy-theoretical version of the section map in Grothendieck's anabelian geometry, and which it also happens to refine. 
We say that $x,y\in X(K)$ are $H$-equivalent if $\iota_{X/K}(x)=\iota_{X/K}(y)$. 

The aim of this paper is to describe the $H$-equivalence relation on $X(K)$ as explicitly as possible for many $K$ and $X$. Let us first turn to the case when $K$ is a finite extension of $\mathbb Q_p$. In this case the map $\iota_{X/K}$ is not surjective in general; every abelian variety of positive dimension is a counterexample (see Proposition \ref{section}). However it is possible to describe the equivalence relation it induces on $X(K)$ in rather concrete terms. For any smooth variety $X$ over any field $K$ of characteristic zero let Br$(X)=H^2(X,\mathbb G_m)$ denote the cohomological Brauer group of $X$. We say that $x,y\in X(K)$ are Brauer equivalent if $x^*(b)=y^*(b)$ for all $b\in\textrm{Br}(X)$. We say that $x,y\in X(K)$ are \'etale--Brauer equivalent if for every finite, \'etale morphism $Y\rightarrow X$ of varieties over $K$ and for each $\widetilde x\in Y(K)$ mapping to $x$ there is a $\widetilde y\in Y(K)$ which maps to $y$ and which is Brauer equivalent to $\widetilde x$. Then we have the following
\begin{thm}\label{etalebrauer} Let $K$ be a finite extension of $\mathbb Q_p$ and let $X$ be a smooth quasi-projective variety over $K$. Then \'etale--Brauer equivalence and $H$-equivalence coincide on $X(K)$.
\end{thm}
It is important to note that this claim is not true for more general fields;
Ch\^atelet surfaces over number fields provide counterexamples (see the remark following Theorem \ref{chatelet} below). The main ingredients of the theorem above, besides obstruction theory, is duality for the Galois cohomology of $K$ and Gabber's theorem on the existence of Azumaya algebras. We also provide examples to show that the theorem cannot be strengthened by substituting Brauer equivalence for \'etale--Brauer equivalence; see Theorem \ref{bielliptic} below. We can also characterise $H$-equivalence for the real number field:
\begin{thm}\label{real} Let $K$ be the real number field $\mathbb R$ and let $X$ be a smooth affine or projective variety over $K$. Then two $K$-rational points of $X$ are $H$-equivalent if and only if they are in the same connected component of the topological space $X(K)$.
\end{thm}
The main tools of the proof of this result are a celebrated theorem of Mah\'e (see \cite{HM} and \cite{Mah}), the theory of Stiefel-Whitney classes for quadratic bundles (see \cite{EKV} and \cite{Mil}), and an equivariant version of a basic comparison result of Artin-Mazur (see \cite{AM}).  The reader should note that in \cite{Qu} Quick developed a general theory of homotopy fixed point spaces for simplicial pro-sets equipped with a continuous action of a profinite group, which can be applied to Friedlander's \'etale topological type functor (see \cite{Fr}). His construction offers an alternative route for the foundations of our investigations.  
\begin{cont} In the next chapter we review the relative \'etale homotopy type and homotopy fixed points of varieties as defined by Harpaz and Schlank and their relation to the Artin-Mazur construction. In the third chapter we introduce a pointed version of the relative \'etale homotopy type and compare it with the previously defined constructions. In the fourth chapter we study the relationship between the \'etale homotopy groups of finite \'etale coverings. In the following chapter we show that the \'etale homotopy types of abelian varieties and smooth curves which are not projective of genus zero are Eilenberg-MacLane spaces over algebraically closed fields of characteristic zero. The fact presented in these two chapters might be well-known to the experts, but I could not find a convenient reference. In the sixth chapter we prove two useful lemmas about lifting a pair of rational points on certain principal bundles. Then we prove the fundamental theorem of obstruction theory for $H$-equivalence in the seventh chapter. We study the analogue of the Manin pairing for homotopy fixed points in the eighth chapter. In the ninth chapter we prove that \'etale--Brauer equivalence is strictly finer than Brauer equivalence on $X(K)$, when $K$ is a $p$-adic field and $X$ is a bielliptic surface  over $K$, using a rather standard set of tools. Theorem \ref{etalebrauer} is proved in the tenth chapter, while in the eleventh we prove Theorem \ref{real}. In the twelfth chapter we introduce a natural homotopy version of Gro\-thendieck's section and the Shafarevich-Tate conjectures over number fields by substituting the arithmetic fundamental group with the relative version of the \'etale homotopy type, which we call the Homotopy Section Principle (HSP), and prove that it is equivalent to its well-established analogues in the special case of curves and abelian varieties. We provide further examples of varieties which satisfy HSP (see Theorems \ref{bloch}, \ref{birational} and \ref{chatelet}) in the final two chapters, including generalised Ch\^atelet surfaces.
\end{cont}
\begin{ackn} This article was inspired by Bertrand To\"en's paper \cite{To} who, according to my best knowledge, was first to suggest the study of rational points via maps like $\iota_{X/K}$. I also wish to thank Jean-Louis Colliot-Th\'el\`ene, Yonathan Harpaz, Tomer Schlank, Alexei Skorobogatov and Akio Tamagawa either for useful discussions or for useful correspondence concerning the contents of this paper, and for the referee for his many comments and corrections. The author was partially supported by the EPSRC grants P19164 and P36794. 
\end{ackn}

\section{Basic definitions}

\begin{defn} Let $\Gamma$ be a profinite group. By a $\Gamma$-set we mean a set with a $\Gamma$-action such that each element has an open stabiliser. Let $\Gamma-Sets$ denote the category whose objects are $\Gamma$-sets and whose morphisms are $\Gamma$-equivariant maps between them. By a simplicial $\Gamma$-set we mean a simplicial object in $\Gamma-Sets$. These form a category $\Gamma-SSets$ in the usual way. Note that for every simplicial $\Gamma$-set $S$ and for every $n\in\mathbb N$ the $n$-skeleton $\textrm{sk}_n(S)$, the $n$-coskeleton $\textrm{cosk}_n(S)$ and the Kan replacement $\textrm{Ex}^{\infty}(S)$ are all naturally equipped with a $\Gamma$-action. Since with respect to this action the stabiliser of each simplex of $\textrm{sk}_n(S)$, $\textrm{cosk}_n(S)$ and $\textrm{Ex}^{\infty}(S)$ is open, these constructions furnish three functors, the $n$-skeleton $\textrm{sk}_n:\Gamma-SSets\rightarrow\Gamma-SSets$, the $n$-coskeleton $\textrm{cosk}_n:\Gamma-SSets\rightarrow\Gamma-SSets$, and the Kan replacement $\textrm{Ex}^{\infty}:\Gamma-SSets\rightarrow \Gamma-SSets$ functors. Moreover let $P_n:\Gamma-SSets\rightarrow\Gamma-SSets$ denote the corresponding analogue of the simplicial version of the $n$-th Postnikov piece given by the rule:
$$P_n(S) =\textrm{cosk}_{n+1}(\textrm{sk}_{n+1}(\textrm{Ex}^{\infty}S))$$
for every simplicial $\Gamma$-set $S$.
\end{defn}
\begin{notn} For every category $\mathcal C$ let Pro$-\mathcal C$ be the category of pro-objects of $\mathcal C$. For every pair of categories $\mathcal C$, $\mathcal D$ let $\mathcal C\times\mathcal D$ denote their direct product and for every category $\mathcal C$ let $\mathcal C^{op}$ denote its opposite category. Clearly there is a natural equivalence between $(\mathcal C\times\mathcal D)^{op}$ and $\mathcal C^{op}\times\mathcal D^{op}$ which for the sake of simplicity we will not distinguish. We will consider every directed set, and in particular every ordered set a category in the usual way. 
\end{notn}
\begin{defn}\label{2.2} In \cite{Go} Goerss constructs a model category structure on $\Gamma-SSets$, called the strict model structure. The corresponding homotopy category will be denoted Ho$(\Gamma-SSets)$ and will be called the homotopy category of simplicial $\Gamma$-sets. Similarly to the construction in chapter 4 of \cite{AM} we may define a Postnikov tower functor
$$(\cdot)^{\natural}:\textrm{Pro}-\Gamma-SSets\longrightarrow
\textrm{Pro}-\Gamma-SSets$$
as follows: if $I$ is a small filtering index category and $\Pi:I^{op}\rightarrow\Gamma-SSets$ is a pro-object of $\Gamma-SSets$ then the functor $\Pi^{\natural}:I^{op}\times\mathbb N^{op}\rightarrow\Gamma-SSets$ is given by the rule:
$$\Pi^{\natural}(\alpha,n)=P_n(\Pi(\alpha))\quad(\forall\alpha\in\textrm{ob}(I),\forall n\in\mathbb N).$$
We will denote by the same symbol the variant of the Postnikov tower functor in the category $\textrm{Pro}-\textrm{Ho}(\Gamma-SSets)$ by the usual abuse of notation.
\end{defn}
\begin{defn}\label{2.3} Next we recall the definition of the relative \'etale homotopy type, following Harpaz and Schlank (see \cite{HS}). Let $K$ be a field, let $\Gamma_K=\textrm{Gal}(\overline K|K)$ denote the absolute Galois group of $K$ and let $Sch_{/K}$ denote the category of locally Noetherian schemes over Spec$(K)$. Let
$$\pi_{0/K}:Sch_{/K}\longrightarrow\Gamma_K-Sets$$
denote the functor which takes the $K$-scheme $X$ to the $\Gamma_K$-set of connected components of $X_{\overline K}$. By applying this functor level-wise and composing it with the localisation functor $\Gamma_K-SSets\rightarrow \textrm{Ho}(\Gamma_K-SSets)$ we get a functor from the category of \'etale hypercoverings of the $K$-scheme $X$ to the homotopy category of simplicial $\Gamma_K$-sets. This construction furnishes, similarly to what is done in chapter 9 of \cite{AM}, another functor:
$$Et_{/K}:Sch_{/K}\longrightarrow\textrm{Pro}-\textrm{Ho}(\Gamma_K-SSets)$$
which we will call the relative \'etale homotopy type of $X$ over $K$. Note that by functoriality we get a natural map:
$$\iota_{X/K}:X(K)\longrightarrow 
[Et_{/K}(\textrm{Spec}(K)),Et_{/K}(X)]
\longrightarrow
[Et_{/K}(\textrm{Spec}(K))^{\natural},Et_{/K}(X)^{\natural}]$$
where the second map is furnished by applying the Postnikov tower functor. We will call the pro-set $[Et_{/K}(\textrm{Spec}(K))^{\natural},Et_{/K}(X)^{\natural}]$ the homotopy fixed points of $X$ and we will denote it by the symbol $X(hK)$.
\end{defn}
Our next aim is to describe the target of this map more explicitly (and to justify the terminology which we have just introduced).
\begin{defn} Let $Sets$ and $SSets$ denote the category of sets and the category of simplicial sets, respectively. Let $\Gamma$ be as above. The category $\Gamma-SSets$ is equipped with a natural concept of homotopy fixed points (see \cite{Go}); we will denote this functor $\Gamma-SSets\rightarrow SSets$ by $\Pi\mapsto\Pi^{h\Gamma}$. For every small filtering index category $I$ and pro-object $\Pi:I^{op}\rightarrow\textrm{Ho}(\Gamma-SSets)$ we define the $\Gamma$-homotopy fixed points set of $\Pi$, denoted by $\Pi(E\Gamma)$, to be
$$\Pi(E\Gamma)=\lim_{\alpha\in\textrm{ob}(I)}\pi_0(\Pi({\alpha})^{h\Gamma}).$$
We will frequently consider the limit $\Pi(E\Gamma)$ as a topological space via its natural pro-discrete topology. This structure is enough to reconstruct the underlying pro-set. By a formula of Harpaz and Schlank there is a natural identification:
$$[(E\Gamma)^{\natural},\Pi^{\natural}]=\Pi(E\Gamma)$$
where $E\Gamma$ is an analogue of the total space of the universal $\Gamma$-bundle in this setting (see Definition 2.3 of \cite{HS}). They also show that when $\Gamma=\Gamma_K$ is the absolute Galois group of a field $K$ then $Et_{/K}(\textrm{Spec}(K))=E\Gamma_K$ hence we have $X(hK)=Et_{/K}(X)(E\Gamma_K)$ for every $K$-scheme $X$, justifying our terminology.
\end{defn}
\begin{notn} Let Ho$(SSets)$ denote the homotopy category of simplicial sets, let $Sch$ denote the category of locally Noetherian schemes, and let 
$$Et:Sch\longrightarrow\textrm{Pro}-\textrm{Ho}(SSets)$$
denote the Artin-Mazur \'etale homotopy type functor. For every $X$ as above and for every $n\in\mathbb N$ let $Et^n(X)$ denote the $n$-th Postnikov piece $P_n(Et(X))$ and let $Et(X)^{\natural}$ denote the Postnikov tower of $Et(X)$, respectively. For every field $K$ and for every scheme $X$ over $K$ let $\overline X$ denote the base change $X_{\overline K}$. Moreover for every such $K$, for every $X\in\textrm{ob}(Sch_{/K})$ and for every $n\in\mathbb N$ let $Et^n_{/K}(X)$ denote the $n$-th Postnikov piece $P_n(Et_{/K}(X))$.
\end{notn} 
\begin{lemma}\label{truncationcomparison} Let $K$ be a field and $X$ be a variety over $K$. Then there are natural isomorphisms: 
$$f^n(X):Et^n(\overline X)\longrightarrow Et^n_{/K}(X)\quad\textrm{and}
\quad f^{\natural}(X):Et(\overline X)^{\natural}
\longrightarrow Et_{/K}(X)^{\natural}$$
in the category $\textrm{\rm Pro}-\textrm{\rm Ho}(SSets)$.
\end{lemma}
\begin{proof} The first half of the claim is Proposition 2.14 of \cite{HS}. The second half is an immediate consequence of the first half and the compatibility of the maps $f^n(X)$.
\end{proof}
\begin{notn} For every $X\in\textrm{ob}(\textrm{\rm Ho}(SSets))$ let $X^{\wedge}\in\textrm{ob}(\textrm{Pro}-\textrm{\rm Ho}(SSets))$ denote its profinite completion. The basic result about the homotopy type of complex algebraic varieties is the following classical theorem of Artin and Mazur:
\end{notn}
\begin{thm}\label{complex} Let $X$ be a geometrically unibranch algebraic variety defined over $\mathbb C$. Then there is a canonical weak homotopy equivalence:
$$\eta_X:Et(X)^{\natural}\longrightarrow(X(\mathbb C)^{\wedge})^{\natural}$$
in $\textrm{\rm Pro}-\textrm{\rm Ho}(SSets)$.
\end{thm}
\begin{proof} This is Corollary 12.10 of \cite{AM} on page 143.
\end{proof}
\begin{prop}\label{connectedcomponents} Assume that $K$ is algebraically closed. Then $\iota_{X/K}$ is surjective. Two points $x,y\in X(K)$ are $H$-equivalent if and only if they lie in the same connected component of $X$ with respect to the Zariski topology. 
\end{prop}
\begin{proof} When $K$ is algebraically closed its absolute Galois group is trivial. Therefore Ho$(\Gamma_K-SSets)$ is the homotopy category of simplicial sets and $Et_{/K}(X)$ is just the usual Artin-Mazur \'etale homotopy type of $X$. Because $K$ is algebraically closed $Et_{/K}(\textrm{Spec}(K))$ is contractible. Therefore there is a natural bijection:
$$X(hK)\cong[*,Et_{/K}(X)^{\natural}]\cong\pi_0(Et_{/K}(X)^{\natural})\cong
\pi_0(X)$$
where the third identification is the consequence of a fundamental comparison theorem of Artin--Mazur (Corollary 10.8 of \cite{AM} on page 122). The claim now follows from the naturality of $\iota_{X/K}$.
\end{proof}

\section{The pointed relative \'etale homotopy type}

\begin{defn} Let $\Gamma$ be as above. By a pointed simplicial $\Gamma$-set we mean a simplicial $\Gamma$-set $S_*$ with a point $p\in S_0$ fixed by $\Gamma$. These form a category $\Gamma-SSets_*$ in the usual way. Note that for every pointed simplicial $\Gamma$-set $S$ and for every $n\in\mathbb N$ the $n$-skeleton $\textrm{sk}_n(S)$, the $n$-coskeleton $\textrm{cosk}_n(S)$, the Kan replacement $\textrm{Ex}^{\infty}(S)$ and $n$-th Postnikov piece $P_n(S)$ are all naturally equipped with a point fixed by $\Gamma$ and we will denote the corresponding four functors, the $n$-skeleton $\textrm{sk}_n:\Gamma-SSets_*\rightarrow\Gamma-SSets_*$, the $n$-coskeleton $\textrm{cosk}_n:\Gamma-SSets_*\rightarrow
\Gamma-SSets_*$, the Kan replacement $\textrm{Ex}^{\infty}:\Gamma-SSets_*\rightarrow \Gamma-SSets_*$ and the $n$-th Postnikov piece $P_n:\Gamma-SSets_*\rightarrow\Gamma-SSets_*$ by the same symbols by a slight abuse of notation.
\end{defn}
\begin{defn} The homotopy category of $\Gamma-SSets_*$ with respect to the pointed version of weak equivalences of Goerss's strict model structure, called the homotopy category of pointed simplicial $\Gamma$-sets, will be denoted by Ho$(\Gamma-SSets_*)$. Similarly to the construction recalled in Definition \ref{2.2} we may define a Postnikov tower functor
$$(\cdot)^{\natural}:\textrm{Pro}-\Gamma-SSets_*\longrightarrow
\textrm{Pro}-\Gamma-SSets_*$$
and we will denote by the same symbol the corresponding the Postnikov tower functor in the category $\textrm{Pro}-\textrm{Ho}(\Gamma-SSets_*)$ by the usual abuse of notation. This is of course justified as the formation of these invariants commute with the forgetful functor $\textrm{Pro}-\textrm{Ho}(\Gamma-SSets_*)\rightarrow\textrm{Pro}-\textrm{Ho}(\Gamma-SSets)$.
\end{defn}
\begin{defn} Since $\Gamma-SSets$ is a model category, it has all colimits, in particular pushouts and equalizers. If $Y_*$ is a simplicial $\Gamma$-set and $X_*\subseteq Y_*$ is a sub simplicial $\Gamma$-set, then let $Y_*/X_*$ denote the simplicial $\Gamma$-set which is the pushout of the inclusion map $X_*\to Y_*$. We call $Y_*/X_*$ the contraction of $Y_*$ by $X_*$. If $X_*,Y_*$ are simplicial $\Gamma$-sets and $f_*:X_*\to Y_*, g_*:X_*\to Y_*$ are maps of simplicial $\Gamma$-sets, then let $Y_*(f_*=g_*)$ denote the simplicial $\Gamma$-set which is the coequalizer of $f_*$ and $g_*$. 
\end{defn}
\begin{defn} Note that for every pair of simplicial $\Gamma$-sets $X_*$ and $Y_*$ the product $X_*\times Y_*$ equipped with the natural (diagonal) $\Gamma$-action is also a simplicial $\Gamma$-set. Similarly the coproduct (disjoint union) $X_*\amalg Y_*$ of simplicial $\Gamma$-sets $X_*$ and $Y_*$ is also a simplicial $\Gamma$-set equipped with the tautological $\Gamma$-action. Let $I$ denote the $1$-simplex $\Delta^1=[0,1]$ with the trivial $\Gamma$-action; this choice makes it into a simplicial $\Gamma$-set. For every morphism $f:X_*\to Y_*$ of simplicial $\Gamma$-sets the mapping cylinder $Cyl(f)$ is the coequalizer of the two maps:
$$f':X\longrightarrow Y_*\amalg X_*\times I\textrm{ and }
p:X\longrightarrow Y_*\amalg X_*\times I,$$
where $f'$ is the composition of $f$ and the tautological inclusion $Y_*\subseteq Y_*\amalg X_*\times I$, and $p$ is the composition of the map identifying $X_*$ with $X_*\times\{1\}\subseteq X_*\times I$ and the tautological inclusion $X_*\times I\subset Y_*\amalg X_*\times I$. We define the mapping cone $Cone(f)$ of an $f:X_*\to Y_*$ as above as the contraction of $Cyl(f)$ by the image of the map:
$$q:X_*\longrightarrow Cyl(f),$$
where $q$ is the composition of the map identifying $X_*$ with $X_*\times\{0\}\subseteq X_*\times I$ with the tautological inclusion $X_*\times I\subset Y_*\amalg X_*\times I$ composed with the natural surjection $Y_*\amalg X_*\times I\mapsto Cyl(f)$.  Note that $Cone(f)$ is canonically a pointed simplicial $\Gamma$-set where the base point is the image of $q(X_0)\subset Cyl(f)_0$ under the contraction map $Cyl(f)\to Cone(f)$.
\end{defn}
\begin{defn}\label{3.6} By a pointed $K$-scheme $(X,x)$ we will mean a locally Noetherian scheme $X$ over $K$ with a $K$-valued point $x:\textrm{Spec}(K)\to X$ on $X$. These form the objects of a category $Sch_{/K*}$ where a morphism $f$ from an object $(X,x)$ to another object $(Y,y)$ is a map $f:X\to Y$ of schemes over $K$ such that $f(x)=y$. Now let $(X,x)$ be a pointed $K$-scheme and let $H_*$ be an \'etale hypercovering of $X$. Then the pull-back $x^*(H_*)$ is an \'etale hypercovering of $\textrm{Spec}(K)$, and the map $x$ induces a morphism $x_*(H_*):\pi_{0/K}(x^*(H_*))\to\pi_{0/K}(H_*)$ of simplicial $\Gamma$-sets. Let $\pi_{0/K}(H_*,x)$ denote the mapping cone of the composition of this map $x_*(H_*)$ and the canonical inclusion $\pi_{0/K}(H_*)\subseteq\textrm{Ex}^{\infty}(\pi_{0/K}(H_*))$; it is a pointed simplicial $\Gamma$-set. A map $f:H_*\to J_*$ between \'etale hypercoverings of $X$ induces a map $\pi_{0/K}(f,x):\pi_{0/K}(H_*,x)\to\pi_{0/K}(J_*,x)$ between pointed simplicial $\Gamma$-sets, and a homotopy between two maps $f:H_*\to J_*,g:H_*\to J_*$ induces a pointed $\Gamma$-equivariant homotopy between $\pi_{0/K}(f,x)$ and $\pi_{0/K}(g,x)$. Therefore we may apply part $(i)$ of Corollary 8.13 of \cite{AM} on page 105 to conclude that the functor
$$H_*\mapsto\pi_{0/K}(H_*,x)$$
above induces an object $Et_{/K}(X,x)$ of $\textrm{Pro}-\textrm{Ho}(\Gamma-SSets_*)$. We will call the latter the pointed relative \'etale homotopy type of $(X,x)$.
\end{defn}
\begin{notn}\label{3.6b} Let $(X,x)$ be a pointed $K$-scheme. For every \'etale hypercovering $H_*$ of $X$ let 
$$i(H_*,x):\pi_{0/K}(H_*)\longrightarrow\pi_{0/K}(H_*,x)$$
be the composition of the natural inclusion map:
$$\pi_{0/K}(H_*)\longrightarrow Cone(x_*(H_*))$$
and the map $Cone(x_*(H_*))\to\pi_{0/K}(H_*,x)$ induced by the functoriality of mapping cones. This is a natural transformation between two functors from the category of \'etale hypercoverings over $X$ into $\Gamma-SSets$, and hence it induces a map:
$$i(X,x):Et_{/K}(X,x)\longrightarrow Et_{/K}(X)$$
of pro-objects of the homotopy category $\textrm{Ho}(\Gamma-SSets)$, where by slight abuse of notation we let $Et_{/K}(X,x)$ also denote the image of the pointed relative \'etale homotopy type of $(X,x)$ with respect to the forgetful functor
$$\textrm{Pro}-\textrm{Ho}(\Gamma-SSets_*)\longrightarrow
\textrm{Pro}-\textrm{Ho}(\Gamma-SSets).$$
The map $i(X,x)$ is obviously natural.
\end{notn}
\begin{prop}\label{comp1} The map $i(X,x):Et_{/K}(X,x)\to Et_{/K}(X)$ induces a bijection on homotopy fixed points.
\end{prop}
\begin{proof} Note that the functor $\pi_{0/K}$ induces an equivalence between the category of \'etale coverings over Spec$(K)$ and $\Gamma-Sets$, so $\pi_{0/K}(x^*(H_*))$ is a contractible simplicial $\Gamma$-set. Therefore $i(H_*,x)$ is a weak equivalence for every \'etale hypercovering $H_*$ of $X$ with respect to Goerss's {\it weak} model structure (see Theorem A on page 189 and Definition 1.11 on page 194 in \cite{Go}), and hence induces a bijection $\pi_0(\pi_{0/K}(H_*)^{h\Gamma})\to\pi_0(Cone(x_*(H_*))^{h\Gamma})$. So the same holds for $i(X,x)$, too.
\end{proof}
\begin{defn} Let $Sch_*$ be the category of pointed locally Noetherian schemes, and as usual, we denote the objects of $Sch_*$ as pairs $(X,x)$, where $X$ is a locally Noetherian scheme and $x$ is a geometric point of $X$. By slight abuse of notation let 
$$Et:Sch_*\longrightarrow\textrm{Pro}-\textrm{Ho}(SSets_*)$$
denote the pointed version of the Artin-Mazur \'etale homotopy type functor. For every object $(X,x)$ of $Sch_*$ let $\pi_n(X,x)$ denote the $n$-th homotopy group of $Et(X,x)$ for every $n\geq1$ when $X$ is connected. 
\end{defn}
\begin{notn}\label{3.9} For every pointed scheme $(X,x)$ and for every $n\in\mathbb N$ let $Et^n(X,x)$ denote the $n$-th Postnikov piece $P_n(Et(X,x))$ and let $Et(X,x)^{\natural}$ denote the Postnikov tower of $Et(X,x)$, respectively. Similarly for every field $K$ and for every pointed $K$-scheme $(X,x)$ and for every $n\in\mathbb N$ let $Et^n_{/K}(X,x)$ denote the $n$-th Postnikov piece $P_n(Et_{/K}(X,x))$ and let $Et_{/K}(X,x)^{\natural}$ denote the Postnikov tower of $Et_{/K}(X,x)$, respectively. Since we fixed a separable closure $\overline K$ of $K$, we may associate to
every $K$-valued point $x:\textrm{Spec}(K)\to X$ of a $K$-scheme $X$ a
$\overline K$-valued point $\overline x:\textrm{Spec}(K)\to X$ which is the composition of the map $\textrm{Spec}(\overline K)\to\textrm{Spec}(K)$ induced by the inclusion $K\subseteq\overline K$ and $x$.
\end{notn}
Similarly to above, by slight abuse of notation we let $Et_{/K}(X,x),Et^n_{/K}(X,x)$ also denote the image of the pointed relative \'etale homotopy type of $(X,x)$ and of the $n$-th truncation of the latter with respect to the forgetful functor
$$\textrm{Pro}-\textrm{Ho}(\Gamma-SSets_*)\longrightarrow
\textrm{Pro}-\textrm{Ho}(SSets_*).$$
\begin{prop}\label{truncationcomparison2} Let $K$ be a field and $(X,x)$ be a pointed $K$-scheme such that $X$ is a variety over $K$. Then there are natural isomorphisms: 
$$f^n(X,x):Et^n(\overline X,\overline x)\longrightarrow Et^n_{/K}(X,x)\quad\textrm{and}
\quad f^{\natural}(X,x):Et(\overline X,\overline x)^{\natural}
\longrightarrow Et_{/K}(X,x)^{\natural}$$
in the category $\textrm{\rm Pro}-\textrm{\rm Ho}(SSets_*)$.
\end{prop}
\begin{proof} Since the second half is an immediate consequence of the first half and the compatibility of the maps $f^n(X,x)$, it will be enough to prove the former. Let $\Delta^0$ denote the $0$-simplex, as usual. Let $(H_*,h)$ be a  pointed \'etale hypercovering of $(X,\overline x)$. By definition $h$ is a map $\Delta^0\to\pi_{0/K}(x^*(H_*))$ of simplicial sets. For every $(H_*,h)$ as above let $\pi_{0/K}(H_*,h,x)$ be the contraction of $\pi_{0/K}(H_*,x)$ by the image of the map
$$c_h:\Delta^1\cong\Delta^0\times\Delta^1\longrightarrow\pi_{0/K}(H_*,x)$$
of simplicial sets, where $c_h$ is the composition $\pi\circ\iota\circ(h\times\textrm{id}_I)$,
where $\iota$ is the inclusion:
$$\pi_{0/K}(x^*(H_*))\times\Delta^1\subset
\textrm{Ex}^{\infty}(\pi_{0/K}(H_*))\amalg\pi_{0/K}(x^*(H_*))\times\Delta^1,$$
and $\pi$ is the canonical surjection
$$\textrm{Ex}^{\infty}(\pi_{0/K}(H_*))\amalg\pi_{0/K}(x^*(H_*))\times\Delta^1\longrightarrow\pi_{0/K}(H_*,x).$$
Let
$$a(H_*,h):\pi_{0/K}(H_*,x)\longrightarrow\pi_{0/K}(H_*,h,x)$$
be the contraction map. We will consider $\pi_{0/K}(H_*,h,x)$ a pointed simplicial set, where its distinguished point $b(h)\in\pi_{0/K}(H_*,h,x)_0$ is the image of the base point of the pointed simplicial $\Gamma$-set under $a(H_*,h)$. Note that this map is a weak equivalence in $\textrm{\rm Pro}-\textrm{Ho}(SSets_*)$ since we contracted a contractible sub simplicial set. Therefore by $(i)$ of Corollary 8.13 of \cite{AM} on page 105 the functor
$$(H_*,h)\mapsto(\pi_{0/K}(H_*,h,x),b(h))$$
induces an object $Et'(\overline X,\overline x)$ of $\textrm{Pro}-\textrm{Ho}(SSets_*)$ which is isomorphic to $Et_{/K}(X,x)$ in this category. For every $(H_*,h)$ as above let $\pi_{0/K}(h)\in\pi_{0/K}(H_0)$ be the point corresponding to $h$ and let
$$b(H_*,h):\pi_{0/K}(H_*)\longrightarrow\pi_{0/K}(H_*,x)$$
be the map of pointed simplicial sets which is the composition of the natural inclusion map:
$$\pi_{0/K}(H_*)\subset\pi_{0/K}(H_*)\amalg\pi_{0/K}(x^*(H_*))
\times\Delta^1$$
with the natural surjection:
$$\pi_{0/K}(H_*)\amalg\pi_{0/K}(x^*(H_*))\times\Delta^1
\longrightarrow Cone(x_*)$$
composed with the map $Cone(x_*)\to\pi_{0/K}(H_*,x)$ induced by the functoriality of mapping cones. Let $c(H_*,h)$ be the composition of $b(H_*,h)$ with $a(H_*,h)$; then this map is a morphism:
$$c(H_*,h):(\pi_{0/K}(H_*),\pi_{0/K}(h))\longrightarrow
(\pi_{0/K}(H_*,h,x),b(h))$$
of pointed simplicial sets. Since $\pi_{0/K}(x^*(H_*))$ is a contractible simplicial $\Gamma$-set, them map  $c(H_*,h)$ is a weak equivalence for every pointed \'etale hypercovering $(H_*,h)$ of $X$. Therefore by $(i)$ of Corollary 8.13 of \cite{AM} on page 105 the functor
$$(H_*,h)\mapsto(\pi_{0/K}(H_*),\pi_{0/K}(h))$$
induces an object $Et''(\overline X,\overline x)$ of $\textrm{Pro}-\textrm{Ho}(SSets_*)$ which is isomorphic to $Et'(\overline X,\overline x)$, and hence to $Et_{/K}(X,x)$ in this category. Therefore it will be sufficient to prove that there are natural isomorphisms: 
$$g^n(X,x):Et_n(\overline X,\overline x)\longrightarrow P_n(Et''(\overline X,
\overline x))$$
in the category $\textrm{\rm Pro}-\textrm{\rm Ho}(SSets_*)$ which are compatible with each other and with truncation.

Note that the inclusion of the indexing category of $Et''(\overline X,\overline x)$ in the indexing category of $Et(\overline X,\overline x)$ furnishes a natural map
$$g(X,x):Et(\overline X,\overline x)\longrightarrow Et''(\overline X,
\overline x).$$
In order to prove that $g(X,x)$ is an isomorphism after taking $n$-th truncations we can argue the same way as Harpaz and Schlank did in the proof of Proposition 2.14 of \cite{HS}.
\end{proof}

\section{Homotopy groups of finite \'etale covers}

Recall that for every object $X$ of $Sch$ with a geometric point $x$, and for every $n\in\mathbb N$ let $\pi_n(X,x)$ denote the homotopy group $\pi_n(Et(X),x)$. In this section assume that $K$ is an algebraically closed field.
\begin{prop}\label{imsmarter} Let $f:(X,x)\rightarrow (Y,y)$ be a finite \'etale map of pointed smooth, connected quasi-projective varieties over $K$. Then the induced map:
$$\pi_n(f):\pi_n(X,x)\longrightarrow\pi_n(Y,y)$$
is an isomorphism for every $n\geq2$.
\end{prop}
\begin{proof} Because $\pi_1(Y,y)$ is topologically finitely generated, its open normal subgroups are cofinal, and hence there is a finite \'etale map $g:(Z,z)\rightarrow (X,x)$ of pointed smooth, connected quasi-projective varieties over $K$ such that the image of $\pi_1(f\circ g):\pi_1(Z,z)\rightarrow\pi_1(Y,y)$ is an open normal subgroup. In this case the image of $\pi_1(f):\pi_1(Z,z)\rightarrow\pi_1(X,x)$ is an open normal subgroup, too. It will be enough to show that the maps $\pi_n(f\circ g)$ and $\pi_n(f)$ are isomorphisms for every $n\geq2$. Since the composition $f\circ g$ is also a finite \'etale map, we've reduced the claim to the special case when the image of  $\pi_1(f)$ is an open normal subgroup. 

In this case $f$ is a finite Galois covering; let $G$ denote the covering group and let $\alpha:G\rightarrow\textrm{Aut}(X)$ be the action corresponding to the deck  transformations. Let $(H_*,h)$ be a pointed hypercovering of $Y$ with respect to the \'etale site of $Y$ pointed with $y$. Let $C_*$ denote the \'etale \v Cech hypercovering:
$$\xymatrix{
 X\ar@/^.5pc/@<1ex>[r] & X\times_YX \ar@<-.5ex>[l] \ar@<.5ex>[l]
 & X\times_YX\times_YX \ar@<-.7ex>[l] \ar[l] \ar@<.7ex>[l]
 & \cdots \ar@<-1ex>[l] \ar@<-.3ex>[l] \ar@<.3ex>[l] \ar@<1ex>[l]}$$
generated by the cover $X\rightarrow Y$ and equip it with a point $c$ with respect to the same pointed site. Let $(I_*,i)$ be the fibre product of $(H_*,h)$ and $(C_*,c)$ over $Y$; this pointed simplicial object is also a hypercovering. Let $(J_*,j)$ be the pull-back $f^*(I_*,i)$ of the pointed hypercovering $(I_*,i)$ onto the \'etale site of $X$ pointed with $x$ with respect to $f$ and let $f_*:(\pi_0(J_*),\pi_0(j))\rightarrow(\pi_0(I_*),\pi_0(i))$ be the map induced by $f$ between the pointed simplicial sets of connected components of $J_*$ and $I_*$.

The action $\alpha$ induces an action of $G$ on  $J_*$, and hence an action of $G$ on $\pi_0(J_*)$. If we equip $\pi_0(I_*)$ with the trivial action then $f_*$ is $G$-equivariant. Because $f$ is a finite \'etale cover the map $f_*$ is surjective. Moreover for every $n$ the \'etale cover $I_n\rightarrow Y$ factors through $f:X\rightarrow Y$, and hence the action of $G$ on the connected components of the base change  $J_n\rightarrow X$ of this map to $X$ is free. We get that $f_*:(\pi_0(J_*),\pi_0(j))\rightarrow(\pi_0(I_*),\pi_0(i))$ is a $G$-cover and hence the induced maps
$$\pi_n(f_*):\pi_n(\pi_0(J_*),\pi_0(j))\longrightarrow\pi_n(\pi_0(I_*),\pi_0(i))$$
are isomorphisms for all $n\geq2$.

Let $(L_*,l)$ be a pointed hypercovering of $X$ with respect to the \'etale site of $X$ pointed with $x$. By composing the structure maps with $f$ we get a pointed hypercovering of $Y$ with respect to the \'etale site of $Y$ pointed with $y$ which we will denote by $(H_*,h)$ by slight abuse of notation. By applying the same construction to $(H_*,h)$ as above we get a pointed hypercovering $(J_*,j)$ of the \'etale site of $X$ pointed with $x$ which dominates $(L_*,l)$. Therefore pointed hypercovers of $X$ of the form as $(J_*,j)$ above are cofinal, so the injectivity of the maps $\pi_n(f)$ for all $n\geq2$ follows.

Let $\gamma$ be an element of $\pi_n(Y,y)$ (where $n\geq2$). For every pointed hypercovering $(L_*,l)$ of $X$ with respect to the \'etale site of $X$ pointed with $x$ we will construct an element $\gamma'_{(L_*,l)}\in\pi_n(\pi_0(L_*),\pi_0(l))$ as follows. Let $(H_*,h)$, $(I_*,i)$ and $(J_*,j)$ be the same as in the paragraph above. Let $\gamma_{(I_*,i)}\in\pi_n(\pi_0(I_*),\pi_0(i))$ be the image of $\gamma$ under the tautological map $\pi_n(Y,y)\to\pi_n(\pi_0(I_*),\pi_0(i))$, let $\pi_n(f_*)^{-1}(\gamma_{(I_*,i)})$ be unique pre-image of $\gamma_{(I_*,i)}$ under the isomorphism:
$$\pi_n(f_*):\pi_n(\pi_0(J_*),\pi_0(j))\longrightarrow\pi_n(\pi_0(I_*),\pi_0(i)),$$
and let $\gamma'_{(L_*,l)}$ be the image of $\pi_n(f_*)^{-1}(\gamma_{(I_*,i)})$ under the natural map
$$\pi_n(\pi_0(J_*),\pi_0(j))\longrightarrow\pi_n(\pi_0(L_*),\pi_0(l)).$$
It is easy to check that the elements $\gamma\_{(L_*,l)}$ glue together to an element $\gamma'$ of $\pi_n(X,x)$ whose image under $\pi_n(f)$ is $\gamma$. The surjectivity of the maps $\pi_n(f)$ for all $n\geq2$ follows. 
\end{proof}
\begin{notn} For every group $\Gamma$ let $\widehat{\Gamma}$ be its profinite completion. For every object $X$ of $Sch$, for every $n\in\mathbb N$ and for every pro-abelian group $A$ let $H_n(X,A)$ denote the homology group $H_n(Et(X),A)$. For every object $(X,x)$ of $Sch_*$ and $n$ as above let:
$$h_n(X,x):\pi_n(X,x)\longrightarrow H_n(X,\widehat{\mathbb Z})$$
denote the Hurewitz map. Let $X$ be a smooth, geometrically irreducible, quasi-projective variety over $K$. Let $x$ be a $K$-valued point of $X$ and let $\textrm{\rm Fet}(X,x)$ denote the category of finite \'etale pointed connected covers $(Y,y)$ of $(X,x)$ such that the image of the induced map $\pi_1(f):\pi_1(Y,y)\rightarrow\pi_1(Y,y)$ is an open characteristic subgroup for every object $f:(Y,y)\rightarrow (X,x)$. Since for every $f:(Y,y)\rightarrow (X,x)$ as above the induced map $\pi_2(f):\pi_2(Y,y)\rightarrow\pi_2(X,x)$ is an isomorphism by Proposition \ref{imsmarter}, the projective limit of the inverses of these maps is an isomorphism:
$$a_X:\pi_2(X,x)\longrightarrow\!\!\!\!\!\!\!\lim_{(Y,y)\in\textrm{\rm Fet}(X,x)}\pi_2(Y,y).$$
Moreover we may take the projective limit of the Hurewitz maps:
$$b_X:\!\!\!\!\!\!\!\lim_{(Y,y)\in\textrm{\rm Fet}(X,x)}\!\!\!\!\!\!\!\pi_2(Y,y)
\longrightarrow
\!\!\!\!\!\!\!\lim_{(Y,y)\in\textrm{\rm Fet}(X,x)}\!\!\!\!\!\!\!H_2(Y,\widehat{\mathbb Z}).$$
\end{notn}
\begin{thm}\label{myhurewitz} The map:
$$b_X\circ a_X:\pi_2(X,x)\longrightarrow\lim_{(Y,y)\in\textrm{\rm Fet}(X,x)}H_2(Y,\widehat{\mathbb Z})$$
is an isomorphism.
\end{thm}
\begin{proof} First we are going to prove that $b_X$ is injective. Let
$\gamma$ be a non-zero element of $\pi_2(X,x)$. Then there is a pointed hypercovering $(H_*,h)$ of $X$ with respect to the \'etale site of $X$ pointed with $x$ such that the image of $\gamma$ under the natural map $\pi_2(X,x)
\rightarrow\pi_2(\pi_0(H_*),\pi_0(h))$ is non-zero. Let $N$ be the kernel of the natural map $\pi_1(X,x)\rightarrow\pi_1(\pi_0(H_*),\pi_0(h))$; it is an open normal subgroup. Because $\pi_1(X,x)$ is topologically finitely generated, its open characteristic subgroups are cofinal, and hence there is an object $f:(Y,y)\rightarrow (X,x)$ of $\textrm{\rm Fet}(X,x)$ such that the image of $\pi_1(f):\pi_1(Y,y)\rightarrow\pi_1(X,x)$ lies in $N$.  

Note that $f$ is a finite Galois covering; let $G$ denote the covering group and let $\alpha:G\rightarrow\textrm{Aut}(X)$ be the action corresponding to the deck  transformations. Let $C_*$ denote the \'etale \v Cech hypercovering generated by the cover $Y\rightarrow X$ and equip it with a point $c$ with respect to the same pointed site. Let $(I_*,i)$ be the pointed simplicial object which is the fibre product of $(H_*,h)$ and $(C_*,c)$ and let $(J_*,j)$ be the pull-back $f^*(I_*,i)$ of the pointed hypercovering $(I_*,i)$ onto the \'etale site of $X$ pointed with $x$ with respect to $f$. Let $f_*:(\pi_0(J_*),\pi_0(j))\rightarrow(\pi_0(I_*),\pi_0(i))$ be the map induced by $f$ between the pointed simplicial sets of connected components of $J_*$ and $I_*$.

We may argue as above to conclude that $f_*:(\pi_0(J_*),\pi_0(j))\rightarrow(\pi_0(I_*),\pi_0(i))$ is a $G$-cover with respect to the action induced by $\alpha$ on  $\pi_0(J_*)$ and the trivial action on $\pi_0(I_*)$. Therefore the induced map
$$\pi_1(f_*):\pi_1(\pi_0(J_*),\pi_0(j))\longrightarrow\pi_1(\pi_0(I_*),\pi_0(i))$$
is injective. There is a commutative diagram:
$$
\CD
\pi_1(Y,y)@>>>\pi_1(\pi_0(J_*),\pi_0(j))\\
@V\pi_1(f)VV@V\pi_1(f_*)VV\\
\pi_1(X,x)
@>>>\pi_1(\pi_0(I_*),\pi_0(i)).
\endCD
$$
By assumption the composition of $\pi_1(f)$ and the lower horizontal map has trivial image. Since the upper horizontal map is surjective by Corollary 10.6 of 
\cite{AM} on pages 121--122, we get that $\pi_1(f_*)$ has trivial image, too, and hence $\pi_1(\pi_0(J_*),\pi_0(j))$ is the trivial group. There is a similar commutative diagram:
$$
\CD
\pi_2(Y,y)@>>>\pi_2(\pi_0(J_*),\pi_0(j))\\
@V\pi_2(f)VV@V\pi_2(f_*)VV\\
\pi_2(X,x)
@>>>\pi_2(\pi_0(I_*),\pi_0(i))
\endCD
$$
for $\pi_2$. The image of $\gamma$, considered as an element of $\pi_2(Y,y)$, is non-zero under the composition of $\pi_2(f)$ and the lower horizontal map by assumption. Therefore its image $\gamma'$ under the upper horizontal map is also non-zero. Since $\pi_1(\pi_0(J_*),\pi_0(j))$ is trivial we get that the image of $\gamma'$ under the Hurewitz map $\pi_2(\pi_0(J_*),\pi_0(j))\rightarrow H_2(\pi_0(J_*),\widehat{\mathbb Z})$ is non-zero. By naturality this implies that the image of $\gamma$ under the Hurewitz map $h_2(Y,y):\pi_2(Y,y)\longrightarrow H_2(Y,\widehat{\mathbb Z})$ is non-zero, too. 

Next we are going to prove that $b_X$ is surjective. For every $(Z,z)\in\textrm{\rm Fet}(X,x)$ let $\mathcal C(Z,z)$ be the category of morphisms
$(Y,y)\rightarrow(Z,z)$ in $\textrm{\rm Fet}(X,x)$ and let
$$c_{(Z,z)}:
\!\!\!\!\!\!\lim_{(Y,y)\rightarrow(Z,z)\in\mathcal C(Z,z)}\!\!\!\!\!\!
H_2(Y,\widehat{\mathbb Z})\rightarrow H_2(Z,\widehat{\mathbb Z})$$
be the tautological map. Because the pre-image $(b_X\circ a_X)^{-1}(\textrm{Im}(c_{(Z,z)}))\subseteq\pi_2(X,x)$ is closed and $\pi_2(X,x)$ is profinite, by compactness it will be enough to show that for every object $(Z,z)$ as above $(b_X\circ a_X)^{-1}(\textrm{Im}(c_{(Z,z)}))$ is non-empty. Now fix an element $\gamma\in\textrm{Im}(c_{(Z,z)})$ and choose a
$$\gamma'\in\!\!\!\!\!\!\lim_{(Y,y)\rightarrow(Z,z)\in\mathcal C(Z,z)}\!\!\!\!\!\!
H_2(Y,\widehat{\mathbb Z})$$
such that $c_{(Z,z)}(\gamma')=\gamma$. For every pointed hypercovering $(H_*,h)$ of $Z$ with respect to the \'etale site of $Z$ pointed with $z$ we are going to construct an element
$\gamma_{(H_*,h)}\in\pi_2(\pi_0(H_*),\pi_0(h))$ as follows. 

Let $N$ be the kernel of the natural map $\pi_1(Z,z)\rightarrow\pi_1(\pi_0(H_*),\pi_0(h))$. Using that $\pi_1(Z,z)$ is topologically finitely generated as above we get that there is a morphism $f:(Y,y)\rightarrow (Z,z)$ of $\textrm{\rm Fet}(X,x)$ such that the image of $\pi_1(f):\pi_1(Y,y)\rightarrow\pi_1(Z,z)$ lies in $N$. Let $C_*$ denote the \'etale \v Cech hypercovering generated by the cover $Y\rightarrow Z$ and equip it with a point $c$ with respect to the same pointed site. Let $(I_*,i)$ be the pointed hypercovering which is the fibre product of $(H_*,h)$ and $(C_*,c)$ and let $(J_*,j)$ be the pull-back $f^*(I_*,i)$ of $(I_*,i)$ onto the \'etale site of $Y$ pointed with $y$ with respect to $f$. Let $f_*:(\pi_0(J_*),\pi_0(j))\rightarrow(\pi_0(I_*),\pi_0(i))$ be the map induced by $f$ between the pointed simplicial sets of connected components of $J_*$ and $I_*$.

As we saw in the proof of injectivity the group $\pi_1(\pi_0(J_*),\pi_0(j))$ is trivial and hence the Hurewitz map $\pi_2(\pi_0(J_*),\pi_0(j))\rightarrow H_2(\pi_0(J_*),\widehat{\mathbb Z})$ is an isomorphism. Therefore there is a unique $\gamma_f\in\pi_2(\pi_0(J_*),\pi_0(j))$ whose image under this Hurewitz map is the image of $\gamma'$ under the composition of $c_{(Y,y)}$ and the natural map $H_2(Y,\widehat{\mathbb Z})\rightarrow H_2(\pi_0(J_*),\widehat{\mathbb Z})$. Let $\gamma_{(H_*,h)}\in\pi_2(\pi_0(H_*),\pi_0(h))$ be the image of $\gamma_f$ under the composition of the functorial maps
$\pi_2(f_*):\pi_2(\pi_0(J_*),\pi_0(j))\rightarrow\pi_2(\pi_0(I_*),\pi_0(i))$ and 
$\pi_2(\pi_0(I_*),\pi_0(i))\rightarrow\pi_2(\pi_0(H_*),\pi_0(h))$.

First we are going to show that $\gamma_{(H_*,h)}$ is independent of the choice of the morphism $f$. Let $f':(Y',y')\rightarrow (Z,z)$ be another morphism of $\textrm{\rm Fet}(X,x)$ such that the image of $\pi_1(f'):\pi_1(Y',y')\rightarrow\pi_1(Z,z)$ lies in $N$. Then the fibre product of $f$ and $f'$ over $(Z,z)$ is another morphism of $\textrm{\rm Fet}(X,x)$ with this property. Moreover it can be factorized into a composition of a morphism $g$ of $\textrm{\rm Fet}(X,x)$ and $f$, and into a composition of a morphism $g'$ of $\textrm{\rm Fet}(X,x)$ and $f'$, too. Therefore it will be enough to show that this construction applied to the morphism $f\circ g$ of $\textrm{\rm Fet}(X,x)$ will give the same element in $\pi_2(\pi_0(H_*),\pi_0(h))$ as $f$, where $g:(V,v)\rightarrow (Y,y)$ is a morphism of $\textrm{\rm Fet}(X,x)$.

Let $(C_*,c)$ be the same Cech hypercovering as above and let $C'_*$ denote the \'etale \v Cech hypercovering generated by the cover $V\rightarrow Z$ and equip it with a point $c'$ with respect to the same pointed site. Note that $g$ furnishes a map $\delta:(C'_*,c')\rightarrow(C_*,c)$ of pointed hypercoverings. Let $(C_*,c),(I_*,i)$ and $(J_*,j)$ be as above, let $(I'_*,i')$ be the fibre product of $(H_*,h)$ and $(C'_*,c')$, and let $(J'_*,j')$ be the pull-back $(f\circ g)^*(I'_*,i')$. Let $f_*:(\pi_0(J_*),\pi_0(j))\rightarrow(\pi_0(I_*),\pi_0(i))$ and $(f\circ g)_*:(\pi_0(J'_*),\pi_0(j'))\rightarrow(\pi_0(I'_*),\pi_0(i'))$ be the map induced by $f$ and by $f\circ g$, respectively. Note that $g$ induces a map $g_*:g^*(\pi_0(J_*),\pi_0(j))\rightarrow (\pi_0(J_*),\pi_0(j))$. Then there is a commutative diagram:
$$
\CD
\pi_2(\pi_0(J'_*),\pi_0(j'))@>\pi_2((f\circ g)_*)>>\pi_2(\pi_0(I'_*),\pi_0(i'))
@>>>\pi_2(\pi_0(H_*),\pi_0(h))\\
@VV\pi_2(g_*\circ(f\circ g)^*(\delta_H))V@VV\pi_2(\delta_H)V@|\\
\pi_2(\pi_0(J_*),\pi_0(j))@>\pi_2(f_*)>>\pi_2(\pi_0(I_*),\pi_0(i))
@>>>\pi_2(\pi_0(H_*),\pi_0(h)),
\endCD
$$
where $\delta_H:(I'_*,i')\rightarrow(I_*,i)$ is the fibre product of $\delta$ with $(H,h)$, the morphism $(f\circ g)^*(\delta_H):(\pi_0(J'_*),\pi_0(j'))\rightarrow g^*(\pi_0(J_*),\pi_0(j))=(f\circ g)^*(\pi_0(I_*),\pi_0(i))$ is the base change of $\delta_H$ with respect to $f\circ g$; while the middle and left hand side vertical maps are induced by $\delta_H$ and $g_*\circ(f\circ g)^*(\delta_H)$, respectively. Since the image of $\gamma_{f\circ g}$ under the left vertical map is $\gamma_f$, the claim above follows.

In order to conclude the proof of the theorem itself we only need to show that for every morphism $\delta:(H_*,h)\rightarrow(H'_*,h')$ of pointed hypercoverings of $Z$ with respect to the \'etale site of $Z$ pointed with $z$ the induced map $\pi_2(\delta):\pi_2(\pi_0(H_*),\pi_0(h))\rightarrow\pi_2(\pi_0(H'_*),\pi_0(h'))$ takes $\gamma_{(H_*,h)}$ to $\gamma_{(H'_*,h')}$. Indeed in this case these $\gamma_{(H_*,h)}$ glue together to an element of $\pi_2(Z,z)$ whose image is $\gamma$ under the Hurewitz map, by construction. Now let $f:(Y,y)\rightarrow (Z,z)$ be a morphism of $\textrm{\rm Fet}(X,x)$ such that the image of $\pi_1(f):\pi_1(Y,y)\rightarrow\pi_1(Z,z)$ lies in the kernel of the natural map $\pi_1(Z,z)\rightarrow\pi_1(\pi_0(H_*),\pi_0(h))$. Let $(C_*,c),(I_*,i)$ and $(J_*,j)$ be as above. Note that the image of $\pi_1(f)$ lies in the kernel of the natural map $\pi_1(Z,z)\rightarrow\pi_1(\pi_0(H'_*),\pi_0(h'))$, too. Let $(I'_*,i')$ be the fibre product of $(H'_*,h')$ and $(C_*,c)$ and let $(J'_*,j')$ be the pull-back $f^*(I'_*,i')$ onto the \'etale site of $Y$ pointed with $y$ with respect to $f$. By slight abuse of notation let $f_*$ denote both maps $(\pi_0(J_*),\pi_0(j))\rightarrow(\pi_0(I_*),\pi_0(i))$ and $(\pi_0(J'_*),\pi_0(j'))\rightarrow(\pi_0(I'_*),\pi_0(i'))$ induced by $f$. Then there is a commutative diagram:
$$
\CD
\pi_2(\pi_0(J_*),\pi_0(j))@>\pi_2(f_*)>>\pi_2(\pi_0(I_*),\pi_0(i))@>>>
\pi_2(\pi_0(H_*),\pi_0(h))\\
@V\pi_2(f^*(\delta_C))VV@V\pi_2(\delta_C)VV@V\pi_2(\delta)VV\\
\pi_2(\pi_0(J'_*),\pi_0(j'))@>\pi_2(f_*)>>\pi_2(\pi_0(I'_*),\pi_0(i'))
@>>>\pi_2(\pi_0(H'_*),\pi_0(h')),
\endCD
$$
where the middle and left hand side vertical maps are induced by the fibre product $\delta_C$ of $\delta$ with $(C,c)$, and the pull-back $f^*(\delta_C)$ of $\delta_C$ with respect to $f$, respectively. By construction the image of $\gamma_f\in\pi_2(\pi_0(J_*),\pi_0(j))$ constructed for $(H,h)$ under $\pi_2(f^*(\delta_C))$ is the element denoted by the same symbol and constructed for $(H',h')$. The claim is now clear.
\end{proof}

\section{The homotopy type of curves and abelian varieties}

\begin{defn} Following Serre (see page 16 of \cite{Se}) we will say that a group $\Gamma$ is good if the homomorphism of cohomology groups $H^n(\widehat{\Gamma}, M)\rightarrow H^n(\Gamma, M)$ induced by the natural homomorphism $\Gamma\rightarrow\widehat{\Gamma}$ is an isomorphism for every finite $\Gamma$-module $M$. For every smooth, connected, quasi-projective variety $X$ over any field $K$ let $\pi_n(X)$ denote the isomorphism class of the $n$-th homotopy group of $\pi_n(X,x)$ for some geometric point $x$ and for every $n\geq1$. As the notation indicates these isomorphism classes do not depend on the choice of the base point.
\end{defn}
\begin{prop}\label{good} Let $X$ be a smooth variety over $\mathbb C$ such that $X(\mathbb C)$ has the homotopy type of the Eilenberg-MacLane space $B\pi_1(X(\mathbb C))$ and the group $\pi_1(X(\mathbb C))$ is good. Then $Et(X)$ is weakly homotopy equivalent to $B\pi_1(X)$.
\end{prop}
\begin{proof} By Corollary 6.6 of \cite{AM} on page 72 the profinite completion $(B\Gamma)^{\wedge}$ of $B\Gamma$ is weakly homotopy equivalent to $B\widehat{\Gamma}$ if and only if $\Gamma$ is a good group. Because we assumed that $X$ is smooth $Et(X)$ is weakly homotopy equivalent to $B\widehat{\pi_1(X(\mathbb C))}$ by Theorem \ref{complex}. On the other hand the profinite completion of $\pi_1(X(\mathbb C))$ is isomorphic to $\pi_1(X)$ by the Grauert-Remmert theorem. The claim is now clear.
\end{proof}
\begin{rem} It is important to note that the condition requiring the fundamental group to be good is not only sufficient, but also necessary. In particular there are algebraic varieties $X$ over $\mathbb C$ such that $X(\mathbb C)$ has the homotopy type of the Eilenberg-MacLane space $B\pi_1(X(\mathbb C))$, but the group $\pi_1(X(\mathbb C))$ is not good, therefore $Et(X)$ is not an Eilenberg-MacLane space. For an important class of examples see Lemma 3.16 of \cite{Mo} on page 146.
\end{rem}
\begin{prop}\label{notproper} Let $X$ be a smooth geometrically irreducible quasi-projective variety over an algebraically closed field $K$ of characteristic zero and let $F$ be another algebraically closed field containing $K$. Then $Et(X)$ is weakly homotopy equivalent to $Et(X_F)$.
\end{prop}
\begin{proof} This claim is a special case of Corollary 12.12 of \cite{AM} on page 144 when $X$ is proper, using also Theorem 11.1 of \cite{AM} on page 124. We only need to add a little bit more when $X$ is not proper. By Hironaka's resolution of singularities there is a projective variety $Y$ over $K$ which contains $X$ as an open subvariety such that the complement $C\subset Y$ is a normal crossings divisor. By the tame invariance theorem the tame fundamental groups $\pi_1^C(Y)$ and $\pi_1^{C_F}(Y_F)$ are isomorphic. But since the base fields have characteristic zero we have $\pi_1^C(Y)\cong\pi_1(X)$ and $\pi_1^{C_F}(Y_F)\cong\pi_1(X_F)$. Therefore $\pi_1(X)\cong\pi_1(X_F)$, so the argument presented in \cite{AM} can be applied in this case, too.
\end{proof}
\begin{cor}\label{curvesandabelians} The following hold:
\begin{enumerate}
\item[$(a)$] Let $X$ be a smooth, geometrically connected curve over an algebraically closed field $K$ of characteristic zero which is not a projective curve of genus zero. Then $Et(X)$ is weakly homotopy equivalent to $B\pi_1(X)$.
\item[$(b)$] Let $X$ be an abelian variety over an algebraically closed field $K$ of characteristic zero. Then $Et(X)$ is weakly homotopy equivalent to $B\pi_1(X)$.
\end{enumerate}
\end{cor}
\begin{proof} Recall that a smooth, geometrically connected curve $Y$ defined over a field has type $(g,d)$ if $g$ is the genus of the smooth projective completion $Y^c$ of $Y$ and $d$ is the number of geometric points in the complement of $Y$ in $Y^c$. Let $X$ be a smooth, geometrically connected curve of type $(g,d)$ such that $(g,d)\neq(0,0)$ over an algebraically closed field $K$ of characteristic zero. There is a subfield $F\subset K$ which is finitely generated over $\mathbb Q$ and $X$ is already defined over $F$, that is, there is a smooth, geometrically connected curve $Y$ of type $(g,d)$ over $F$ whose base change to $K$ is $X$.

Choose an embedding $i:\overline F\rightarrow\mathbb C$ of fields. Then the base change $Y_{\mathbb C}$ of the curve $Y_{\overline F}$ to $\mathbb C$ with respect to this embedding is also a smooth, geometrically connected curve of type $(g,d)$. The topological space $Y_{\mathbb C}(\mathbb C)$ has the homotopy type of the Eilenberg-MacLane space $B\pi_1(Y_{\mathbb C}(\mathbb C))$. The topological fundamental group of a smooth, connected complex curve is good (this fact follows at once from part $(a)$ of problem 1 of \cite{Se} on page 15) so we get from Proposition \ref{good} that $Et(Y_{\mathbb C})$ is weakly homotopy equivalent to $B\pi_1(Y_{\mathbb C})$. By a repeated application of Proposition \ref{notproper} we get that $Et(X)=Et(Y_K)$ is weakly homotopy equivalent to $Et(Y_{\mathbb C})$ and hence $\pi_1(X)\cong\pi_1(Y_{\mathbb C})$ and so $Et(X)$ is weakly homotopy equivalent to $B\pi_1(X)$.

The proof of claim $(b)$ is essentially the same as the proof of claim $(a)$; we only need to add that finitely generated free abelian groups are good (see part $(d)$ of problem 2 of \cite{Se} on page 16) so $Et(A)$ is weakly homotopy equivalent to $B\pi_1(A)$ for every abelian variety $A$ defined over $\mathbb C$ by Propositions \ref{good} and \ref{notproper}.
\end{proof}

\section{Grothendieck's short exact sequence}

\begin{notn}\label{4.3} Let $X$ be a geometrically connected variety defined over $K$. Let $\eta$ be a $\overline K$-valued point of $X$. Then Grothendieck's short exact sequence of \'etale fundamental groups for $X$ is:
\begin{equation}\label{fundamentalsequence}
\CD1@>>>\pi_1(\overline X,\eta)@>>>
\pi_1(X,\eta)@>>>\Gamma_K@>>>1,\endCD
\end{equation}
which is an exact sequence of profinite groups in the category of topological groups. Every $K$-rational point $x\in X(K)$ induces a section $\Gamma_K\rightarrow\pi_1(X,\eta)$ of the sequence (\ref{fundamentalsequence}), well-defined up to conjugation. Let $\textrm{Sec}(X/K)$ denote the set of conjugacy classes of sections of (\ref{fundamentalsequence}) (in the category of profinite groups where morphisms are continuous homomorphisms). Then we have a map:
$$s_{X/K}:X(K)\longrightarrow\textrm{Sec}(X/K)$$
which sends every point $x\in X(K)$ to the corresponding conjugacy class of sections.
\end{notn}
\begin{defn}\label{6.2} For every characteristic open subgroup $N$ of $\pi_1(\overline X,\eta)$ consider the short exact sequence:
\begin{equation}\label{Nsequence}
\CD1@>>>\pi_1(\overline X,\eta)/N@>>>
\pi_1(X,\eta)/N@>>>\Gamma_K@>>>1,\endCD
\end{equation}
we get by dividing out (\ref{fundamentalsequence}) by $N$. Let $\textrm{Sec}(X/K,N)$ denote the set of conjugacy classes of sections of (\ref{Nsequence}). Let
$$s_{X/K,N}:X(K)\longrightarrow\textrm{Sec}(X/K,N)$$
denote the composition of $s_{X/K,N}$ and the natural forgetful map
$$\phi_{X/K,N}:\textrm{Sec}(X/K)\longrightarrow\textrm{Sec}(X/K,N).$$ 
Note that for every pair of characteristic open subgroups $N'\subseteq N$ of
$\pi_1(\overline X,\eta)$ the composition of $\phi_{X/K,N'}$ and the forgetful map
$\textrm{Sec}(X/K,N')\rightarrow\textrm{Sec}(X/K,N)$ is $\phi_{X/K,N}$. Therefore we may take the projective limit of the maps $\phi_{X/K,N}$ to get a map:
$$\phi_{X/K}=\lim_N\phi_{X/K,N}:\textrm{Sec}(X/K)\longrightarrow
\lim_N\textrm{Sec}(X/K,N).$$ 
where the limit is over the set of  characteristic open subgroups of
$\pi_1(\overline X,\eta)$ directed with respect to reverse inclusion.
\end{defn}
\begin{prop}\label{section_limit} The map $\phi_{X/K}$ is a bijection.
\end{prop}
\begin{proof} Let $r,s$ be two sections $\Gamma_K\to\pi_1(X,\eta)$ such that for every characteristic open subgroup $N$ of $\pi_1(\overline X,\eta)$ the compositions of $r$ and $s$ with the quotient map $\pi_1(X,\eta)\to\pi_1(X,\eta)/N$ are conjugates. Then for every such $N$ the set
$$C_N=\{g\in\pi_1(\overline X,\eta)|g^{-1}r(h)gs(h)^{-1}\in N\ (\forall h\in\Gamma_K)\}$$
is non-empty. Since the sets $C_N$ are closed in the compact topological space
$\pi_1(\overline X,\eta)$, their intersection:
$$\bigcap_{N}C_N=
\{g\in\pi_1(\overline X,\eta)|g^{-1}r(h)gs(h)^{-1}\in\bigcap_NN\ (\forall h\in\Gamma_K)\}$$
is also non-empty. Because in a topologically finitely generated profinite group, such as $\pi_1(\overline X,\eta)$, the intersection of all characteristic open subgroups is the identity element, we get that $r$ and $s$ are conjugates. Therefore $\phi_{X/K}$ is injective.

Now let $r$ be an element of $\lim_N\textrm{Sec}(X/K,N)$ and for every $N$ as above let $r_N\in\textrm{Sec}(X/K,N)$ be the image of $r$ under the projection $\lim_M\textrm{Sec}(X/K,M)\to\textrm{Sec}(X/K,N)$. For every positive integer $m$ let $N(m)$ be the intersection of all open subgroups of $\pi_1(\overline X,\eta)$ of index at most $m$. We are going to construct a section $s_m:\Gamma_K\to\pi_1(X,\eta)/N(m)$ whose conjugacy class is $r_{N(m)}$ for every $m$ by induction as follows. When $m=1$ then this section is just the identity. Assume now that $s_{m-1}$ is already constructed. Let $s_m'$ be a section $\Gamma_K\to\pi_1(X,\eta)/N(m)$ whose conjugacy class is $r_{N(m)}$. Because $r_{N(m)}$ maps to $r_{N(m-1)}$ under the forgetful map $\textrm{Sec}(X/K,N(m))\rightarrow\textrm{Sec}(X/K,N(m-1))$ we get that there is a $g\in\pi_1(\overline X,\eta)/N(m)$ such that the composition of $g^{-1}s_m'g$ and the quotient map $\pi_1(X,\eta)/N(m)\rightarrow\pi_1(X,\eta)/N(m-1)$ is $s_{m-1}$. Let $s_m$ be $g^{-1}s_m'g$. These sections are compatible and their limit is a section:
$$s:\Gamma_K\longrightarrow\lim_{m=1}^{\infty}\pi_1(X,\eta)/N(m)=
\pi_1(X,\eta)$$
whose image is $r$ under $\phi_{X/K}$. So the latter is surjective, too.
\end{proof}
\begin{defn}\label{6.2b}We say that $X$ is well-equipped with $K$-rational points if the map $s_{X/K,N}$ is surjective for every characteristic open subgroup $N$ of $\pi_1(\overline X,\eta)$. Note that for a different choice of a base point $\eta'$ there is an isomorphism between $\pi_1(X,\eta)$ and $\pi_1(X,\eta')$ which maps $\pi_1(\overline X,\eta)$ onto $\pi_1(\overline X,\eta')$, canonical up to conjugacy. Therefore the sets $\textrm{Sec}(X/K)$ and $\textrm{Sec}(X/K,N)$ are independent of the choice of the base point $\eta$, as the notation indicates.
\end{defn}
\begin{prop}\label{manypoints} The algebraic groups $GL_n$ and $PGL_n$ are well-equipped with $K$-rational points over any characteristic zero field $K$ and positive integer $n$.
\end{prop}
\begin{proof} Let $1$ denote the unit of $GL_n(\overline K)$ and $PGL_n(\overline K)$, too. The quotient map $p:GL_n\rightarrow PGL_n$ by the centre of $GL_n$ induces a surjection
$$\pi_1(p):\pi_1(\overline{GL}_n,1)\longrightarrow
\pi_1(\overline{PGL}_n,1),$$
and hence it will be enough to prove the claim for $GL_n$ only. Let $i:GL_1\rightarrow GL_n$ be the map which embeds $GL_1$ into $GL_n$ as scalar matrices with ones on the diagonal except on the upper left corner. This map induces an isomorphism:
$$\pi_1(i):\pi_1(\overline{GL}_1,1)\longrightarrow\pi_1(\overline{GL}_n,1).$$ Therefore it will be enough to prove the claim for $GL_1$ only. There is a natural isomorphism:
$$\pi_1(\overline{GL}_1,1)\cong\widehat{\mathbb Z},$$
and for every $k\in\mathbb N$ there is natural bijection: 
$$\textrm{\rm Sec}(GL_1/K,k\widehat{\mathbb Z})\cong
H^1(K,\mu_k),$$
where $\mu_k\subseteq\overline K^*$ is the module of $k$-th roots of unity. Moreover under this identification $s_{GL_1/K,k\widehat{\mathbb Z}}$ corresponds to the coboundary map furnished by Kummer theory. Since $H^1(K,\overline K^*)$ is zero by Hilbert's Theorem 90, the claim now follows.
\end{proof}
\begin{prop}\label{p1exactforbundles} Assume that $K$ is an algebraically closed field of characteristic zero. Let $G$ be a geometrically connected algebraic group over $K$ and let $f:X\rightarrow Y$ be a principal $G$-bundle over a geometrically connected smooth variety $Y$ over $K$. Let $x\in X$ and set $y=f(x)$. Then the sequence:
\begin{equation}\label{fibrationsequence}
\CD\pi_1(f^{-1}(y),x)@>>>
\pi_1(X,x)@>\pi_1(f)>>\pi_1(Y,y)@>>>1\endCD
\end{equation}
is exact, where the first map is induced by the inclusion $f^{-1}(y)
\subseteq X$.
\end{prop}
\begin{proof} There is a subfield $F\subset K$ which is finitely generated over $\mathbb Q$ and $G,X,Y$ and $f$ are already defined over $F$. Therefore by the invariance theorem for the (tame) \'etale fundamental group it will be sufficient to prove the claim for $\overline F$. By the axiom of choose there is an embedding $\overline F\rightarrow\mathbb C$ of fields, and hence we may assume that $K$ is $\mathbb C$ without the loss of generality, again by the invariance theorem. Because the map $X(\mathbb C)\rightarrow Y(\mathbb C)$ induced by $f$ is a Serre fibration, there is a short exact sequence
\begin{equation}\label{topseq}
\CD\pi_1(f^{-1}(y)(\mathbb C),x)@>>>
\pi_1(X(\mathbb C),x)@>>>\pi_1(Y(\mathbb C),y)@>>>1\endCD
\end{equation}
of topological fundamental groups of complex analytic spaces. The profinite completion functor is right exact, so the completion of (\ref{topseq}) is also exact. By the Grauert--Remmert theorem the latter is the sequence (\ref{fibrationsequence}).
\end{proof}
\begin{prop}\label{liftforbundles} Let $G$ be a geometrically connected algebraic group well-equipped with $K$-rational points over a field $K$ of characteristic zero. Let $f:X\rightarrow Y$ be a principal $G$-bundle over a smooth variety $Y$ over $K$ and let $x,y\in Y(K)$ be such that
\begin{enumerate}
\item[$(i)$] we have $s_{Y/K}(x)=s_{Y/K}(y)$,
\item[$(ii)$] the sets $f^{-1}(x)(K)$ and $f^{-1}(y)(K)$ are non-empty.
\end{enumerate}
Then for every characteristic open subgroup $N$ of $\pi_1(\overline X,\eta)$ (where $\eta\in X(\overline K)$ is arbitrary), there are two points $x_N\in f^{-1}(x)(K)$ and $y_N\in f^{-1}(y)(K)$ such that $s_{X/K,N}(x_N)=s_{X/K,N}(y_N)$.
\end{prop}
\begin{proof} Pick two points $x'\in f^{-1}(x)(K)$, $y'\in f^{-1}(y)(K)$ and let $\eta\in f^{-1}(x)(\overline K)$. Let $r\in s_{X/K}(x'), s\in s_{X/K}(y')$ be two sections of the short exact sequence (\ref{fundamentalsequence}). Because both $x'$ and $\eta$ lie in $f^{-1}(x)(K)$, we may assume that the image of $r$ is in $\pi_1(\overline {f^{-1}(x)},\eta)$ without the loss of generality. Now let $\theta=f(\eta)$ and let $r_0,s_0$ be the composition of $r, s$ and $\pi_1(f):\pi_1(X,\eta)\rightarrow\pi_1(Y,\theta)$, respectively. As $r_0\in s_{Y/K}(x), s_0\in s_{Y/K}(y)$ we get that these sections are conjugate. By Proposition \ref{p1exactforbundles} the map $\pi_1(f)$ is surjective, therefore we may assume that $r_0$ and $s_0$ are the same, by conjugating $s$, if it is necessary. This implies that $s$, just as $r$, lies in $\pi_1(f^{-1}(x),\eta)$ by Proposition \ref{p1exactforbundles}. Let $N$ be now a characteristic open subgroup of $\pi_1(\overline X,\eta)$ and let $M'$ be the pre-image of $N$ with respect to the map $\pi_1(\overline {f^{-1}(x)},\eta)\to\pi_1(\overline X,\eta)$. As $M'$ is open and $\pi_1(\overline {f^{-1}(x)},\eta)$ is topologically finitely generated there is a characteristic open subgroup $M$ of $\pi_1(\overline {f^{-1}(x)},\eta)$ lying in $M'$. Because $f^{-1}(x)$ has a $K$-rational point, it is isomorphic to $G$. Therefore it is well-equipped with $K$-rational points, so there is an $x_N\in f^{-1}(x)(K)$ such that the composition of $s$ and the quotient map $\pi_1(f^{-1}(x),\eta)\rightarrow\pi_1(f^{-1}(x),\eta)/M$ lies in $s_{f^{-1}(x)/K,M}(x_N)$. If we set $y_N=y'$ then it is clear that the pair $x_N,y_N$ satisfies the required properties. 
\end{proof}
\begin{prop}\label{liftforbundles2} Let $G$ be a geometrically connected algebraic group well-equipped with $K$-rational points over a field $K$ of characteristic zero such that $\pi_1(\overline G)$ is finite. Let $f:X\rightarrow Y$ be a principal $G$-bundle over a smooth variety $Y$ over $K$ and let $x,y\in Y(K)$ be such that
\begin{enumerate}
\item[$(i)$] we have $s_{Y/K}(x)=s_{Y/K}(y)$,
\item[$(ii)$] the sets $f^{-1}(x)(K)$ and $f^{-1}(y)(K)$ are non-empty.
\end{enumerate}
Then there are two points $x'\in f^{-1}(x)(K)$ and $y'\in f^{-1}(y)(K)$ such that $s_{X/K}(x')=s_{X/K}(y')$.
\end{prop}
\begin{proof} The proof is the same as above, except that we look at sections of the full Grothendieck short exact sequence (\ref{fundamentalsequence}) for $X$. We leave the details to the reader. 
\end{proof}
\begin{rem}\label{pointsremark} By Proposition \ref{manypoints} we may apply Proposition \ref{liftforbundles} to principal $GL_n$-bundles. Note that $SL_n\rightarrow PSL_n$ is a finite \'etale cover and $SL_n$ is simply connected. As $PSL_n$ and $PGL_n$ are isomorphic, we get that $\pi_1(\overline{PGL}_n)$ is finite. Therefore by Proposition \ref{manypoints} we may apply Proposition \ref{liftforbundles2} to principal $PGL_n$-bundles.
\end{rem}

\section{Basic consequences of obstruction theory}

\begin{defn}\label{7.1} For every $n\in\mathbb N$ by functoriality we get a natural map:
$$\iota^n_{X/K}:X(K)\rightarrow 
[Et_{/K}(\textrm{Spec}(K)),Et_{/K}(X)]
\rightarrow
[Et_{/K}(\textrm{Spec}(K))^{\natural},Et^n_{/K}(X)]$$
where the second map is furnished by applying the Postnikov tower functor and composing with the $n$-th truncation map $Et_{/K}(X)^{\natural}\rightarrow Et^n_{/K}(X)$. We will denote $[Et_{/K}(\textrm{Spec}(K))^{\natural},Et^n_{/K}(X)]$ by the symbol $X^n(hK)$. The $n$-th truncation map $S^{\natural}\rightarrow P_n(S)$ is natural, and hence it induces a natural map
$$h^n_{X/K}:X(hK)\longrightarrow X^n(hK)$$
such that $\iota^n_{X/K}=h^n_{X/K}\circ\iota_{X/K}$. For every positive integer $n$ let $\sim_n$ denote the following equivalence relation on $X(hK)$: for every pair $x,y\in X(hK)$ we have $x\sim_ny$ if and only if $h^n_{X/K}(x)=h^n_{X/K}(y)$. It is clear that the equivalence relation $\sim_{n+1}$ is finer than the equivalence relation $\sim_n$. 
\end{defn}
\begin{defn}\label{7.2} For every pointed $K$-scheme $(X,x)$ let
$$(X,x)(hK)=Et_{/K}(X,x)(E\Gamma_K).$$
Let $\chi(X,x):(X,x)(hK)\to X(hK)$ denote the bijection induced by the map $i(X,x)$ (see Proposition \ref{comp1}). Then we have a unique natural map $\iota_{(X,x)/K}:X(K)\to(X,x)(hK)$ such that the diagram: 
$$\xymatrix{
  &  (X,x)(hK)\ar[dd]^{\chi(X,x)} \\
  X(K)\ar[ru]^{\iota_{(X,x)/K}}\ar[rd]_{\iota_{X/K}}  & \\
  & X(hK)}$$
is commutative. 
\end{defn}
\begin{defn}\label{7.3} We are going to need a variant of the equivalence relations $\sim_n$ for pointed relative \'etale homotopy types, too. Let $(X,x)$ be a pointed $K$-scheme as above. We will denote $Et^n_{/K}(X,x)(E\Gamma_K)$ by the symbol $(X,x)^n(hK)$. The $n$-th truncation map furnishes a natural map
$$h^n_{(X,x)/K}:(X,x)(hK)\longrightarrow(X,x)^n(hK).$$
By slight abuse of notation for every positive integer $n$ let $\sim_n$ denote the following equivalence relation on $(X,x)(hK)$: for every pair $x,y\in(X,x)(hK)$ we have $x\sim_ny$ if and only if $h^n_{(X,x)/K}(x)=h^n_{(X,x)/K}(y)$. 
\end{defn}
This notation is justified because
\begin{lemma}\label{comp1b} Under the map $\chi(X,x)$ the equivalence relation $\sim_n$ on $(X,x)(hK)$ corresponds to the the equivalence relation $\sim_n$ on $X(hK)$.
\end{lemma}
\begin{proof} For every \'etale hypercovering $H_*$ of $X$ there is a natural commutative diagram: 
$$\xymatrixcolsep{5pc}\xymatrix{
\pi_{0/K}(H_*)\ar[r]^{i(H_*,x)}\ar[d]  &  \pi_{0/K}(H_*,x)\ar[d] \\
P_n(\pi_{0/K}(H_*))\ar[r]^{P_n(i(H_*,x))}  & P_n(\pi_{0/K}(H_*,x)),}$$
where we use the notation of \ref{3.6b}, and where the vertical maps are induced by truncation. Since the upper horizontal map is a {\it weak} equivalence, so is the lower horizontal map, so it induces a bijection on homotopy fixed points. The claim now follows by taking the limit, similarly to the proof of Proposition \ref{comp1}.
\end{proof}
\begin{notn}\label{convention} For every profinite group $\Gamma$ and every pro-discrete $\Gamma$-module $M$ let $H^k(\Gamma,M)$ denote the projective limit of the (continuous) cohomology groups $H^k(\Gamma,N)$ where $N$ runs through the directed system of discrete quotients of $M$.  For the sake of simple notation for the rest of the paper for every field $K$ and for every pro-discrete $\Gamma_K$-module $M$ let $H^k(K,M)$ denote the group $H^k(\Gamma_K,M)$ introduced above. Note that these groups commute with projective limits. 
\end{notn}
\begin{thm}\label{obstructiontheory} Assume that $X$ is a smooth geometrically connected variety over $K$ and assume that $X(K)$ is non-empty. For every $x,y\in X(hK)$ the following holds:
\begin{enumerate}
\item[$(i)$] we have $x=y$ if and only if $x\sim_ny$ for every $n\in\mathbb N$,
\item[$(ii)$] there is a natural bijection:
$$j_{X/K}:X^1(hK)\longrightarrow\textrm{\rm Sec}(X/K)$$
such that for every $p\in X(K)$ we have $j_{X/K}(\iota^1_{X/K}(p))=s_{X/K}(p)$,
\item[$(iii)$] for every positive integer $n$ if $x\sim_ny$ then there exists a natural obstruction class $\delta^X_n(x,y)\in H^{n+1}(K,\pi_{n+1}(\overline X))$ such that
$x\sim_{n+1}y$ if and only if $\delta^X_n(x,y)=0$.
\end{enumerate}
\end{thm}
\begin{proof} The first claim is an immediate consequence of the definition of $Et_{/K}(X)^{\natural}$. Next we are going to prove $(ii)$. Let $f:Y\to X$ be a torsor under a finite \'etale group $\mathcal G$ over $K$. Then $\overline K$-valued points $G=\mathcal G(\overline K)$ of $\mathcal G$ form a finite group equipped with an action of $\Gamma_K$.  Let $C_*$ denote the \'etale \v Cech hypercovering:
$$\xymatrix{
 Y\ar@/^.5pc/@<1ex>[r] & Y\times_XY \ar@<-.5ex>[l] \ar@<.5ex>[l]
 & Y\times_XY\times_XY \ar@<-.7ex>[l] \ar[l] \ar@<.7ex>[l]
 & \cdots \ar@<-1ex>[l] \ar@<-.3ex>[l] \ar@<.3ex>[l] \ar@<1ex>[l]}$$
generated by the cover $Y\rightarrow X$. It is explained at the beginning of section 9 of \cite{HS} that there is a natural map:
$$c_Y:\pi_{0/K}(C_*)(E\Gamma_K)\longrightarrow H^1(\Gamma_K,G).$$
Moreover by Lemma 9.1 of \cite{HS} for every $p\in X(K)$ the image of the corresponding homotopy fixed point in $\pi_{0/K}(C_*)(E\Gamma_K)$ with respect to $c_Y$ is the element which classifies the $\mathcal G$-torsor $Y_p=f^{-1}(p)$. When $Y$ is geometrically connected there is a {\it weak} equivalence $\pi_{\Gamma}(C_*)\to BG$, where we equip the latter with the tautological $\Gamma_K$-action, and hence the map $c_Y$ is a bijection (see the discussion after Lemma 9.7 of \cite{HS}). Moreover the composition of the natural map $X(hK)\to\pi_{0/K}(C_*)(E\Gamma_K)$ and $c_Y$ factors through $h^1_{X/K}$.

Let $\eta\in X(K)$ be a $K$-rational point and let $\overline{\eta}$ denote the $\overline K$-valued point associated to $\eta$ (see Notation \ref{3.9}). Fix an element $s$ of $s_{X/K}(\eta)$. For every characteristic open subgroup $N$ of $\pi_1(X,\overline x)$ let $N'\subseteq\pi_1(X,\overline x)$ be the subgroup generated by $N$ and the image of $s$. Since $N'$ is an open subgroup there is a connected finite \'etale cover $f_N:Y_N\to X$ such that the image of $\pi_1(Y_N)$ with respect to $\pi_1(f_N)$ is $N'$. Let $\mathcal G_N$ denote the unique finite \'etale group over $K$ such that
$\mathcal G_N(\overline K)$ is
$$\pi_1(\overline X,\overline{\eta})/N\cong\pi_1(X,\overline{\eta})/N'$$
equipped with its natural $\Gamma_K$-action induced by conjugation. Then $f_N:Y_N\to X$ is a torsor under $\mathcal G_N$, and hence by applying the construction above to the \'etale \v Cech hypercovering generated by $f_N$ we get a map
$$j_{X/K,N}:X^1(hK)\longrightarrow\textrm{\rm Sec}(X/K,N)$$
such that for every $p\in X(K)$ we have $j_{X/K,N}(\iota^1_{X/K}(p))=s_{X/K,N}(p)$. According to Proposition \ref{section_limit} by taking the limit over every characteristic open subgroup $N$ of $\pi_1(X,\overline{\eta})$ we get a continuous map
$$j_{X/K}:X^1(hK)\longrightarrow\textrm{\rm Sec}(X/K)$$
between compact Hausdorff topological spaces such that for every $p\in X(K)$ we have $j_{X/K}(\iota^1_{X/K}(p))=s_{X/K}(p)$. By Lemma 9.11 of \cite{HS} this map is a bijection.

For every connected pointed $\Gamma_K$-space $S$ which has only finitely many non-trivial homotopy groups there is a natural spectral sequence for homotopy groups of homotopy fixed points: 
$$E^{s,t}_2=H^s(K,\pi_t(S))
\Rightarrow\pi_{t-s}(S^{h\Gamma_K})$$
constructed by Goerss in \cite{Go} (see Theorem B on page 189). Therefore for objects of this category there are a bijection and a natural obstruction class of the type described in the last two claims. So claim $(iii)$ follows from Lemma \ref{comp1b} and Proposition \ref{truncationcomparison2} by applying Goerss's results to $Et_{/K}(X,\eta)$. 
\end{proof}
\begin{lemma}\label{sectionanabelian} Assume that $X$ is a smooth geometrically connected variety over $K$ and $Et(\overline X)$ is weakly homotopy equivalent to $B\pi_1(\overline X)$. Also suppose that $X(K)\neq\emptyset$. Then for every $x,y\in X(hK)$ we have $x=y$ if and only if $j_{X/K}(x)=j_{X/K}(y)$.
\end{lemma}
\begin{proof} All the higher homotopy groups of $\overline X$ vanish, so the claim is immediate from the theorem above.
\end{proof}
We say that two points $x,y\in X(K)$ are directly $A$-equivalent if there is a map $f:\mathbb A_K^1\rightarrow X$ of $K$-varieties such that $f(0)=x$ and $f(1)=y$. The $A$-equivalence on $X(K)$ is the equivalence relation generated by direct $A$-equivalence. 
\begin{prop}\label{a1} Assume that $K$ has characteristic zero. Then for every $X$ over $K$ the map
$$\iota_{X/K}:X(K)\longrightarrow X(hK)$$
factors through $A$-equivalence.
\end{prop}
\begin{proof} It will be sufficient to show that for every two points $x,y\in X(K)$ which are directly $A$-equivalent we have $\iota_{X/K}(x)=\iota_{X/K}(y)$. Let $f:\mathbb A_K^1\rightarrow X$ be a morphism of $K$-varieties such that $f(0)=x$ and $f(1)=y$. Both $\iota_{X/K}(x)$ and $\iota_{X/K}(y)$ lie in the image of the map $f_*:\mathbb A^1_K(hK)\rightarrow X(hK)$. By Corollary \ref{curvesandabelians} and Lemma \ref{sectionanabelian} the set $\mathbb A^1_K(hK)$ consists of one element, since $\mathbb A^1_{\overline K}$ has trivial \'etale fundamental group. The claim is now clear.
\end{proof}
\begin{rem}\label{tamagawa} The validity of a such a claim was already suggested by To\"en in \cite{To}, but his original claim is not true as stated in positive characteristic. In the special case when $K$ is a finite field it was observed in \cite{Ta} (see Proposition 2.8 on pages 151-152) that the map $s_{\mathbb A^1_K/K}$, and hence the map $\iota_{\mathbb A^1_K/K}$, is injective. Therefore it is not true in general that for two points $x,y\in X(K)$ which are directly $A$-equivalent we have $\iota_{X/K}(x)=\iota_{X/K}(y)$ when $K$ has positive characteristic.
\end{rem}
  
\section{The Manin pairing}

\begin{notn} For every object $X:I\rightarrow\textrm{Ho}(\Gamma-SSets)$ of $\textrm{Pro}-\textrm{Ho}(\Gamma-SSets)$ such that $I$ is countable let $X_{h\Gamma}\in\textrm{ob}(\textrm{Pro}-\textrm{Ho}(SSets))$ denote the pro-homotopy quotient defined at the beginning of section 6.2 of \cite{HS}. Note that this construction can be applied to the objects $Et_{/K}(X),Et^n_{/K}(X)$ and $Et_{/K}(X)^{\natural}$ when $\Gamma=\Gamma_K$ is the absolute Galois group of the field $K$.
\end{notn}
\begin{prop} Let $K$ be a field and $X$ be a variety over $K$. Then there are natural isomorphisms: 
$$Et_{/K}(X)_{h\Gamma_K}\cong Et(X),\quad Et^n_{/K}(X)_{h\Gamma_K}
\cong Et^n(X),\quad Et_{/K}(X)^{\natural}_{h\Gamma_K}\cong Et(X)^{\natural}$$
in the category $\textrm{\rm Pro}-\textrm{\rm Ho}(SSets)$.
\end{prop}
\begin{proof} The second isomorphism is the content of Proposition 6.14 in \cite{HS}. The third isomorphism follows from the naturality of this isomorphism, and the first isomorphism is shown in the proof of the proposition mentioned above.
\end{proof}
\begin{defn} Let $X$ be again a variety over $K$. Note that by functoriality we get a natural map:
$$\lambda_{X/K}:X(K)\longrightarrow 
[Et(\textrm{Spec}(K)),Et(X)]
\longrightarrow
[Et(\textrm{Spec}(K))^{\natural},Et(X)^{\natural}]$$
where the second map is furnished by applying the Postnikov tower functor. By applying the pro-homotopy quotient functor and the proposition above we get that there is a natural map:
$$\kappa_{X/K}:X(hK)\longrightarrow[Et(\textrm{Spec}(K))^{\natural},Et(X)^{\natural}]$$
such that $\lambda_{X/K}=\kappa_{X/K}\circ\iota_{X/K}$.
\end{defn}
\begin{defn} By Corollary 10.8 of \cite{AM} on pages 122-123 there is a natural equivalence between the category of locally constant \'etale sheaves of finite abelian groups on $X$ and local coefficient systems of finite abelian groups on $Et(X)^{\natural}$ and under this equivalence the \'etale cohomology of $X$ with coefficients in a locally constant \'etale sheaf $\mathcal F$ of finite abelian groups is the same as the cohomology of $Et(X)^{\natural}$ with coefficients in the local coefficient system corresponding to $\mathcal F$. We will not distinguish these two categories in all that follows. In particular for a finite, \'etale group scheme $G$ over $K$ we will identify $H^*(Et(X)^{\natural},G)$ and $H^*(X,G)$.
\end{defn}
A basic, but important corollary of these observations is the following. Let $G$ and $X$ be as above, and let $c\in H^i(X,G)$ be a cohomology class for some $i\in\mathbb N$. 
\begin{lemma}\label{basicinvariant} Assume that $x,y\in X(K)$ are $H$-equivalent. Then the cohomology classes $x^*( c ),y^*( c )\in H^i(K,G)$ are equal.
\end{lemma}
\begin{proof} The lemma follows from the commutativity of the diagram:
$$\CD
X(K)@.\times@.H^i(X,G)@>>>H^i(K,G)\\
@V\lambda_{X/K}VV@.@|@|\\
[Et(\textrm{Spec}(K))^{\natural},Et(X)^{\natural}]
@.\times@.H^i(Et(X)^{\natural},G)
@>>>H^i(Et(\textrm{Spec}(K))^{\natural},G)
\endCD$$
where the horizontal maps are the pairings furnished by pull-back.
\end{proof}
For every $n\in\mathbb N$ not divisible by the characteristic of $K$ let $(\cdot,\cdot)_n$ denote the pairing:
$$(\cdot,\cdot)_n:X(hK)
\times H^2(X,\mu_n)
\rightarrow H^2(K,\mu_n)$$
given by the rule $(x,c)_n=\kappa_{X/K}(x)^*(c)$.
\begin{lemma} The diagram:
\begin{equation}\label{5.5.1}
\CD
X(K)@.\times@.\textrm{\rm Br}(X)@>(\cdot,\cdot)>>\textrm{\rm Br}(K)\\
@V\iota_{X/K}VV@.@AAA@AA\alpha_nA\\
X(hK)
@.\times@.H^2(X,\mu_n)
@>(\cdot,\cdot)_n>>\textrm H^2(K,\mu_n)
\endCD
\end{equation}
commutes where the middle and right vertical arrows are induced by the inclusion $\mu_n\subset\mathbb G_m$ of sheaves.
\end{lemma}
\begin{proof} The lemma is immediate from the functoriality of the constructions involved.
\end{proof}
\begin{lemma}\label{5.6} Let $X$ be a geometrically connected variety over $K$ and $x,y\in X(K)$. Then $\iota_{X/K}(x)\sim_1\iota_{X/K}(y)$ if and only if for every 
finite, \'etale map $f:Y\rightarrow X$ of geometrically connected varieties over $K$ such that there is an $\widetilde x\in Y(K)$ with the property $f(\widetilde x)=x$ there is a $\widetilde y\in Y(K)$ such that $f(\widetilde y)=y$.
\end{lemma}
\begin{proof} By part $(ii)$ of Theorem \ref{obstructiontheory} we have $\iota_{X/K}(x)\sim_1\iota_{X/K}(y)$ if and only if $s_{X/K}(x)=s_{X/K}(y)$. It is well-known that the latter condition is equivalent to the second condition of the claim (see for example \cite{Ta}).
\end{proof}
\begin{lemma}\label{5.7} Let $f:Y\rightarrow X$ be a finite, \'etale map of varieties over $K$. Assume that $x,y\in X(K)$ are $H$-equivalent and there is an $\widetilde x\in Y(K)$ such that $f(\widetilde x)=x$. Then there is a $\widetilde y\in Y(K)$ such that $\widetilde x,\widetilde y$ are $H$-equivalent and $f(\widetilde y)=y$.
\end{lemma}
\begin{proof} The connected component $X'$ of $X$ on which $x$ lies is geometrically connected. By Proposition \ref{connectedcomponents} the point $y$ must lie on the same component. The connected component $Y$ of $Y$ on which $\widetilde x$ lies is also geometrically connected and the restriction $f|_{Y'}:Y'\rightarrow X'$ is a finite, \'etale map. Therefore we may assume without the loss of generality that $X$ and $Y$ are geometrically connected. Hence by Lemma \ref{5.6} there is a $\widetilde y\in Y(K)$ such that $f(\widetilde y)=y$, and we may even assume that $\iota_{X/K}(\widetilde x)\sim_1\iota_{X/K}(\widetilde y)$. 

It will be enough to show that $\widetilde x$ and $\widetilde y$ are $H$-equivalent. We will prove that $\iota_{X/K}(\widetilde x)\sim_n\iota_{X/K}(\widetilde y)$ for every $n\geq2$ by induction. Since the map $\pi_n(f):\pi_n(\overline Y)\rightarrow\pi_n(\overline X)$ is an isomorphism for every $n\geq2$ by Proposition \ref{imsmarter} we get that the induced map:
$$H^n(\pi_n(f)):H^n(K,\pi_n(\overline Y))\rightarrow H^n(K,\pi_n(\overline X))$$
is also an isomorphism. By naturality of the obstruction classes we have:
$$H^n(\pi_n(f))(\delta^Y_n(\widetilde x,\widetilde y))=\delta^X_n(x,y).$$
Since the right-hand side is zero we get that $\delta^Y_n(\widetilde x,\widetilde y)$ is also zero. 
\end{proof}
\begin{prop} Let $X$ be a regular variety over $K$ and assume that $x,y\in X(K)$ are $H$-equivalent. Then $x$ and $y$ are \'etale-Brauer equivalent.
\end{prop}
\begin{proof} Let $Y\rightarrow X$ be a finite, \'etale morphism of varieties over $K$ such that there is an $\widetilde x\in Y(K)$ mapping to $x$. By Lemma \ref{5.7} there is a $\widetilde y\in Y(K)$ such that $\widetilde x,\widetilde y$ are $H$-equivalent and $f(\widetilde y)=y$. It will be enough to show that $\widetilde x,\widetilde y$ are Brauer equivalent. Because $Y$ is the finite \'etale cover of a regular variety, it is also regular, so it will be enough show that given a regular variety $X$ over $K$ for every pair of points $x,y\in X(K)$ which are $H$-equivalent are also Brauer equivalent. Because $X$ is regular, the group $H^2(X,\mathbb G_m)$ is torsion (see Proposition 1.4 of \cite{Gr1} on page 291). Therefore for every $b\in H^2(X,\mathbb G_m)$ there is a natural number $n\in\mathbb N$ and a $c\in H^2(X,\mu_n)$ such that $b$ is the image of $c$ under the natural map $H^2(X,\mu_n)\rightarrow H^2(X,\mathbb G_m)$. The claim now follows from the commutativity of the diagram (\ref{5.5.1}).
\end{proof}

\section{Brauer equivalence versus \'etale-Brauer equivalence}

\begin{defn} Let $E$, $E'$ be two elliptic curves defined over a field $K$ whose characteristic is not two, and let $t\in E'(K)$ be a point of order two. Let 
$\rho:E\rightarrow E$ be the multiplication by $-1$ map, and let
$\sigma:E'\rightarrow E'$ be the translation by $t$. Let $X$ denote the quotient of $E\times E'$ by the fixed point free involution $(\rho,\sigma)$. Then $X$ is a smooth projective geometrically irreducible surface over $K$. We call such surfaces bielliptic. 
\end{defn}
\begin{prop}\label{bi1} Let $K$ be a finite extension of $\mathbb Q_p$ and let $X$ be a bielliptic surface over $K$. Then the map $s_{X/K}$ is injective.
\end{prop}
\begin{proof} Let $x,y\in X(K)$ be two different points and let $\eta$ be
$\overline K$-valued point of $X$. Let $\mathcal K\subseteq\pi_1(\overline X,\eta)$ be the characteristic subgroup such that the quotient
$\pi_1(\overline X,\eta)/\mathcal K$ is the maximal $2$-torsion abelian quotient of $\pi_1(\overline X,\eta)$. Fix an element $s$ of $s_{X/K}(x)$ and
let $\mathcal K'\subseteq\pi_1(X,\eta)$ be the subgroup generated by $\mathcal K$ and the image of $s$; this is an open subgroup. Let $f:Y\rightarrow X$ be the connected finite \'etale cover such that the image of $\pi_1(Y)$ with respect to $\pi_1(f)$ is $\mathcal K'$. Since $\overline Y$ is a finite \'etale cover of the abelian variety $\overline{E\times E'}$, where we use the notation of the definition above, it is also an abelian variety. Therefore $Y$ is a principal homogenous space over an abelian variety defined over $K$. 

By construction there is a $x'\in Y(K)$ such that $f(x')=x$. So $Y$ has a $K$-rational point, and hence it is also an abelian variety over $K$. If $y$ does not have a lift to a $K$-valued point of $Y$ then $s_{X/K}(x)\neq s_{X/K}(y)$. So we may assume that there is a $y'\in Y(K)$ such that $f(y')=y$. Let $G=\pi_1(X,\eta)/\pi_1(Y,\eta)$. Since we may take a finite extension of $K$ during the proof of injectivity of $s_{X/K}$, we may assume that the action of $\Gamma_K$ on $G$ is trivial without the loss of generality. In this case this finite group is the Galois group of the connected finite \'etale cover $f:Y\rightarrow X$. It also acts on $\textrm{Sec}(Y/K)$ and two elements of $\textrm{Sec}(Y/K)$ are in the same $G$-orbit if and only if they have the same image under the map $\textrm{Sec}(Y/K)\to\textrm{Sec}(X/K)$ induced by $f$. Moreover the section map $s_{Y/K}:Y(K)\to\textrm{Sec}(Y/K)$ is $G$-equivariant. 
Therefore if $s_{X/K}(x)=s_{X/K}(y)$ then $s_{Y/K}(x')=s_{Y/K}(g(y'))$ for some $g\in G$. As we already noted for abelian varieties over $K$ the section map is injective, so $x'=g(y')$ in this case, which implies that $x=y$. This is a contradiction, therefore $s_{X/K}(x)$ is different from $s_{X/K}(y)$.
\end{proof}
We continue to use the notation which we introduced above. Let $D$ denote the quotient of $E'$ by the fixed point free involution $\sigma$ and let $g:X\rightarrow D$ be the quotient map.
\begin{prop}\label{bi3} Assume that $K$ is a finite extension of $\mathbb Q_p$. Then the map $g^*:H^2(D,\mathbb G_m)\rightarrow H^2(X,\mathbb G_m)$ induced by $g$ has finite cokernel.
\end{prop}
\begin{proof} For every variety $Y$ over $K$ the Hochschild-Serre spectral sequence:
$$E^2_{p,q}=H^p(K,H^q(\overline Y,\mathbb G_m))\Rightarrow
H^{p+q}(Y,\mathbb G_m)$$
furnishes on $H^2(Y,\mathbb G_m)$ a natural filtration:
\begin{equation}
0=E^2_3\subseteq E^2_2
\subseteq E^2_1\subseteq
E^2_0=H^2(Y,\mathbb G_m)
\end{equation}
such that
$$E^{\infty}_{p,2-p}\cong E^2_p/E^2_{p+1},\quad p=0,1,2.$$
The members $E^2_2$ and $E^2_1$ are usually denoted by $\textrm{Br}_0(Y)$ and $\textrm{Br}_1(Y)$, respectively. Because $H^3(K,\mathbb G_m)=0$ (see Proposition 15 of \cite{Se} on page 93), the coboundary map:
$$d^2_{1,1}:E^2_{1,1}=
H^1(K,H^1(\overline Y,\mathbb G_m))
\longrightarrow H^3(K,\mathbb G_m)$$
is zero and therefore:
$$E^{\infty}_{1,1}=E^3_{1,1}=\textrm{Ker}(d^2_{1,1})=H^1(K,H^1(\overline Y,\mathbb G_m)).$$
In short we have a natural exact sequence:
$$\CD0@>>>
\textrm{Br}_0(Y)@>>>\textrm{Br}_1(Y)@>>>
H^1(K,H^1(\overline Y,\mathbb G_m))@>>>0.\endCD$$
The group $\textrm{Br}_0(Y)$ is the image of the natural map $\textrm{Br}(K)\rightarrow\textrm{Br}(Y)$, therefore the map $g^*:\textrm{Br}_0(D)\rightarrow
\textrm{Br}_0(X)$ is surjective. The Hochschild-Serre spectral sequence
furnishes a natural injection:
$$\CD\textrm{Br}(Y)/\textrm{Br}_1(Y)@>>>
\textrm{Br}(\overline Y).\endCD$$
Since  $\textrm{Br}(\overline X)$ is dual to the torsion subgroup of the N\'eron-Severi group $\textrm{NS}(\overline X)$ in our case (see \cite{Sk2} on page 403), which is finite, it will be enough to show that the map
$$\CD g_*:H^1(K,H^1(\overline D,\mathbb G_m))
@>>>H^1(K,H^1(\overline X,\mathbb G_m))\endCD$$
induced by $g$ has finite cokernel. Since $g$ is the Albanese map for $X$ (see {\it loc. \!\!\!sit.}), the map $\textrm{Pic}^0(\overline D)\rightarrow
\textrm{Pic}^0(\overline X)$ induced by $g$ is an isomorphism. Therefore, by looking at the cohomological long exact sequence associated 
to the short exact sequence of $\Gamma_K$-modules:
$$\CD0@>>>\textrm{Pic}^0(\overline X)(\overline K)
@>>> H^1(\overline X,\mathbb G_m)
@>>>\textrm{NS}(\overline X)@>>>0,\endCD$$
we are reduced to show that $H^1(K,\textrm{NS}(\overline X))$ is finite. The abelian group $\textrm{NS}(\overline X)$ is finitely generated, so there is a finite Galois extension $L|K$ such that the action of Gal$(\overline L|L)$ on $\textrm{NS}(\overline X)$ is trivial. The abelianization of Gal$(\overline L|L)$ is isomorphic to the profinite completion of $L^*$, so it is topologically finitely generated. Therefore
$$H^1(L,\textrm{NS}(\overline X))\cong
\textrm{Hom}(\textrm{Gal}(\overline L|L),\textrm{NS}(\overline X))$$
is finite. Therefore the inflation map:
$$H^1(\textrm{Gal}(L|K),\textrm{NS}(\overline X))
\longrightarrow H^1(K,\textrm{NS}(\overline X))$$
has finite cokernel. 
Since $\textrm{Gal}(L|K)$ is finite we get that $H^1(K,\textrm{NS}(\overline X))$ is finite.
\end{proof}
Let $X$ be a smooth variety over a field $K$ and let $b\in H^2(X,\mathbb G_m)$. We say that $x,y\in X(K)$ are $b$-equivalent if $x^*(b)=y^*(b)$. This defines an equivalence relation of $X(K)$, which we will call $b$-equivalence.
\begin{prop}\label{bi2} Assume that $K$ is a finite extension of $\mathbb Q_p$. Then $b$-equivalence classes are open in the $p$-adic topology.
\end{prop}
\begin{proof} Let $x\in X(K)$. It will be enough to show that $x$ has a $p$-adically open neighbourhood $U$ in $X(K)$ such that $x$ and $y$ is $b$-equivalent for every $y\in U$. We may assume that $X$ is affine by taking a Zariski-open neighbourhood of $x$. By a theorem of Gabber (see Theorem 1.1 of \cite{Jo}) there is an Azumaya algebra $A$ of some rank $n$ on $X$ which represents $b$. Let $\pi:Y\rightarrow X$ be the $PGL_n$-torsor corresponding to $A$, that is, the torsor whose class in $H^1(X,PGL_n)$ is the same as the class of $A$. Let $\sigma\in H^1(K,PGL_n)$ be the class of the fibre of $Y$ over $x$ and let $\pi^{\sigma}:Y^{\sigma}\rightarrow X$ be the twist of $\pi:Y\rightarrow X$ by $\sigma$. Then the fibre of $Y^{\sigma}$ over $x$ is a trivial $PGL_n$-torsor. Therefore there is a $x'\in Y^{\sigma}(K)$ such that $\pi^{\sigma}(x')=x$. Because $\pi^{\sigma}$ is a submersion there is an $p$-adically open neighbourhood $U$ of $x$ in $X(K)$ and a $p$-adically analytical section $s:U\rightarrow Y^{\sigma}(K)$ of $\pi^{\sigma}$ with $s(x)=x'$. Therefore for every $y\in U$ the fibre of $Y^{\sigma}$ over $y$ has a $K$-rational point, so it is a trivial $PGL_n$-torsor. The claim is now clear.
\end{proof}
\begin{thm}\label{bielliptic} Let $K$ be as above and let $X$ be a bielliptic surface over $K$. Then \'etale--Brauer equivalence is strictly finer than Brauer equivalence on $X(K)$.
\end{thm}
Of course our choice of example is motivated by Skorobogatov's classical paper \cite{Sk2}, and the result above can be considered its natural local counterpart. (Also compare with \cite{Ha} which uses similar ideas.)
\begin{proof} Let $X$ be a bielliptic surface over $K$ and let $g:X\rightarrow D$ be the map introduced above. By Proposition \ref{bi1} the \'etale-Brauer equivalence-classes of $X(K)$ consists of points, so it will be enough to show that the Brauer equivalence-classes of $X(K)$ are infinite. Let $r\in D(K)$. Note that every $x,y\in g^{-1}( r  )(K)$ are $b$-equivalent for every $b\in H^2(X,\mathbb G_m)$ in the image of the map $g^*:H^2(D,\mathbb G_m)\rightarrow H^2(X,\mathbb G_m)$. Therefore finitely many Brauer equivalence classes intersect $g^{-1}( r  )(K)$ by Propositions \ref{bi3} and \ref{bi2}. By the inverse function theorem both $D$ and $g^{-1}( r  )$ have infinitely many $K$-valued points, so the claim holds.
\end{proof}
In the rest of this section we study the somewhat independent question of the surjectivity of $\iota_{X/K}$ over $p$-adic fields.
\begin{notn} Let $Groups$ denote the category of groups. Let $\mathbb N_d$ denote the category whose objects are positive integers and for every pair of objects $m$, $n\in\textrm{ob}(\mathbb N_d)$ the set of morphisms from $m$ to $n$ consists of the ordered pair $\phi_{m,n}=(m,n)$, if $n|m$, and it is empty, otherwise. For every abelian group $A$ and natural number $n$ let $A[n]$ denote the subgroup of $n$-torsion elements of $A$ and let $A^{tor}$ denote the pro-group $A^{tor}:\mathbb N_d\rightarrow
Groups$ given by the rule:
$$A^{tor}(n)=A[n]$$
such that for every pair of positive integers $m,n$ such that $n|m$ the homomorphism:
$$A^{tor}(\phi_{m,n}):A[m]\rightarrow A[n]$$
is the multiplication by $m/n$ map.
\end{notn}
\begin{prop}\label{section} Let $K$ be a finite extension of $\mathbb Q_p$. The following holds:
\begin{enumerate}
\item[$(a)$] for every smooth, geometrically connected projective curve $X$ of genus at least two over $K$ the map $\iota_{X/K}$ is injective, and it is surjective if the local version of Grothendieck's section conjecture holds for $X$. 
\item[$(b)$] For every abelian variety $X$ over $K$ the map $\iota_{X/K}$ is injective, and it is surjective if and only if $X$ is zero dimensional. 
\end{enumerate}
\end{prop}
\begin{proof} First assume that $X$ is a smooth, geometrically connected projective curve of genus at least two over $K$. Recall that the local version of Grothendieck's section conjecture claims that the map $s_{X/K}$ is a bijection. We also know that in this case $s_{X/K}$ is injective. Therefore claim $(a)$ follows at once from Corollary \ref{curvesandabelians} and Lemma \ref{sectionanabelian}. Assume now that $X$ is an abelian variety over $K$. By Corollary \ref{curvesandabelians} and Lemma \ref{sectionanabelian} the map $j_{X/K}$ is a bijection. Moreover there is a natural bijection 
\begin{equation}\label{4.9.1a}
\textrm{\rm Sec}(X/K)\cong H^1(K,\prod_{\textrm{$l$ is prime}}T_l(X))
\end{equation}
where $T_l(X)$ denotes the $l$-th Tate module of $X$, and under this identification $s_{X/K}$ corresponds to the coboundary map furnished by Kummer theory. In particular $\iota_{X/K}$ is injective, and the cokernel of $s_{X/K}$ is $H^1(K,X)^{tor}$. By Corollary 3.4 of \cite{Mi2} on page 53 the groups $H^1(K,X)$ and $X^{\vee}(K)^{\wedge}$ are isomorphic, where the latter is the Pontryagin dual of the compact group $X^{\vee}(K)$ of $K$-valued points of the dual $X^{\vee}$ of $X$. Let $\mathcal O_K$ denote the valuation ring of $K$. The profinite group $X^{\vee}(K)$ is the direct sum of a finite group and $\dim(X^{\vee})=
\dim(X)$ copies of $\mathcal O_K$ by the inverse function theorem. Therefore
the group $H^1(K,X)^{tor}$ is zero if and only if $X$ is zero dimensional. Hence claim $(b)$ is true. 
\end{proof}

\section{\'Etale-Brauer equivalence versus $H$-equivalence}

\begin{defn}\label{10.1} Assume now that $K$ is a $p$-adic field and for every $n\in\mathbb N$ let $c_{K,n}:H^2(K,\mu_n)\rightarrow
\mathbb Z/n\mathbb Z$ be the isomorphism furnished by local Tate duality. For every geometrically irreducible variety $X$ defined over $K$ let
$$\{\cdot,\cdot\}_n:H^2(K,H_2(\overline X,\mathbb Z/n\mathbb Z))
\times H^0(K,H^2(\overline X,\mu_n))
\longrightarrow\mathbb Z/n\mathbb Z$$
be the bilinear pairing given by the rule:
$$\{x,y\}_n=c_{K,n}(x\cup y),$$
(for every $x\in H^2(K,H_2(\overline X,\mathbb Z/n\mathbb Z)),
y\in H^0(K,H^2(\overline X,\mu_n))$), where the cup product
$$\cup:H^2(K,H_2(\overline X,\mathbb Z/n\mathbb Z))\times
H^0(K,H^2(\overline X,\mu_n))
\longrightarrow H^2(K,\mu_n)$$
is induced by the evaluation pairing:
$$H_2(\overline X,\mathbb Z/n\mathbb Z)\times
H^2(\overline X,\mu_n)\longrightarrow\mu_n.$$
\end{defn}
The following lemma is immediate from the functoriality of the constructions involved:
\begin{lemma}\label{pairingnatural} Let $f:X\rightarrow Y$ be a morphism of geometrically
irreducible varieties over $K$. Then the diagram:
\begin{equation}
\CD
H^2(K,H_2(\overline X,\mathbb Z/n\mathbb Z))@.\times@.
H^0(K,H^2(\overline X,\mu_n))
@>\{\cdot,\cdot\}_n>>\mathbb Z/n\mathbb Z\\
@V H_2(f)_*VV@.@A(f^*)_*AA@|\\
H^2(K,H_2(\overline Y,\mathbb Z/n\mathbb Z))@.\times@.
H^0(K,H^2(\overline Y,\mu_n))
@>\{\cdot,\cdot\}_n>>\mathbb Z/n\mathbb Z\endCD
\end{equation}
commutes.\qed
\end{lemma}
\begin{notn}\label{10.3} Let $X$ be a geometrically irreducible variety over $K$ with $X(K)\neq\emptyset$. For every positive integer $n$ and for every $x,y\in X(hK)$ such that $x\sim_1y$ let $\delta^X_1(x,y)_n$ denote the image of the obstruction class $\delta^X_1(x,y)$ under the composition of the natural map:
$$H^2(K,\pi_2(\overline X))\rightarrow H^2(K,\pi_2(\overline X)/n\pi_2(\overline X)),$$
and  the homomorphism
$$H_{*,n}:H^2(K,\pi_2(\overline X)/n\pi_2(\overline X))
\longrightarrow H^2(K,H_2(\overline X,\mathbb Z/n\mathbb Z))$$
induced by the Hurewitz map:
$$H_n:\pi_2(\overline X)/n\pi_2(\overline X)
\longrightarrow H_2(\overline X,\mathbb Z/n\mathbb Z).$$
Let
$$\alpha_n:H^2(X,\mu_n)\rightarrow H^2(X,\mathbb G_m)$$
denote the map induced by the inclusion $\mu_n\rightarrow\mathbb G_m$, and finally let
$$\phi_n:H^2(X,\mu_n)\rightarrow H^0(K,H^2(\overline X,\mu_n))$$
be the map induced by base change. 
\end{notn}
\begin{lemma}\label{linebundle} Let $X$ be a geometrically irreducible smooth quasi-projective variety over $K$. For every $n\in\mathbb N$, for every $c\in H^2(X,\mu_n)$ and for every $x,y\in X(K)$ such that $\alpha_n(c)=0$ and
$\iota_{X/K}(x)\sim_1\iota_{X/K}(y)$ we have:
$$\{\delta^X_1(\iota_{X/K}(x),\iota_{X/K}(y))_n,\phi_n(c)\}_n=0.$$
\end{lemma}
\begin{proof} We will need to introduce the analogues of the concepts in Definition \ref{10.1} and Notation \ref{10.3} for hypercoverings. For every geometrically irreducible variety $Z$ defined over $K$ and for every \'etale hypercovering $H_*$ of $Z$ let
$$\{\cdot,\cdot\}_n:H^2(K,H_2(\pi_{0/K}(H_*),\mathbb Z/n\mathbb Z))
\times H^0(K,H^2(\pi_{0/K}(H_*),\mu_n))
\longrightarrow\mathbb Z/n\mathbb Z$$
also denote the bilinear pairing given by the rule:
$$\{a,b\}_n=c_{K,n}(a\cup b),$$
(for every $a\in H^2(K,H_2(\pi_{0/K}(H_*),\mathbb Z/n\mathbb Z)),
b\in H^0(K,H^2(\pi_{0/K}(H_*),\mu_n))$), where the cup product
$$\cup:H^2(K,H_2(\pi_{0/K}(H_*),\mathbb Z/n\mathbb Z))\times
H^0(K,H^2(\pi_{0/K}(H_*),\mu_n))
\longrightarrow H^2(K,\mu_n)$$
is induced by the evaluation pairing:
$$H_2(\pi_{0/K}(H_*),\mathbb Z/n\mathbb Z)\times
H^2(\pi_{0/K}(H_*),\mu_n)\longrightarrow\mu_n.$$

Assume now that $Z(K)\neq\emptyset$ and pick a point $z\in Z(K)$. In Definition \ref{3.6} we introduced a pointed simplicial $\Gamma_K$-set $\pi_{0/K}(H_*,z)$ such that there is a natural map $\pi_{0/K}(H_*)\to\pi_{0/K}(H_*,z)$ which is a weak equivalence. Let $a,b\in\pi_0(\pi_{0/K}(H_*)^{\Gamma_K})$ be such that their image is the same under the map
$$\pi_0(\pi_{0/K}(H_*)^{\Gamma_K})\longrightarrow\pi_0(P_1(\pi_{0/K}(H_*))^{\Gamma_K})$$
induced by the truncation morphism $\pi_{0/K}(H_*)\to P_1(\pi_{0/K}(H_*))$. By Theorem B on page 189 of \cite{Go} there is an obstruction class $\delta^{H_*}_1(a,b)\in H^2(K,\pi_2(\pi_{0/K}(H_*,z)))$. For every positive integer $n$ and $a,b$ as above let $\delta^{H_*}_1(a,b)_n$ denote the image of $\delta^{H_*}_1(a,b)$ under the composition of the natural map:
$$H^2(K,\pi_2(\pi_{0/K}(H_*,z)))\rightarrow 
H^2(K,\pi_2(\pi_{0/K}(H_*,z))/n\pi_2(\pi_{0/K}(H_*,z))),$$
the homomorphism
$$H^2(K,\pi_2(\pi_{0/K}(H_*,z))/n\pi_2(\pi_{0/K}(H_*,z)))
\longrightarrow H^2(K,H_2(\pi_{0/K}(H_*,z),\mathbb Z/n\mathbb Z))$$
induced by the Hurewitz map:
$$\pi_{0/K}(H_*,z))/n\pi_2(\pi_{0/K}(H_*,z))
\longrightarrow H_2(\pi_{0/K}(H_*,z),\mathbb Z/n\mathbb Z),$$
and the inverse of the isomorphism:
$$H^2(K,H_2(\pi_{0/K}(H_*),\mathbb Z/n\mathbb Z))
\longrightarrow H^2(K,H_2(\pi_{0/K}(H_*,z),\mathbb Z/n\mathbb Z))$$
induced by the weak equivalence $\pi_{0/K}(H_*)\to\pi_{0/K}(H_*,z)$. Finally for every $a\in X(hK)$ let $a^{H_*}\in\pi_0(\pi_{0/K}(H_*)^{\Gamma_K})$ denote its image under the canonical map $X(hK)\to\pi_0(\pi_{0/K}(H_*)^{\Gamma_K})$.

Now let us start the proof in earnest. Because $\alpha_n(c)=0$ there is a line bundle $\mathcal L$ on $X$ such that the image of its isomorphism class under the coboundary map $\textrm{Pic}(X)=H^1(X,\mathbb G_m)\rightarrow H^2(X,\mu_n)$ is $c$. Let $\pi:Y\rightarrow X$ denote the total space of $\mathcal L$ with the zero section removed; then $Y$ is a $\mathbb G_m$-torsor over $X$ whose class in $H^1(X,\mathbb G_m)$ is the same as the class of $\mathcal L$. It will be enough to show that 
$$\{\delta^{H_*}_1(\iota_{X/K}(x)^{H_*},\iota_{X/K}(y)^{H_*})_n,d\}_n$$
is zero for every \'etale hypercovering $H_*$ of $X$ and for every
$$d\in H^0(K,H^2(\pi_{0/K}(H_*),\mu_n))$$
whose image is $\phi_n(c)$ with respect to the homomorphism:
$$H^0(K,H^2(\pi_{0/K}(H_*),\mu_n))
\longrightarrow H^0(K,H^2(\overline X,\mu_n))$$
induced by the pull-back map $H^2(\pi_{0/K}(H_*),\mu_n)\to H^2(\overline X,\mu_n)$.

Fix such a hypercovering $H_*$ and let $\pi^*(H_*)$ denote the pull-back of $H_*$ onto $Y$ with respect to $\pi$. Note that $\phi_n(\pi^*(c))=0$, in fact even $\pi^*(c)=0$. Indeed the latter follows as the pull-back of the torsor $\pi:Y\rightarrow X$ onto $Y$ with respect to $\pi$ is trivial: the diagonal $Y\rightarrow Y\times_XY$ is a section. Let
$$d'\in H^0(K,H^2(\pi_{0/K}(\pi^*(H_*)),\mu_n))$$
be the image of $d$ with respect to the homomorphism:
$$H^0(K,H^2(\pi_{0/K}(H_*),\mu_n))
\longrightarrow H^0(K,H^2(\pi_{0/K}(\pi^*(H_*)),\mu_n))$$
induced by $\pi$. As the image of $d'$ is $\phi_n(\pi^*(c))$ under the
pull-back:
$$H^0(K,H^2(\pi_{0/K}(\pi^*(H_*)),\mu_n))\longrightarrow
H^0(K,H^2(\overline Y,\mu_n))$$
by naturality we get that there is a morphism $f:I_*\to\pi^*(H_*)$ of \'etale hypercoverings of $Y$ such that $f^*(d')\in H^0(K,H^2(\pi_{0/K}(I_*),\mu_n))$ is zero. 

Let $z\in Y(K)$ be arbitrary. (There are such points, for example the fibre above $x$ contains $K$-rational points.) Let $\overline z\in Y(\overline K)$ denote the geometric point lying above $z$, corresponding to the choice of algebraic closure $\overline K\supset K$. Note that the image of the canonical map from $\pi_1(\overline Y,\overline z)\to\pi_{0/K}(I_*,z)$, which is well-defined by Proposition \ref{truncationcomparison2}, is finite. Therefore by Remark \ref{pointsremark} there are two points $x'$, $y'$ in $Y(K)$ whose image under $\pi$ is $x$, $y$, respectively, such that $\iota_{X/K}(x')^{I_*},\iota_{X/K}(y')^{H_*}\in\pi_0(\pi_{0/K}(I_*)^{\Gamma_K})$ have the same image under the map
$$\pi_0(\pi_{0/K}(I_*)^{\Gamma_K})\longrightarrow
\pi_0(P_1(\pi_{0/K}(I_*))^{\Gamma_K})$$
induced by truncation. So by the above the obstruction class
$$\delta^{I_*}_1(\iota_{Y/K}(x')^{I_*},\iota_{Y/K}(y')^{I_*})\in H^2(K,\pi_2(\pi_{0/K}(I_*,z)))$$
is well-defined. By naturality of obstruction classes it will be enough to show that
$$\{\delta^{I_*}_1(\iota_{Y/K}(x')^{I_*},\iota_{Y/K}(y')^{I_*})_n,f^*(d')\}_n=0.$$
But this is clear since $f^*(d')$ is zero.
\end{proof}
\begin{lemma}\label{gabbertrick} Let $X$ be a geometrically irreducible smooth quasi-projective variety over $K$. For every $n\in\mathbb N$, for every $c\in H^2(X,\mu_n)$ and for every $x,y\in X(K)$ such that $\iota_{X/K}(x)\sim_1\iota_{X/K}(y)$ we have:
$$(\iota_{X/K}(x),c)_n=(\iota_{X/K}(y),c)_n
\quad\Rightarrow\quad
\{\delta^X_1(\iota_{X/K}(x),\iota_{X/K}(y))_n,\phi_n(c)\}_n=0.$$
\end{lemma}
\begin{proof} By a theorem of Gabber (see Theorem 1.1 of \cite{Jo}) there is an Azumaya algebra $A$ on $X$ which represents $\alpha_n(c)\in H^2(X,\mathbb G_m)$. Without the loss of generality we may assume that $A$ has rank $n$ by enlarging $n$ if it is necessary, since for every pair of positive integers $n|m$ the map $H^2(K,\mu_n)\rightarrow H^2(K,\mu_m)$ induced by the inclusion map $\mu_n\subseteq\mu_m$ is injective. Let $\pi:Y\rightarrow X$ be the $PGL_n$-torsor corresponding to $A$, that is, the torsor whose class in $H^1(X,PGL_n)$ is the same as the class of $A$. Let $\sigma\in H^1(K,PGL_n)$ be the class of the fibre of $Y$ over $x$ and let $\pi^{\sigma}:Y^{\sigma}\rightarrow X$ be the twist of $\pi:Y\rightarrow X$ by $\sigma$. Then the fibre of $Y^{\sigma}$ over $x$ is a trivial $PGL_n$-torsor. Because $(x,c)_n=(y,c)_n$ and the natural map $H^1(K,PGL_n)\rightarrow H^2(K,\mu_n)$ is injective we get that the fibre of $Y^{\sigma}$ over $y$ is also a trivial $PGL_n$-torsor. So there are points $x'$, $y'$ in $Y^{\sigma}(K)$ whose image under $\pi^{\sigma}$ is $x$, $y$, respectively. By Remark \ref{pointsremark} we may even assume that $\iota_{Y^{\sigma}/K}(x')\sim_1\iota_{Y^{\sigma}/K}(y')$. So by Lemma \ref{pairingnatural} and the naturality of obstruction classes it will be enough to show that
$$\{\delta^{Y^{\sigma}}_1(\iota_{Y^{\sigma}/K}(x'),\iota_{Y^{\sigma}/K}(y'))_n,
\phi_n((\pi^{\sigma})^*(c))\}_n=0.$$
In order to do so it will be enough to show that
$$\alpha_n((\pi^{\sigma})^*(c))=(\pi^{\sigma})^*(\alpha_n( c ))=0$$
by Lemma \ref{linebundle}. But the latter is clear since the pull-back of the torsor
$\pi:Y^{\sigma}\rightarrow X$ onto $Y^{\sigma}$ with respect to $\pi^{\sigma}$ is trivial: the diagonal $Y^{\sigma}\rightarrow Y^{\sigma}\times_XY^{\sigma}$ is a section.
\end{proof}
\begin{notn} For every quasi-projective variety $Y$ over $K$ consider the Hoch\-schild-Serre spectral sequence
$$E^2_{p,q}=H^p(K,H^q(\overline Y,\mu_n))\Rightarrow
H^{p+q}(Y,\mu_n).$$
Because $H^3(K,\mu_n)=0$ (for example by local Tate duality), the coboundary map:
$$d^3_{0,2}:E^3_{0,2}\longrightarrow E^3_{3,0}\subseteq H^3(K,\mu_n)$$
is zero and therefore:
$$E^{\infty}_{0,2}=E^3_{0,2}=\textrm{Ker}(d^2_{0,2}).$$
Therefore the spectral sequence furnishes an exact sequence:
$$\CD H^2(Y,\mu_n)@>\phi_n>>
H^0(K,H^2(\overline Y,\mu_n))@> d_n >> 
H^2(K,H^1(\overline Y,\mu_n))\endCD$$
which is functorial. Here $d_n$ is the coboundary map $d^2_{0,2}:E^2_{0,2}\to E^2_{2,1}$ and $\phi_n$ is the quotient map $H^2(Y,\mu_n)\to E^{\infty}_{0,2}\cong\textrm{Ker}(d^2_{0,2})$ by the highest step in the filtration on $H^2(Y,\mu_n)$ induced by the spectral sequence.
\end{notn}
\begin{lemma}\label{goodcover} Let $X$ be a geometrically irreducible smooth quasi-projective variety over $K$ and let $x,y\in X(K)$ be \'etale-Brauer equivalent. Then there is a connected finite \'etale cover $f:Y\rightarrow X$ such that
\begin{enumerate}
\item[$(i)$] there are $\widetilde x,\widetilde y\in Y(K)$ such that $f(\widetilde x)=x,f(\widetilde y)=y$, and $\widetilde x,\widetilde y$ are Brauer-equivalent,
\item[$(ii)$] for every $c\in H^0(K,H^2(\overline X,\mu_n))$ we have $d_n(f^*( c ))=0$. 
\end{enumerate}
\end{lemma}
\begin{proof} Let $\eta$ be a $\overline K$-valued point of $X$. Let $\mathcal K$ be the intersection of the kernels of all continuous homomorphisms
$\pi_1(\overline X,\eta)\rightarrow\mathbb Z/n\mathbb Z$. Because
$\pi_1(\overline X,\eta)$ is topologically finitely generated, the subgroup
$\mathcal K$ is open and characteristic. Let $\mathcal K'\subseteq
\pi_1(X,\eta)$ be the subgroup generated by $\mathcal K$ and by the image of an element of $s_{X/K}(x)$; this is an open subgroup. Let $f:Y\rightarrow X$ be the connected finite \'etale cover such that the image of $\pi_1(Y)$ with respect to $\pi_1(f)$ is $\mathcal K'$. By construction there is a $\widetilde x\in Y(K)$ such that $f(\widetilde x)=x$. Because $x$ and $y\in X(K)$ are \'etale-Brauer equivalent, there is a $\widetilde y\in Y(K)$ such that $f(\widetilde y)=y$, and $\widetilde x,\widetilde y$ are Brauer-equivalent. Recall that by the universal coefficient theorem:
$$\CD0\rightarrow\textrm{Ext}^1(H_n(V,\widehat{\mathbb Z}),\mathbb Z/n\mathbb Z) 
\rightarrow
H^n(V,\mathbb Z/n\mathbb Z)
@>e_n>>\textrm{Hom}(H_n(V,\widehat{\mathbb Z}),\mathbb Z/n\mathbb Z))
\rightarrow0\endCD$$
for every variety $V$ over $\overline K$, where the map $e_n$ is induced by the evaluation pairing and the Ext groups are for the category of pro-groups. In particular there is a natural isomorphism $H^1(V,\mathbb Z/n\mathbb Z)
\cong\textrm{Hom}(H_1(V,\widehat{\mathbb Z}),\mathbb Z/n\mathbb Z))$. Therefore the pull-back map $f^*:H^1(\overline X,\mathbb Z/n\mathbb Z)\rightarrow H^1(\overline Y,\mathbb Z/n\mathbb Z)$ is zero. Because the map $d_n$ is functorial, the claim now follows.
\end{proof}
Assume that $X(K)\neq\emptyset$. For every $x,y\in X(hK)$ such that $x\sim_1y$ let $\epsilon^X_1(x,y)$ denote the image of the obstruction class $\delta^X_1(x,y)$ under the homomorphism:
$$H_{*}:H^2(K,\pi_2(\overline X))
\longrightarrow H^2(K,H_2(\overline X,\widehat{\mathbb Z}))$$
induced by the Hurewitz map:
$$H:\pi_2(\overline X)
\longrightarrow H_2(\overline X,\widehat{\mathbb Z}).$$
\begin{prop}\label{cohomological} Let $X$ be a geometrically irreducible smooth quasi-projective variety over $K$ and let $x,y\in X(K)$ be \'etale-Brauer equivalent. Then we have:
$$\epsilon^X_1(\iota_{X/K}(x),\iota_{X/K}(y))=0.$$
\end{prop}
Note that the claim is meaningful because $\iota_{X/K}(x)\sim_1\iota_{X/K}(y)$ by Lemma \ref{5.6}.
\begin{proof} It will be enough to show that 
$$\delta^X_1(\iota_{X/K}(x),\iota_{X/K}(y))_n=0.$$
for every $n\in\mathbb N$. By the universal coefficient theorem quoted above and by local Tate duality every element of $H^2(K,H_2(\overline X,\mathbb Z/n\mathbb Z))$ annihilated by the pairing $\{\cdot,\cdot\}_n$ must be zero. So it will be enough to show that
$$\{\delta^X_1(\iota_{X/K}(x),\iota_{X/K}(y))_n,c\}_n=0$$
for every $c\in H^0(K,H^2(\overline X,\mu_n))$. Let $f:Y\rightarrow X$ and $\widetilde x,\widetilde y\in Y(K)$ be as in Lemma \ref{goodcover}. By Lemma \ref{pairingnatural} and the naturality of obstruction classes it will be enough to show that
$$\{\delta^Y_1(\iota_{Y/K}(x'),\iota_{Y/K}(y'))_n,f^*(c)\}_n=0$$
for every $c\in H^0(K,H^2(\overline X,\mu_n))$. Because $d_n(f^*( c ))=0$ this claim follows from Lemma \ref{gabbertrick}.
\end{proof}
The following result is Theorem \ref{etalebrauer} of the introduction.
\begin{thm} Let $K$ be a finite extension of $\mathbb Q_p$ and let $X$ be a smooth quasi-projective variety over $K$. Then \'etale--Brauer equivalence and $H$-equivalence coincide on $X(K)$.
\end{thm}
\begin{proof} Let $x,y\in X(K)$ be \'etale-Brauer equivalent. We need to show that they are $H$-equivalent. We may assume without the loss of generality that $X$ is geometrically irreducible. We already noted that $\iota_{X/K}(x)\sim_1\iota_{X/K}(y)$. Also note that by
Theorem \ref{obstructiontheory} it will be enough to show that $\delta^X_1(\iota_{X/K}(x),\iota_{X/K}(y))=0$ since the cohomological dimension of $K$ is two, so in this case the obstruction classes $\delta^X_n(\iota_{X/K}(x),\iota_{X/K}(y))$ will be zero for every $n\geq2$, too.

Fix an element $s$ of $s_{X/K}(x)$. For every open characteristic subgroup $\mathcal K$ of $\pi_1(\overline X,\eta)$ let $\mathcal K'\subseteq
\pi_1(X,\eta)$ be the subgroup generated by $\mathcal K$ and the image of $s$; this is an open subgroup. Moreover for every such $\mathcal K$ let $f_{\mathcal K}:Y_{\mathcal K}\rightarrow X$ be the connected finite \'etale cover such that the image of $\pi_1(Y)$ with respect to $(f_{\mathcal K})*$ is $\mathcal K'$. By construction there is a $x_{\mathcal K}\in Y_{\mathcal K}(K)$ such that $f_{\mathcal K}(x_{\mathcal K})=x$. Because $x$ and $y$ are \'etale-Brauer equivalent, there is a $y_{\mathcal K}\in Y_{\mathcal K}(K)$ such that $f_{\mathcal K}(y_{\mathcal K})=y$, and $x_{\mathcal K},y_{\mathcal K}$ are \'etale-Brauer equivalent. By the naturality of the obstruction classes the cohomology classes
$$\epsilon^{Y_{\mathcal K}}_1(x_{\mathcal K},y_{\mathcal K})\in H^2(K,H_2(
\overline Y_{\mathcal K},\widehat{\mathbb Z}))$$
furnish an element of $\lim_{(Y,y')\in\textrm{\rm Fet}(\overline X,x')}
H^2(K,H_2(Y,\widehat{\mathbb Z}))$, where $x'$ is a $\overline K$-valued point of $X$, which is the image of $\delta^X_1(x,y)$ with respect to the map:
$$H^2(K,\pi_2(\overline X))\cong
H^2(K,
\!\!\!\!\!\!\!\!\lim_{(Y,y')\in\textrm{\rm Fet}(\overline X,x')}\!\!\!\!\!\!\!\!
H_2(Y,\widehat{\mathbb Z}))
\longrightarrow\!\!\!\!\!\!\!\!
\lim_{(Y,y')\in\textrm{\rm Fet}(\overline X,x')}\!\!\!\!\!\!\!\!
H^2(K,H_2(Y,\widehat{\mathbb Z}))$$
furnished by the map $b_X\circ a_X$ which is an isomorphism by Theorem \ref{myhurewitz}. By Proposition \ref{cohomological} the classes $\epsilon^{Y_{\mathcal K}}_1(x_{\mathcal K},y_{\mathcal K})$ are zero, and hence the theorem holds.
\end{proof}

\section{$H$-equivalence over the real number field}

\begin{defn} Let $X$ be any scheme. Recall that a quadratic space over $X$ is a vector bundle $\mathcal E$ over $X$, that is, a locally free $\mathcal O_X$-module of finite rank, and an isomorphism $h:\mathcal E\rightarrow
\mathcal E^*$, where $\mathcal E^*$ denotes the dual of $\mathcal E$, which is symmetric, that is, the composition
$$\CD\mathcal E@>>>\mathcal E^{**}@>h^*>>\mathcal E^*\endCD$$
is equal to $h$, where the first map is the natural isomorphism of $\mathcal E$ with its bidual, and $h^*$ is the dual of $h$. In the special case when $X= \textrm{Spec}(K)$, where $K$ is a field, this concept is the same as the concept of a non-degenerate quadratic form over $K$.
\end{defn}
\begin{defn} Consider the case when $K=\mathbb R$. By Sylvester's theorem every non-degenerate quadratic form $q$ over $\mathbb R$ is isomorphic to a diagonal form:
$$\langle\underbrace{1,1,\ldots,1}_m,\underbrace{-1,\ldots,-1}_n\rangle,$$
and the ordered pair $(m,n)$ only depends on the isomorphism class of $q$. Let $\rho(q)=m+n$ denote the rank of $q$, and let $\sigma(q)=m-n$ denote the signature of $q$, respectively. By the above non-degenerate forms over $\mathbb R$ are classified by their rank and signature. Let $X$ be a smooth variety over $\mathbb R$, let $U\subset X(\mathbb R)$ be a connected component, and let $q=(\mathcal E,h)$ be a quadratic space over $X$. Then for every $x\in U$ the pull-back $x^*(q)$ has the same signature, which we will call the signature of $q$ on $U$. 
\end{defn}
\begin{thm}[Mah\'e, Houdebine-Mah\'e] Let $X$ be a smooth variety over $\mathbb R$ which is either affine or projective. Let $U,V\subset X(\mathbb R)$ be two different connected components. Then there is a quadratic space $q=(\mathcal E,h)$ over $X$ such that the signature of $q$ on $U$ is zero, and the signature of $q$ on $V$ is non-zero. 
\end{thm}
\begin{proof} See Theorem 1.1.1 in \cite{Mah} in the affine case, and the main result of \cite{HM} for the projective case.
\end{proof}
\begin{defn} Let $K$ be a field whose characteristic is not $2$. Every non-degenerate quadratic form $q$ is isomorphic to a diagonal form:
$$\langle a_1,a_2,\ldots,a_n\rangle.$$
The Stiefel-Whitney classes of the form $q$ above are defined as the cup product (see \cite{Mil}):
$$w(q)=1+w_1(q)+\cdots+w_n(q)=(1+\delta(a_1))(1+\delta(a_2))
\cdots(1+\delta(a_n)),$$
where $w_i(q)\in H^i(K,\mathbb Z/2\mathbb Z)$, and
$$\delta:K^*\longrightarrow H^1(K,\mathbb Z/2\mathbb Z)$$
is the boundary map of the Kummer exact sequence:
$$\CD0@>>>\mathbb Z/2\mathbb Z@>>>
\mathbb G_m@>{x\mapsto x^2}>>\mathbb G_m@>>>0.\endCD$$
The Stiefel-Whitney classes are independent of the diagonalisation of $q$. 
\end{defn}
\begin{rem}\label{computation} Assume again that $K$ is the real number field. Then as a graded algebra:
$$H^*(\mathbb R,\mathbb Z/2\mathbb Z)\cong\mathbb F_2[t],$$
where $t$ is the generator of the group $H^1(\mathbb R,\mathbb Z/2\mathbb Z)$ of order two. Let $q$ be a non-degenerate quadratic form $q$ over $\mathbb R$ isomorphic to a diagonal form:
$$\langle\underbrace{1,1,\ldots,1}_m,\underbrace{-1,\ldots,-1}_n\rangle,$$
then by the above:
$$w(q)=(1+t)^n=1+nt+\cdots+t^n.$$
\end{rem} 
In \cite{EKV} the authors construct Stiefel-Whitney classes for any quadratic space $q=(\mathcal E,h)$ over a $\mathbb Z[1/2]$-scheme $X$, which lives in mod 2 \'etale cohomology:
$$w_i(q)\in H^i(X,\mathbb Z/2\mathbb Z),$$
is functorial over the category of $\mathbb Z[1/2]$-schemes, and specialises to the construction above when $X$ is the spectrum of a field. We will use these classes to separate connected components of real points of varieties defined over $\mathbb R$.
\begin{prop}\label{separation} Let $X$ be a smooth variety over $\mathbb R$ which is either affine or projective. Let $U,V\subset X(\mathbb R)$ be two different connected components. Then there is a natural number $i$ and a cohomology class $c\in H^i(X,\mathbb Z/2\mathbb Z)$ over $X$ such that for
every $x\in U$ the pull-back $x^*( c )\in H^i(\mathbb R,\mathbb Z/2\mathbb Z)$ is zero, and for every $x\in V$ the pull-back $x^*( c )\in H^i(\mathbb R,\mathbb Z/2\mathbb Z)$ is non-zero. 
\end{prop}
\begin{proof} By Theorem \ref{complex} and Proposition \ref{connectedcomponents} there is a $c\in H^0(X,\mathbb Z/2\mathbb Z)$ such that $x^*( c )\in H^0(\mathbb R,\mathbb Z/2\mathbb Z)$ is zero for every $x\in U$ and $y^*( c )\in H^0(\mathbb R,\mathbb Z/2\mathbb Z)$ is non-zero for every $y\in V$ if $U$ and $V$ lie on two different connected components of $X(\mathbb C)$. Therefore we may assume that $X$ is geometrically connected without the loss of generality. Let $q=(\mathcal E,h)$ be a quadratic space over $X$ such that the signature of $q$ on $U$ is zero, and the signature of $q$ on $V$ is non-zero. We may assume that the signature of $q$ on $V$ is negative by taking $(\mathcal E,-h)$ instead, if it is necessary. Because $X$ is connected, the rank of the vector bundle $\mathcal E$ is constant on $X$. This rank is even, say $2m$, because the signature of $q$ on $U$ is zero. Then the signature of $q$ on $V$ is $2m-2n$, where $n$ is a positive integer bigger than $m$. By Remark \ref{computation} above we have $x^*(w_n(q))=
w_n(x^*(q))=0$ for every $x\in U$, and $x^*(w_n(q))=w_n(x^*(q))=
t^n$ for every $x\in V$. The claim follows. 
\end{proof}
\begin{defn} For every morphism of sites $m:\mathcal C\rightarrow\mathcal C'$ let $m^*:\mathcal C'\rightarrow\mathcal C$ denote the functor underlying $m$. Let $\mathcal C$ be a Grothendieck site. A left action $\alpha$ of a group $\Gamma$ on $\mathcal C$ is a morphism $\alpha(g):\mathcal C\rightarrow\mathcal C$ of sites for each $g\in\Gamma$ such that $\alpha(1)$ is the identity map of $\mathcal C$ and $\alpha(gh)=\alpha(g)\circ\alpha(h)$ for every $g,h\in\Gamma$. When $\Gamma$ is profinite we say that the action $\alpha$ is continuous if for every morphism $h:U\rightarrow V$ in $\mathcal C$ there is an open subgroup $\Delta$ of $\Gamma$ such that
${\alpha(g)}^*(U)=U, {\alpha(g)}^*(V)=V$ and
${\alpha(g)}^*(h)=h$ for every $g\in\Delta$. Assume now that $\Gamma$ is a profinite group. By a $\Gamma$-site
$(\mathcal C,\alpha)$ we mean a Grothendieck site $\mathcal C$ with a continuous left action $\alpha$ of $\Gamma$ on $\mathcal C$. As usual we will drop $\alpha$ from the notation whenever this is convenient. 
\end{defn}
\begin{example}\label{example1} A basic example of a $\Gamma$-site is the Grothendieck site $\Gamma-Sets$, where the coverings are surjective maps, equipped with the left action $\alpha$ such that ${\alpha(g)}^*(U)=U$ for every object $U$ of $\Gamma-Sets$ and for every $g\in\Gamma$, and for every morphism $h:U\rightarrow V$ and $g\in\Gamma$ the map ${\alpha(g)}^*(h):U\rightarrow V$ is given by the rule $x\mapsto gh(x)$. By a slight abuse of notation we will let $\Gamma-Sets$ denote this $\Gamma$-site, too.
\end{example}
\begin{defn} Let $(\mathcal C,\alpha)$ be a $\Gamma$-site. A $\Gamma$-invariant object of $(\mathcal C,\alpha)$ is an object $U$ of $\mathcal C$ such that ${\alpha(g)}^*(U)=U$ for every $g\in\Gamma$. A $\Gamma$-equivariant morphism of $(\mathcal C,\alpha)$ is a morphism $h:U\rightarrow V$ of $\mathcal C$ such that $U$ and $V$ are $\Gamma$-invariant objects and
$h\circ{\alpha(g)}^*(\textrm{id}_U)=
{\alpha(g)}^*(\textrm{id}_V)\circ h$ for every $g\in\Gamma$. Let $\mathcal C^{\Gamma}$ denote category whose objects are $\Gamma$-invariant objects of $\mathcal C$ and whose morphisms are $\Gamma$-equivariant maps between these. Since the composition of $\Gamma$-equivariant morphisms are $\Gamma$-equivariant, with these morphisms $\mathcal C^{\Gamma}$ is indeed a category. Let $\mathbf T$ denote the Grothendieck topology of the site $\mathcal C$, that is, for every object $U$ of $\mathcal C$ let $\mathbf T(U)$ denote the collection of covering sieves of $U$. We say that a sieve $S$ on $U\in\textrm{ob}(\mathcal C^{\Gamma})$ is $\Gamma$-invariant if for every $V\in\textrm{ob}(\mathcal C)$ and for every $h\in S(V)$ there is a $W\in\textrm{ob}(\mathcal C^{\Gamma})$, a morphism $h'\in S(W)$ which is $\Gamma$-equivariant, and a morphism $h'':V\to W$ such that $h=h'\circ h''$. For every $\Gamma$-invariant $S$ as above let $S^{\Gamma}$ denote the sieve on $U$ in the category $\mathcal C^{\Gamma}$ given by the rule:
$$S^{\Gamma}(V)=S(V)\cap\textrm{Hom}_{\mathcal C^{\Gamma}}(V,U).$$
For every $U\in\textrm{ob}(\mathcal C^{\Gamma})$ let $\mathbf T^{\Gamma}(U)$ denote the following collection of sieves $S$ in the category $\mathcal C^{\Gamma}$:
$$\mathbf T^{\Gamma}(U)=\{S^{\Gamma}|\textrm{ $S$ is in $\mathbf T(U)$ and it is $\Gamma$-invariant}\}.$$
\end{defn}
\begin{example}\label{example2} Let $K$ be a field and let $\Gamma_K= \textrm{Gal}(\overline K|K)$ denote the absolute Galois group of $K$ as above. Let $X$ be a locally Noetherian scheme over $K$ and let $\mathcal C$ denote the small \'etale site of the base change of $X$ to $K$. Then $\mathcal C$ is naturally equipped with the structure of a $\Gamma_K$-site, induced by the action of $\Gamma_K$ on $\overline K$. By \'etale decent the category $\mathcal C^{\Gamma_K}$ is equivalent to the small \'etale category of $X$ and $\mathbf T^{\Gamma_K}$ is the \'etale topology of $X$ on it. In particular $\mathbf T^{\Gamma_K}$ is a Grothendieck topology. 
\end{example}
\begin{defn} Assume now that $\mathcal C$ satisfies the conditions in chapters 8 and 9 of \cite{AM} and $\mathbf T^{\Gamma}$ is a Grothendieck topology on $\mathcal C^{\Gamma}$. In particular we suppose that $\mathcal C$ is closed under finite coproducts and it is locally connected in the sense of 9.1--9.2 of \cite{AM}. Let $U$ be a $\Gamma$-invariant object of $(\mathcal C,\alpha)$. Then $\Gamma$ acts on the set $\pi_0(U)$ of connected components of $U$ and this action makes $\pi_0(U)$ into a $\Gamma$-set. Let $\pi_{\Gamma}(U)$ denote this $\Gamma$-set. We say that a simplicial object of $\mathcal C$ is $\Gamma$-invariant if it is a simplicial object of $\mathcal C^{\Gamma}$.  If $X_*$ is a $\Gamma$-invariant simplicial object of $(\mathcal C,\alpha)$ then the face and degeneracy maps of $X_*$ induce $\Gamma$-equivariant maps between the $\Gamma$-set $\pi_{\Gamma}(X_n)$ which makes the collection
$\{\pi_{\Gamma}(X_n)\}_{n=1}^{\infty}$ into an object of $\Gamma-SSets$ which we will denote by $\pi_{\Gamma}(X_*)$. Since the site $(\mathcal C^{\Gamma},\mathbf T^{\Gamma})$ inherits the good properties of the site $\mathcal C$, we may apply part $(i)$ of Corollary 8.13 of \cite{AM} on page 105 to conclude that the functor
$$X_*\mapsto \pi_{\Gamma}(X_*)$$
above induces an object of $\textrm{Pro}-\textrm{Ho}(\Gamma-SSets)$. We will call the later the $\Gamma$-equivariant homotopy type of $\mathcal C$ and denote it by $\Pi(\mathcal C)$.
\end{defn}
\begin{rem} Let $K,\Gamma_K$ and $X$ be as in Example \ref{example2}, and let $\mathcal C$ denote the small \'etale site of the base change of $X$ to $\overline K$, as above. Then the $\Gamma_K$-equivariant homotopy type of $\mathcal C$ is just the relative \'etale homotopy type $Et_{/K}(X)$ of $X$ as defined by Harpaz and Schlank.
\end{rem}
\begin{example}\label{example3} Assume now that $\Gamma$ is a finite group and let $X$ be a locally connected, Hausdorff, paracompact topological space equipped with a continuous left $\Gamma$-action. Let $\mathcal C$ be the ordinary site on the coproducts of open subsets of $X$. Then $\mathcal C$ is naturally equipped with the structure of a $\Gamma$-site, induced by the action of $\Gamma$ on $X$. Moreover the quotient $\Gamma\backslash X$ of $X$ by the action of $\Gamma$ is also a locally connected, Hausdorff, paracompact topological space, and the category $\mathcal C^{\Gamma}$ is the category of coproducts of open subsets of $\Gamma\backslash X$ and $\mathbf T^{\Gamma}$ is the ordinary site of $\Gamma\backslash X$ on it. In particular $\mathbf T^{\Gamma}$ is a Grothendieck topology. 
\end{example}
\begin{defn} Let $\Gamma$ and $X$ be as above. We say that $X$ is
$\Gamma$-contractible if the is a subgroup $\Delta$ such that the coset $\Delta\backslash\Gamma$, equipped with the discrete topology and the natural left $\Gamma$-action, is $\Gamma$-equivariantly homotopy-equivalent to $X$.  Assume now that every open subset of $X$ is paracompact, and that $X$ is locally $\Gamma$-equivariantly contractible, that is, for every (finite) orbit $O\in X$ and for every $\Gamma$-invariant open $U\subseteq X$ containing $O$ there an open $\Gamma$-invariant and
$\Gamma$-contractible $V\subseteq X$ such that $V\subseteq U$ and $O\subseteq V$. 
\end{defn}
Note that when $\Gamma$ is finite then $\Gamma-SSets$ is just the usual category of simplicial sets with a $\Gamma$-action, moreover Goerss's notion of homotopy fixed point spaces coincides with the usual one. Also note that the singular complex $S_*X$ of $X$ is equipped with an action of
$\Gamma$ which makes it an object of $\Gamma-SSets$.
\begin{thm} Let $U_*$ be a $\Gamma$-invariant hypercovering of $\mathcal C$ such that every $\Gamma$-orbit of connected components of every $U_n$ is $\Gamma$-contractible. Then the simplicial $\Gamma$-set $\pi_{\Gamma}(U_*)$ is isomorphic to the simplicial $\Gamma$-set $S_*X$ in $\textrm{\rm Ho}(\Gamma-SSets)$.
\end{thm}
\begin{proof} Let $S_*U_n$ denote the singular complex of $U_n$. Then $S_*U_*$ is a bisimplicial object in $\Gamma-Sets$. We denote by $(DU)_*$ its diagonal simplicial $\Gamma$-set $(DU)_n=S_nU_n$. Then we have obvious maps of simplicial $\Gamma$-sets:
$$\xymatrixcolsep{5pc}\xymatrix{
&(DU)_*\ar[ld]_{\alpha} \ar[rd]^{\beta} & \\
\pi_{\Gamma}(U_*) & & S_*X,}$$
and we claim that these two maps are homotopy equivalences in $\textrm{\rm Ho}(\Gamma-SSets)$, which will prove the theorem.

For every subgroup $\Delta\leq\Gamma$ and for every $\Gamma$-set $Y$  let $Y^{\Delta}$ denote the subset of $Y$ fixed by $\Delta$. Similarly for every simplicial $\Gamma$-set $Y_*$ let $Y^{\Delta}_*$ denote the simplicial set 
$\{Y_n^{\Delta}\}_{n=1}^{\infty}$ such that the face and degeneracy maps are the restrictions of such maps of the simplicial set $Y_*$. By the definition of the strict model structure we need to show that the maps:
$$\xymatrixcolsep{5pc}\xymatrix{
&(DU)_*^{\Delta}\ar[ld]_{\alpha|_{(DU)_*^{\Delta}}} 
\ar[rd]^{\beta|_{(DU)_*^{\Delta}}} & \\
\pi_{\Gamma}(U_*)^{\Delta} & &(S_*X)^{\Delta}}$$
 of simplicial sets are homotopy equivalences in $\textrm{\rm Ho}(SSets)$ for every $\Delta$ as above.

Note that for every $\Gamma$-invariant hypercovering $V_*$ of $X$ the simplicial object $V_*^{\Gamma}$ in the category of disjoint union of open sets of the closed subspace $X^{\Delta}$ of $X$ is a hypercovering, too, since it is
the pull-back of the hypercovering $V_*$ onto $X^{\Delta}$ with respect to the
inclusion map. Because the $U_n$ are $\Gamma$-contractible we have $\pi_{\Gamma}(U_*)^{\Delta}=\pi_0(U_*^{\Delta})$. Moreover every $\Delta$-invariant singular simplex of $X$ must lie in $X^{\Delta}$, so $(S_*X)^{\Delta}$ is the singular complex $S_*X^{\Delta}$ of $X^{\Delta}$. Similarly $S_*U^{\Delta}_n$ denote the singular complex of $U_n^{\Delta}$. Then $S_*U^{\Delta}_*$ is a bisimplicial set. We denote by
$(DU^{\Delta})_*$ its diagonal simplicial set $(DU^{\Delta})_n=S_nU^{\Delta}_n$. Then we have $(DU)_*^{\Delta}=(DU^{\Delta})_*$, so we only need to show that the analogues: 
$$\xymatrixcolsep{5pc}\xymatrix{
&(DU^{\Delta})_*\ar[ld]_{\alpha^{\Delta}} \ar[rd]^{\beta^{\Delta}} & \\
\pi_0(U_*^{\Delta}) & & S_*X^{\Delta}}$$
of the maps $\alpha$ and $\beta$ for the topological space $X^{\Delta}$ are homotopy equivalences in $\textrm{\rm Ho}(SSets)$ for every $\Delta$ as above. Since every $\Gamma$-orbit of connected components of every $U_n$ is $\Gamma$-contractible, the connected components of $U_n^{\Delta}$ are contractible. Similarly $X^{\Delta}$ is locally contractible, since $X$ is locally
$\Gamma$-contractible. Because $X^{\Delta}$ is a closed subspace of a paracompact topological space, it is also paracompact. The claim now follows from Theorem 12.1 of \cite{AM} on page 129.
\end{proof}
\begin{cor}\label{11.16} Let $\mathcal C$ be the ordinary site on the coproducts of open subsets of a Hausdorff topological space $X$ equipped with a continuous left $\Gamma$-action. Assume that every open subset of $X$ is paracompact, and that $X$ is locally $\Gamma$-equivariantly contractible. Then the pro-object $\Pi(\mathcal C)$ is canonically isomorphic to the element $S_*X$ in $\textrm{\rm Pro}-\textrm{\rm Ho}(\Gamma-SSets)$.
\end{cor}
\begin{proof} Because those $\Gamma$-invariant hypercoverings $U_*$ of $\mathcal C$ such that every connected component of every $U_n$ is $\Gamma$-contractible are cofinal by assumption, the claim follows immediately.
\end{proof}
\begin{defn} Let again $\Gamma$ be an arbitrary profinite group. A morphism $m:(\mathcal C,\alpha)\rightarrow(\mathcal C',\alpha')$ of $\Gamma$-sites is a morphism of sites $m:\mathcal C\rightarrow\mathcal C'$ such that for every $g\in\Gamma$ and for every morphism $h:U\rightarrow V$ of $\mathcal C'$ we have ${\alpha}(g)^*(m^*(U))=m({\alpha'(g)}^*(U)),{\alpha(g)}^*(m^*(V))=m^*({\alpha'(g)}^*(V))$ and ${\alpha(g)}^*(m^*(h))=m^*({\alpha'(g)}^*(h))$. For every such $m$ the underlying functor $m^*$ carries $\Gamma$-invariant hypercoverings to $\Gamma$-invariant hypercoverings, and so it furnishes a map
$$\Pi(m):\Pi(\mathcal C)\longrightarrow\Pi(\mathcal C')$$
in $\textrm{\rm Pro}-\textrm{\rm Ho}(\Gamma-SSets)$ (when these are defined). The map $\Pi(m)$ in turn induces a map:
$$\Pi(m)(E\Gamma):\Pi(\mathcal C)(E\Gamma)
\longrightarrow\Pi(\mathcal C')(E\Gamma)$$
of homotopy fixed points. 
\end{defn}
\begin{defn}\label{11.18b} A $\Gamma$-invariant point (or more conveniently $\Gamma$-point) of a $\Gamma$-site $(\mathcal C,\alpha)$ is a morphism $p:\Gamma-Sets\rightarrow\mathcal C^{\Gamma}$ of $\Gamma$-sites. Note that the composition of $p^*$ and the functor $\Gamma-Sets\rightarrow Sets$ forgetting the $\Gamma$-action is a point of the site $\mathcal C$ in the sense of 8.3 of \cite{AM}, which perhaps justifies our terminology. We will let $\mathcal C(\Gamma)$ denote the set of $\Gamma$-points of $(\mathcal C,\alpha)$. Note that the homotopy type $\Pi(\Gamma-Sets)$ is contractible (this is clear from Corollary \ref{11.16}, too), so the set $\Pi(\Gamma-Sets)(E\Gamma)$ has one element. For every $p\in\mathcal C(\Gamma)$ let $\iota_{\mathcal C}(p)\in\Pi(\mathcal C)(E\Gamma)$ denote the image of $\Pi(\Gamma-Sets)(E\Gamma)$ with respect to $\Pi(p)$. Clearly
$$\iota_{\mathcal C}:\mathcal C(\Gamma)\longrightarrow
\Pi(\Gamma-Sets)(E\Gamma)$$
is a natural transformation.
\end{defn}
\begin{example} Let $K$, $\Gamma_K= \textrm{Gal}(\overline K|K)$, $X$ and
$\mathcal C$ be as in Example \ref{example2}. Since the small \'etale site of Spec$(\overline K)$ is isomorphic to $\Gamma_{\Gamma}-Sets$ as a
$\Gamma_K$-site, every $K$-valued point of $X$ supplies a $\Gamma_K$-point of the site $\mathcal C$. Therefore the map $\iota_{\mathcal C}$ introduced above is a generalisation of the map $\iota_{X/K}$. Similarly when $\Gamma$ is a finite group, $X$ is a locally connected, Hausdorff, paracompact topological space equipped with a continuous left $\Gamma$-action and $\mathcal C$ is as in Example \ref{example3} then every point of $X$ fixed by $\Gamma$ furnishes a $\Gamma$-point of the site $\mathcal C$, and hence the restriction of $\iota_{\mathcal C}$ onto $X^{\Gamma}$ is a map $X^{\Gamma}\to\Pi(\mathcal C)(E\Gamma)$.
\end{example}
\begin{prop}\label{11.18} Let $K$ be the real number field $\mathbb R$ and let $X$ be a variety over $K$. Then two $K$-rational points of $X$ are $H$-equivalent if they are in the same connected component of the topological space $X(K)$.
\end{prop}
\begin{proof} Let $\Gamma=\Gamma_K=\textrm{Gal}(\overline K|K)$ be the group of two elements, and let $\mathcal C'$ denote the small \'etale site of the base change of $X$ to $\overline K$. Moreover let $\mathcal C$ be the ordinary site on the coproducts of open subsets of $X(\overline K)$ with respect to its usual topology. In addition to these $\Gamma$-sites we also introduce the $\Gamma$-site $\mathcal C''$ whose objects are topological spaces $X'$ lying over the topological space of $X(\overline K)$ such that the map $X' \rightarrow X$ is a local isomorphism, i.e., that every point $x\in X'$ has a neighborhood which is isomorphic onto its image. Since any \'etale map of schemes $X' \rightarrow X$ over $\overline K$ is a local isomorphism on the underlying topological spaces, and since an open set is in $\mathcal C''$, we have morphisms of $\Gamma$-sites:
$$\xymatrixcolsep{5pc}\xymatrix{
&\mathcal C''\ar[ld]_{m} \ar[rd]^{m'} & \\
\mathcal C & & \mathcal C'.}$$
Now it is clear from the definition of local isomorphisms that every $\Gamma$-invariant hypercovering of $\mathcal C''$ is dominated by a $\Gamma$-invariant hypercovering of $\mathcal C$. Thus the map $\Pi(m):
\Pi(\mathcal C'')\rightarrow\Pi(\mathcal C)$ is a homotopy equivalence in $\textrm{\rm Pro}-\textrm{\rm Ho}(\Gamma-SSets)$, and so
$$\Pi(m)(E\Gamma):\Pi(\mathcal C'')(E\Gamma)\rightarrow
\Pi(\mathcal C)(E\Gamma)$$
is a bijection. Therefore by the naturality of the maps $\iota_{\mathcal C},\iota_{\mathcal C'}$ and $\iota_{\mathcal C''}$ it will be enough to show that for every pair $x,y$ of $K$-rational points of $X$ lying in the same connected component of the topological space $X(K)$ we have $\iota_{\mathcal C}(x)=
\iota_{\mathcal C}(y)$. Let $f:[0,1]\to X(K)$ be a continuous path connecting $x$ with $y$, that is we have $f(0)=x$ and $f(1)=y$. By naturality again it will be enough to show that $\iota_{\mathcal D}(0)=
\iota_{\mathcal D}(1)$ where $\mathcal D$ is the ordinary site on the coproducts of open subsets of $[0,1]$ with respect to its usual topology, equipped with the trivial $\Gamma$-action. But the interval $[0,1]$ is contractible to a point
$\Gamma$-equivariantly, so $\Pi(\mathcal D)(E\Gamma)$ is a one element set by Corollary \ref{11.16}.
\end{proof}
\begin{rem} It is not difficult to push the arguments of this proof a little bit further to prove an equivariant version of Theorem 12.9 of \cite{AM} on page 142, using an equivariant analogue of the profinite completion functor, but we will not pursue this further, since it would take us too far away from our original project. However in a forthcoming publication we will in fact prove such a claim in a much more general context.  
\end{rem}
The following result is Theorem \ref{real} of the introduction.
\begin{thm} Let $K$ be the real number field $\mathbb R$ and let $X$ be a smooth affine or projective variety over $K$. Then two $K$-rational points of $X$ are $H$-equivalent if and only if they are in the same connected component of the topological space $X(K)$.
\end{thm}
\begin{proof} By Lemma \ref{basicinvariant} and Proposition \ref{separation} two $H$-equivalent real points of $X$ must be in the same component. On the other hand by Proposition \ref{11.18} two real points of $X$ in the same connected component must be $H$-equivalent.
\end{proof}

\section{The Homotopy Section Principle}

\begin{notn} Recall that for every field of characteristic zero the topological Gal$(\overline K|K)$-module $\widehat{\mathbb Z}(1)$ is defined as the projective limit $\varprojlim_{n\in\mathbb N}\mu_n$ where the directed set structure on $\mathbb N$ is furnished by divisibility and for every $m,n\in\mathbb N$ such that $m|n$ the transition map $\mu_n\rightarrow\mu_m$ is multiplication by $n/m$. For every number field $K$ let $|K|$ denote the set of places of $K$, and for every $v\in|K|$ let $K_v$ denote the completion of $K$ with respect to $v$. For every $v\in|K|$ fix an embedding $j_v:\overline K\rightarrow\overline K_v$ of $K$-extensions. For every $k\in\mathbb N$ and for every discrete Gal$(\overline K|K)$-module $M$ let $\Sh^k(K,M)$ denote the subgroup:
$$\Sh^k(K,M)=\textrm{Ker}\left(
\prod_{v\in|K|}j_{v*}:H^k(K,M)\rightarrow\prod_{v\in|K|}H^k(
K_v,M)\right)$$
of $H^k(K,M)$ where $j_{v*}$ denotes the restriction
map induced by $j_v$ for every $v\in|K|$. For every topological Gal$(\overline K|K)$-module $M$ which is a projective limit $\varprojlim_{i\in I}M_i$ of discrete Gal$(\overline K|K)$-modules let $\Sh^k(K,M)$ denote $\varprojlim_{i\in I}\Sh^k(K,M_i)$.
\end{notn}
\begin{lemma}\label{sha2} We have $\Sh^2(K,\widehat{\mathbb Z}(1))=0$ for every number field $K$.
\end{lemma}
\begin{proof} By part $(a)$ of Theorem 4.10 of \cite{Mi2} on page 70 the group $\Sh^2(K,\widehat{\mathbb Z}(1))$ is zero if and only if $\Sh^1(K,\mathbb Q/\mathbb Z)$ is, where we equip $\mathbb Q/\mathbb Z$ with the discrete topology and the trivial Gal$(\overline K|K)$-action. We may identify $\Sh^1(K,\mathbb Q/\mathbb Z)$ with the kernel of the map:
$$\textrm{Hom}(\Gamma_K,\mathbb Q/\mathbb Z)\longrightarrow\prod_{v\in|K|}\textrm{Hom}(\Gamma_{K_v},\mathbb Q/\mathbb Z)$$
furnished by restriction onto the family of subgroups $\Gamma_{K_v}\ (\forall v\in|K|)$ of $\Gamma_K$. The claim now follows from the Chebotarev density theorem.
\end{proof}
\begin{prop}\label{no3dim} Assume that $K$ is a number field and $X(K)$ is non-empty. For every $x,y\in X(hK)$ such that $x\sim_2y$ we have:
$$x=y\ \Leftrightarrow\ r_{v*}(x)=r_{v*}(y)
\quad(\textrm{$v\in|K|$, $v$ is real}).$$
\end{prop}
\begin{proof} The product of the restriction maps:
$$\prod_{\substack{v\in|K|\\\textrm{$v$ is real}}}
j_{v*}:H^n(K,M)\longrightarrow
\prod_{\substack{v\in|K|\\\textrm{$v$ is real}}}
H^n(K_v,M)$$ 
is injective for every integer $n\geq3$ and for every discrete Galois module $M$ over $K$ by part $(c)$ of Theorem 4.10 of \cite{Mi2} on page 70. Hence for every $n\geq2$ and for every pair of sections $x,y\in X(hK)$ such that $x\sim_ny$ we get that
$$\delta^X_n(x,y)=0\ \Leftrightarrow
\ \delta^{X_v}_n(r_{v*}(x),r_{v*}(y))=0
\quad(\textrm{$v\in|K|$, $v$ is real})$$
from the naturality of the obstruction classes. The claim now follows from Theorem \ref{obstructiontheory}.
\end{proof}
\begin{defn} Let $K$ be for a moment any field of characteristic zero. Recall that two points $x,y\in X(K)$ are called directly $R$-equivalent if there is a rational map $f:\mathbb P_K^1\dashrightarrow X$ of $K$-varieties such that $f(0)=x$ and $f(\infty)=y$. The $R$-equivalence on $X(K)$ is the equivalence relation generated by direct $R$-equivalence. Let $X(K)/R$ denote the equivalence classes of this relation. Note that $A$-equivalence coincides with $R$-equivalence when $X$ is projective by the valuative criterion of properness. In this case let 
$$\iota_{X/K,R}:X(K)/R\longrightarrow X(hK)$$
be the map furnished by Proposition \ref{a1}. 
\end{defn}
\begin{notn} It is particularly interesting to study $X(K)/R$ through the map $\iota_{X/K,R}$  when $K$ is a number field. For every variety $X$ defined over $K$ and for every $v\in|K|$ let $X_v$ denote the base change of $X$ to Spec$(K_v)$. For every $v\in|K|$ the embedding $j_v$ furnishes a map:
$$r_{v*}:X(hK)\longrightarrow X_v(hK_v)$$
by functoriality. Let $\mathbb A_K=\prod_{v\in|K|}'K_v$ denote the ring of ad\`eles of $K$. Let $X(h\mathbb A_K)$ denote the image of $X(\mathbb A_K)$ with respect to the map
$$\prod_{v\in|K|}\iota_{X_v/K_v}:X_v(K_v)\longrightarrow\prod_{v\in|K|}
X_v(hK_v).$$
We define the set $\textrm{Sel}(X/K)$ of Selmer homotopy fixed points of $X$ to be
$$\textrm{Sel}(X/K)=\big(\prod_{v\in|K|}r_{v*}\big)^{-1}
\big(X(h\mathbb A_k)\big)\subseteq X(hK).$$
\end{notn}
We are interested in the following natural generalisation of the Shafarevich-Tate conjecture:
\begin{prin}\label{hsp} Assume that $X$ is smooth and projective. Then the map:
\begin{equation}\label{HSP}
\CD X(K)/R@>\iota_{X/K,R}>>\textrm{\rm Sel}(X/K)
\endCD
\end{equation}
is injective and its image is dense with respect to the pro-discrete topology of $\textrm{\rm Sel}(X/K)$.
\end{prin}
The claim above is obviously true if $X$ does not have local points everywhere, and hence HSP should be considered as a new form of the local-global principle. The next proposition shows that HSP is indeed a generalisation of standard conjectures of this sort:
\begin{prop}\label{simpleexamples} Let $K$ be a number field.
\begin{enumerate}
\item[$(a)$] HSP holds for Brauer--Severi varieties and for non-singular quadratic hypersurfaces $H\subset\mathbb P^n_K$ of positive dimension.
\item[$(b)$] Let $X$ be a smooth, geometrically connected projective curve $X$ of genus at least two over $K$. Then HSP holds for $X$ if and only if a weak (local-global) form of Grothendieck's section conjecture holds for $X$. 
\item[$(c)$] Let $X$ be an abelian variety over $K$. Then HSP holds for $X$ if and only if the Shafarevich-Tate conjecture holds for $X$. 
\end{enumerate}
\end{prop}
\begin{proof} First assume that $X$ is either a Brauer--Severi variety or a non-singular quadratic hypersurface of positive dimension. When $\textrm{\rm Sel}(X/K)$ is empty there is nothing to prove. Assume now that $\textrm{\rm Sel}(X/K)$ is non-empty: then $X(\mathbb A_K)$ is non-empty, too. Because the local-global principle holds for $X$ we get that $X(K)$ is also non-empty. In this case $X(K)/R$ consists of one element, and hence it will be enough to show that $\textrm{\rm Sel}(X/K)$ also has one element. Let $x,y\in\textrm{\rm Sel}(X/K)$. Because for every $v\in|K|$ the set $X_v(K_v)/R$ has one element, we get that $r_{v*}(x)=r_{v*}(y)$ for every such $v$.

By Corollary 12.13 of \cite{AM} on page 144 we know that $X_{\overline K}$ is weakly  homotopy equivalent to $X_{\mathbb C}$. Because $\pi_1(X(\mathbb C))=\{1\}$ (either because $X_{\mathbb C}$ is isomorphic to $\mathbb P^n_{\mathbb C}$ for $n=\dim(X)$ or by the Lefschetz hyperplane section theorem), we get $\pi_1(X_{\overline K})=\{1\}$ by Theorem \ref{complex}. Therefore $x\sim_1y$. Moreover
$$\pi_2(X_{\overline K})=
H_2(X_{\overline K},\widehat{\mathbb Z})=
\textrm{Hom}(H^2(X_{\overline K},\widehat{\mathbb Z}),\widehat{\mathbb Z})
$$
by Corollary 6.2 of \cite{AM} on page 70. When $X$ is a Brauer--Severi variety, since it has a rational point it is isomorphic to $\mathbb P^n_K$. Therefore $H^2(X_{\overline K},\widehat{\mathbb Z})\cong
\widehat{\mathbb Z}(-1)$, and hence $\pi_2(X_{\overline K})\cong
\widehat{\mathbb Z}(1)$. Because $r_{v*}(x)=r_{v*}(y)$ for every $v\in|K|$, we have $\delta^X_2(x,y)\in\Sh^2(K,\pi_2(X_{\overline K}))$, so this obstruction class vanishes by Lemma \ref{sha2}. So we get that $x\sim_2y$.

When $X$ is a quadratic hypersurface of dimension at least $3$, its embedding $X\rightarrow\mathbb P^{n+1}_K$ as such a hypersurface induces an isomorphism:
$$H^2(X_{\overline K},\widehat{\mathbb Z})\cong
H^2(\mathbb P^{n+1}_{\overline K},\widehat{\mathbb Z}),\widehat{\mathbb Z})
\cong
\widehat{\mathbb Z}(-1)$$
by the Lefschetz hyperplane section theorem, and hence we may conclude as above that $x\sim_2y$. The only remaining case is of a quadratic surface. In this case either
$$H^2(X_{\overline K},\widehat{\mathbb Z})
\cong\widehat{\mathbb Z}(-1)\oplus\widehat{\mathbb Z}(-1),$$
when both pencils of lines on $X$ are defined over $K$, or it is the induction of $\widehat{\mathbb Z}(-1)$ from a quadratic extension of $K$, otherwise. Clearly in the first case the group $\Sh^2(K,\pi_2(X_{\overline K}))$ still vanishes by Lemma \ref{sha2}, while in the second case this claim follows from Shapiro's lemma and Lemma \ref{sha2}. Again we get that $x\sim_2y$. Claim $(a)$ now follows from Proposition \ref{no3dim}.

Assume now that $X$ is a smooth, geometrically connected projective curve of genus at least two over $K$. Then there is a commutative diagram:
$$\CD X(K)@>s_{X/K}>>\textrm{Sec}(X/K)\\
@VVV@VVV\\
\prod_{v\in|K|}X_v(K_v)@>\prod_{v\in|K|}s_{X_v/K_v}>>
\prod_{v\in|K|}\textrm{Sec}(X_v/K_v)\endCD$$ 
where the vertical maps are the obvious maps. The weak local-global form of Grothendieck's section conjecture asserts that the diagram above is cartesian. We also know that in this case $s_{X/K}$ is injective and by Faltings's theorem $X(K)$ is finite. In particular in this case the map in (\ref{HSP}) has dense image if and only if it is surjective. 

Let's assume first that the weak local-global form of Grothendieck's section conjecture holds for $X$, and show that HSP holds for $X$. This is trivial when
$\textrm{Sel}(X/K)$ is empty, so we may assume that $\textrm{Sel}(X/K)$ is non-empty. Then $\textrm{Sec}(X/K)$ has an element whose image in $\prod_{v\in|K|}\textrm{Sec}(X_v/K_v)$ lies in the image of $\prod_{v\in|K|}s_{X_v/K_v}$, and hence $X(K)$ is non-empty, by our assumption. The claim now follows from Corollary \ref{curvesandabelians} and Lemma \ref{sectionanabelian}. Let us prove the converse. We may assume that $\textrm{Sec}(X/K)$ has an element whose image in $\prod_{v\in|K|}\textrm{Sec}(X_v/K_v)$ lies in the image of $\prod_{v\in|K|}s_{X_v/K_v}$ without the loss of generality. By the main theorem of Harpaz and Schlank (see \cite{HS}) for smooth projective varieties the set Sel$(X/K)$ is non-empty if and only if the \'etale-Brauer set of $X$ is. In our case the latter is non-empty (see \cite{St}), so we get that Sel$(X/K)$ is non-empty, too. So by our hypothesis $X(K)$ is non-empty, so the claim follows from Corollary \ref{curvesandabelians} and Lemma \ref{sectionanabelian}. Claim $(b)$ is settled.

Finally consider the case when $X$ is an abelian variety over $K$. By Corollary \ref{curvesandabelians} and Lemma \ref{sectionanabelian} the map $j_{X/K}$ is a bijection. Moreover there is a natural bijection 
\begin{equation}\label{4.9.1}
\textrm{\rm Sec}(X/K)\cong H^1(K,\prod_{\textrm{$l$ is prime}}T_l(X))
\end{equation}
where $T_l(X)$ denotes the $l$-th Tate module of $X$, and under this identification $s_{X/K}$ corresponds to the coboundary map furnished by Kummer theory. In particular $\iota_{X/K}$ is injective.

The image of $\textrm{Sel}(X/K)\subseteq X(hK)$ with respect to the composition of $j_{X/K}$ and the isomorphism of (\ref{4.9.1}) is $\textrm{Sel}(K,X)^{tor}$ where $\textrm{Sel}(K,X)$ is the Selmer group of $X$ over $K$. The quotient of $\textrm{Sel}(K,X)^{tor}$ by the closure of the image of $X(K)$ under the coboundary map is $\Sh(K,X)^{tor}$ where $\Sh(K,X)$ is the Tate-Shafarevich group of $X$ over $K$. Since the group $\Sh(K,X)[n]$ is finite for every positive integer $n$ the group $\Sh(K,X)^{tor}$ is trivial if and only if $\Sh(K,X)$ is finite. So claim $(c)$ holds.
\end{proof}
\begin{rem} It is tempting to believe that HSP should hold for every smooth projective variety because of its very general form, but this is not true. The fundamental reason is the Harpaz--Schlank theorem quoted above which implies that if HSP holds for $X$ then the Brauer-Manin obstruction applied to \'etale covers is the only obstruction for the Hasse principle. Since now there are many counter-examples to this claim (see \cite{Po}, \cite{HSk} and \cite{CTPS}) we get that there are two- and three-dimensional counterexamples to HSP.  However we can offer some positive results; see Theorems \ref{bloch}, \ref{birational} and \ref{chatelet} in the next two sections.
\end{rem}

\section{Geometrically rational and birational surfaces}

\begin{lemma}\label{sha2finite} Let $X$ be a smooth, projective, geometrically rational surface defined over a number field $K$. Then the group
$\Sh^2(K,\pi_2(\overline X))$ is finite.
\end{lemma}
\begin{proof} There is a finite Galois extension $L|K$ such
that the action of Gal$(\overline L|L)$ on Pic$(\overline X)$ is trivial.
Hence $\Sh^1(L,\textrm{Pic}(\overline X)\otimes\mathbb Q/\mathbb Z)=0$ by the Chebotarev density theorem (see the proof of Lemma \ref{sha2}), where we equip
$\mathbb Q/\mathbb Z$ with the discrete topology and the trivial Gal$(\overline K|K)$-action. As the image of $\Sh^1(K,\textrm{Pic}(\overline X)\otimes
\mathbb Q/\mathbb Z)$ lies in $\Sh^1(L,\textrm{Pic}(\overline X)\otimes
\mathbb Q/\mathbb Z)$ with respect to
the restriction map:
$$H^1(K,\textrm{Pic}(\overline X)\otimes\mathbb Q/\mathbb Z)\rightarrow
H^1(L,\textrm{Pic}(\overline X)\otimes\mathbb Q/\mathbb Z),$$
the group $\Sh^1(K,\textrm{Pic}(\overline X)\otimes\mathbb Q/\mathbb Z)$ lies in the image of the inflation map:
$$H^1(\textrm{Gal}(L|K),
\textrm{Pic}(\overline X)\otimes\mathbb Q/\mathbb Z)\rightarrow H^1(K,
\textrm{Pic}(\overline X)\otimes\mathbb Q/\mathbb Z).$$
Since $H^1(\textrm{Gal}(L|K),\textrm{Pic}(\overline X)\otimes\mathbb Q/\mathbb Z)$ is finite we get that $\Sh^1(K,\textrm{Pic}(\overline X)\otimes\mathbb Q/\mathbb Z)$ is finite, too. We have $\pi_1(\overline X)=\pi_1(\mathbb P^2_{\overline K})=\{1\}$ because the \'etale fundamental group is a birational invariant. Hence
$$\pi_2(\overline X)=H_2(\overline X,\widehat{\mathbb Z})=
\textrm{Hom}(H^2(\overline X,\mathbb Q/\mathbb Z),\mathbb Q/\mathbb Z)
=\textrm{Hom}(\textrm{Pic}(\overline X)\otimes\mathbb Q/\mathbb Z,
\mathbb Q/\mathbb Z\otimes_{\widehat{\mathbb Z}}\widehat{\mathbb Z}(1))$$
by Corollary 6.2 of \cite{AM} on page 70 and the fact that for a geometrically
rational surface $X$ the \'etale Chern class map
$$c_1:\textrm{Pic}(\overline X)\otimes\mathbb Q/\mathbb Z
\longrightarrow
H^2(\overline X,
\mathbb Q/\mathbb Z\otimes_{\widehat{\mathbb Z}}\widehat{\mathbb Z}(1))$$
is an isomorphism. Hence by part $(a)$ of Theorem 4.10 of \cite{Mi2} on page 70 there is a perfect duality between $\Sh^2(K,\pi_2(\overline X))$ and $\Sh^1(K,\textrm{Pic}(\overline X)\otimes\mathbb Q/\mathbb Z)$. The claim is now clear.
\end{proof}
Let $\textrm{\rm Sel}_0(s)\subseteq 
\textrm{\rm Sel}(X/K)$ denote the pre-image of $\prod_{v\in|K|}r_{v*}(s)$ with respect to the map $\prod_{v\in|K|}r_{v*}$ for every $s\in\textrm{\rm Sel}(X/K)$.
\begin{lemma}\label{sha_0bound} Let $X$ be a smooth, projective, geometrically rational surface defined over a number field $K$ such that $X(K)\neq\emptyset$ and let $s\in\textrm{\rm Sel}(X/K)$. Then the cardinality of $\textrm{\rm Sel}_0(s)$ is at most the order of $\Sh^2(K,\pi_2(\overline X))$.
\end{lemma}
\begin{proof} By Proposition \ref{no3dim} it will be enough to show that the number of equivalence classes of the relation $\sim_2$ in $\textrm{\rm Sel}_0(s)$ is at most  the order of $\Sh^2(K,\pi_2(\overline X))$. For every $x,y\in\textrm{\rm Sel}_0(s)$ at least we have $x\sim_1y$ by Theorem \ref{obstructiontheory}, since $\pi_1(\overline X)=\{1\}$, and so $|\textrm{Sec}(X/K)|=1$. Moreover $\delta_1^X(x,y)\in\Sh^2(K,\pi_2(\overline X))$ by the naturality of obstruction classes. But
$$\delta_1^X(x,y)=\delta_1^X(x,s)-\delta_1^X(y,s)$$
so the claim follows from the pigeonhole principle.
\end{proof}
\begin{thm}\label{bloch} Let $X$ be a smooth, projective, geometrically rational surface defined over a number field $K$ such that $X(K)\neq\emptyset$. Then
$\textrm{\rm Sel}(X/K)$ is finite.
\end{thm}
Because we expect that $X(K)/R$ is finite for such an $X$ and $K$, this result should be also expected, assuming that HSP holds for $X$.
\begin{proof} Because for every archimedean $v\in|K|$ the topological space $X_v(K_v)$ has only finitely many connected components the set $X_v(K_v)/H$ is finite for such $v$ by Theorem \ref{real}. Since $\pi_1(\overline X)=\{1\}$ we get that Brauer equivalence and $H$-equivalence coincide on $X_v(K_v)$ for every non-archimedean $v\in|K|$ by Theorem \ref{etalebrauer}. Hence by Corollary A.2 of \cite{Bl} on page 55 the set $X_v(K_v)/H$ is finite for every $v\in|K|$. Moreover $X_v(K_v)/H$ has at most one element when $v$ is non-archimedean and $X$ has good reduction at $v$ by Corollary A.3 of \cite{Bl} on page 55. Hence we may conclude that the set $\prod_{v\in|K|}X_v(K_v)/H$ is finite. So we only need to show that the map:
$$\prod_{v\in|K|}r_{v*}:
\textrm{\rm Sel}(X/K)\longrightarrow
\prod_{v\in|K|} X_v(K_v)/H$$ 
is finite to one. This follows from Lemmas \ref{sha2finite} and \ref{sha_0bound}.
\end{proof}
\begin{lemma}\label{blowup1} Let $X$ and $Y$ be smooth projective surfaces over a field $K$ and let $\pi:X\rightarrow Y$ be a composition of monoidal transformations over $K$. Then the map $\pi_*:X(K)/R\rightarrow Y(K)/R$ induced by $\pi$ is a bijection.
\end{lemma}
\begin{proof} We may immediately reduce to the case when $\pi$ is the blow-up of an irreducible subvariety $S\subset Y$ of dimension zero by induction on the number of blow-ups in some sequence of contractions $X \rightarrow X_1 \rightarrow\cdots\rightarrow X_n=Y$ whose composition is $\pi$. If $S$ has no $K$-valued points then $\pi^{-1}(S)$ has no $K$-valued points either and the map $\pi_*$ is obviously a bijection. Otherwise $S$ consists of one $K$-valued point. In this case $\pi^{-1}(S)$ is isomorphic to $\mathbb P^1_K$ and the claim is clear.
\end{proof}
\begin{lemma}\label{blowup2} Let $X$ and $Y$ be smooth projective surfaces over the field $\mathbb R$ and let $\pi:X\rightarrow Y$ be a composition of monoidal transformations over $\mathbb R$. Then the map $\pi_*:\pi_0(X(\mathbb R))\rightarrow\pi_0(Y(\mathbb R))$ induced by $\pi$ is a bijection.
\end{lemma}
\begin{proof} The argument is the same as above.
\end{proof}
\begin{prop}\label{blowup3} Let $X$ and $Y$ be smooth geometrically irreducible projective surfaces over a number field $K$ and let $\pi:X\rightarrow Y$ be a composition of monoidal transformations over $K$. Assume that $Y$ is simply connected and $Y(K)\neq\emptyset$. Then the map $\pi_*:\textrm{\rm Sel}(X/K)\rightarrow\textrm{\rm Sel}(Y/K)$ induced by $\pi$ is injective.
\end{prop}
\begin{proof} Note that $X(K)\neq\emptyset$, since $X$ is birational to $Y$. Again we may assume without the loss of generality that $\pi$ is the blow-up of an irreducible subvariety of dimension zero. Let $x,y\in\textrm{Sel}(X/K)$ be such that $\pi_*(x)=\pi_*(y)$. Because the map induced by $\pi$ between the fundamental groups is an isomorphism, we get that $x\sim_1y$. Also note that by Theorem \ref{obstructiontheory} it will be enough to show that $\delta^X_1(x,y)=0$ since by Lemma \ref{blowup2} and Theorem \ref{real} we have $r_{v*}(x)=r_{v*}(y)$ for every real place $v$ of $K$, so the higher obstruction classes $\delta^X_n(x,y)$ will vanish for every $n\geq2$ by Proposition \ref{no3dim}. 

By the Hurewitz theorem:
$$\pi_2(\overline X)\cong
H_2(\overline X,\widehat{\mathbb Z})\cong
H_2(\overline Y,\widehat{\mathbb Z})\oplus 
\textrm{Ker}(H_2(\pi))\cong\pi_2(\overline Y)\oplus\textrm{Ker}(H_2(\pi)),$$
where
$$
H_2(\pi):H_2(\overline X,\widehat{\mathbb Z})
\longrightarrow H_2(\overline Y,\widehat{\mathbb Z})
$$
is induced by $\pi$. Because $\pi$ is the blow-up of a single closed point, the module $\textrm{Ker}(H_2(\pi))$ is the induction of $\widehat{\mathbb Z}(1)$ from a finite extension of $K$. Therefore $\Sh^2(K,\textrm{Ker}(H_2(\pi))$ vanishes by Lemma \ref{sha2} and Shapiro's lemma. Therefore $\delta^X_1(x,y)$ is zero by the naturality of obstruction classes.
\end{proof}
\begin{thm}\label{birational} Let $\pi:X\rightarrow Y$ be a composition of monoidal transformations between geometrically irreducible smooth projective surfaces over $K$. Assume that $\textrm{\rm Sel}(Y/K)$ is finite, $Y$ is simply connected, the set $Y(K)$ is non-empty, and HSP holds for $Y$ over $K$. Then $\textrm{\rm Sel}(X/K)$ is finite and HSP holds for $X$ over $K$, too.
\end{thm}
This result can be used to supply many examples of surfaces satisfying HSP, for example blow-ups of Ch\^atelet surfaces; see Theorem \ref{chatelet} below.
\begin{proof} By Proposition \ref{blowup3} the set
$\textrm{\rm Sel}(X/K)$ injects into $\textrm{\rm Sel}(Y/K)$, so it is finite. Let $x,y\in X(K)$ be $H$-equivalent. Then $\pi(x)$ and $\pi(y)$ are $H$-equivalent elements of $Y(K)$, so they are $R$-equivalent, since HSP holds for $Y$. Therefore $x$ and $y$ are also $R$-equivalent by Lemma \ref{blowup1}. We get that the map
\begin{equation}\label{8.4.1}
\iota_{X/K}:X(K)/R\longrightarrow\textrm{\rm Sel}(X/K)
\end{equation}
is injective. Let $s$ be an element of $\textrm{\rm Sel}(X/K)$. Because
$\textrm{\rm Sel}(Y/K)$ is finite, the topology on it is discrete. Therefore there is a $y\in Y(K)$ such that $\iota_{Y/K}(y)=\pi_*(s)$. Let $x\in X(K)$ be such that $\pi(x)=y$. Then $\pi_*(\iota_{X/K}(x))=\iota_{Y/K}(y)=\pi_*(s)$, so by Proposition \ref{blowup3} we get that $\iota_{X/K}(x)=s$, so the map in (\ref{8.4.1}) is also surjective.
\end{proof}

\section{Generalised Ch\^atelet surfaces}

\begin{notn} For every torus $S$ defined over a field $K$ of characteristic zero let $C(S)$ denote the Gal$(\overline K|K)$-module of its cocharacters and for every scheme $X$ over $K$ let
$$\delta:H^1(X,S)\rightarrow H^2(X,C(S)\otimes\widehat{\mathbb Z}(1))$$
be the projective limit of the coboundary maps
$$\delta_n:H^1(X,S)\rightarrow H^2(X,S[n])$$
furnished by the corresponding Kummer exact sequences where $S[n]$ denotes the $n$-torsion subgroup scheme of $S$. Note that $\delta$ maps $\Sh^1(K,S)$ into $\Sh^2(K,C(S)\otimes\widehat{\mathbb Z}(1))$ when $K$ is a number field. Let
$$\delta_0:\Sh^1(K,S)\rightarrow\Sh^2(K,C(S)
\otimes\widehat{\mathbb Z}(1))$$ 
be the restriction of $\delta$ onto $\Sh^1(K,S)$.
\end{notn}
\begin{lemma}\label{ordersofshas} Let $S$ be a torus defined over a field $K$ of characteristic zero. Then the following hold:
\begin{enumerate}
\item[$(i)$] the map $\delta:H^1(K,S)\rightarrow H^2(K,C(S)
\otimes\widehat{\mathbb Z}(1))$ is injective,
\item[$(ii)$] the map $\delta_0:\Sh^1(K,S)\rightarrow\Sh^2(K,C(S)
\otimes\widehat{\mathbb Z}(1))$ is bijective when $K$ is a number field. 
\end{enumerate}
\end{lemma}
\begin{proof} First note that for every field $K$ of characteristic
zero and for every torus $S$ defined over $K$ the cohomology group $H^1(K,S)$
has finite exponent. In fact there is a finite Galois extension $L|K$ such
that the action of Gal$(\overline L|L)$ on $C(S)$ is trivial. Hence $H^1(L,S)=0$ by Hilbert's theorem 90. Therefore $H^1(K,S)$ is the image of  
$H^1(\textrm{Gal}(L|K),S)$ under the inflation map. But the group $H^1(\textrm{Gal}(L|K),S)$
is annihilated by the order of Gal$(L|K)$ so the first claim is now clear. 

Now we prove that $\delta_0$ is also surjective when $K$ is a number field. Let $L|K$ be a finite Galois extension of the type as above and assume that the degree of this extension is $m$. By the Albert--Brauer--Hasse--Noether theorem the group $\Sh^2(L,S)$ is trivial, and hence $\Sh^2(K,S)$ is annihilated by multiplication by $m$, since it is a sub-quotient of $H^2(\textrm{Gal}(L|K),S)$. The cokernel of the restriction
\begin{equation}\label{hahaha}
\delta_n|_{\Shi^1(X,S)}:\Sh^1(X,S)\rightarrow\Sh^2(X,S[n])
\end{equation}
of $\delta_n$ onto $\Sh^1(X,S)$ is a subgroup of $\Sh^2(X,S)$. Therefore the map in (\ref{hahaha}) surjects onto $m\Sh^2(X,S[n])$ by the above. Claim $(ii)$ now follows by taking the limit.
\end{proof}
\begin{defn} We say that a smooth geometrically irreducible
projective surface $X$ defined over a field $K$ of characteristic zero is a generalised Ch\^atelet surface if there is an $a\in K^*$ and a separable polynomial $P\in K[x]$ of degree $4$ such that $X$ is a smooth compactification of the affine surface given by the equation:
$$y^2-az^2=P(x)$$
over $K$. For the sake of brevity we will frequently drop the adjective 'generalised' when we talk about generalised Ch\^atelet surfaces. 
\end{defn}
\begin{prop}\label{approximate} Let $X$ be a Ch\^atelet surface defined over a number field $K$. Then for every $s\in\textrm{\rm Sel}(X/K)$ there is an $x\in X(K)$ such that $\iota_{X_v/K_v}(i_v(x))=
r_{v*}(s)$ for every $v\in|K|$.
\end{prop}
\begin{proof} Let $S$ be a finite subset of $|K|$ which contains every archimedean place of $K$ and every non-archimedean place of $K$ where $X$ does not have good reduction. For every $v\in|K|$ choose an $x_v\in X_v(K_v)$ such that $\iota_{X_v/K_v}(x_v)=r_{v*}(s)$. By part $(b)$ of Theorem 8.6 of \cite{CTSSD2} on page 87 for every $v\in|K|$ there is an open neighbourhood $U_v\subseteq X_v(K_v)$ of $x_v$ contained by the $R$-equivalence class of $x_v$ in $X_v(K_v)$. Because $\prod_{v\in|K|}x_v\in X(\mathbb A_K)^{\textrm{Br}}$ by the easy direction of the Harpaz-Schlank theorem, there is an $x\in X(K)$ such that $i_v(x)\in U_v$ for every $v\in S$ by part $(c)$ of Theorem 8.11 of \cite{CTSSD2} on page 92. Clearly $\iota_{X_v/K_v}(i_v(x))=r_{v*}(s)$ for every $v\in S$. Because for every $v\in|K|-S$ the set $X(K)/R$ consists of one element by part $(c)$ of Theorem 8.6 of \cite{CTSSD2} on page 87 we get that $\iota_{X_v/K_v}(i_v(x))=r_{v*}(s)$ for every $v\in|K|-S$ as well.
\end{proof}
\begin{defn} For every smooth projective geometrically irreducible variety $X$ defined over a field $K$ let $CH_0(X)$ denote the Chow group of zero dimensional cycles on $X$ and let $A_0(X)$ denote the kernel of the degree map $\deg:CH_0(X)\rightarrow\mathbb Z$. Fix a point $x\in X(K)$. Then there is a map $\Psi_x:X(K)/R\rightarrow A_0(X)$ which for every $y\in X(K)$ maps the $R$-equivalence class of $y$ to $[y]-[x]$. When $K$ is a number field let $\Sh^1A_0(X)$ denote the subgroup of those elements $c$ of $A_0(X)$ such that the base change of $c$ to $K_v$ is the zero element of $A_0(X_v)$ for every $v\in|K|$. Now let $X$ be a Ch\^atelet surface, let $S$ be a torus over $K$ whose group of characters is isomorphic to Pic$(\overline X)$ as a Gal$(\overline K|K)$-module and let $\mathcal T$ be a universal torsor over $X$. Then we have a map:
$$\rho_{\mathcal T}:X(K)\longrightarrow H^1(K,S)$$
which associates to every $P\in X(K)$ the class of the fibre of $\mathcal T$ at $P$. Moreover there is a unique homomorphism $\Phi_{\mathcal T}:A_0(X)\rightarrow H^1(K,S)$ such that for every $y\in X(K)$ we have $\Phi_{\mathcal T}([y]-[x])=\rho_{\mathcal T}(y)-\rho_{\mathcal T}(x)$ (see page 88 of \cite{CTSSD2}). In particular the map $\rho_{\mathcal T}$ factors through $R$-equivalence; let $\rho_{\mathcal T,R}:X(K)/R\rightarrow
H^1(K,S)$ be the map which sends the $R$-equivalence class of every $y\in X(K)$ to $\rho_{\mathcal T}(y)$. When $K$ is a number field $\Phi_{\mathcal T}$ maps $\Sh^1A_0(X)$ into $\Sh^1(K,S)$. Let $\Phi_0:\Sh^1A_0(X)\rightarrow \Sh^1(K,S)$ be the restriction of $\Phi_{\mathcal T}$ onto $\Sh^1A_0(X)$.
\end{defn}
\begin{thm}[Colliot-Th\'el\`ene--Sansuc--Swinnerton-Dyer]\label{CTSSW} Let $X$ be a Ch\^atelet surface, let $\mathcal T$ be a universal torsor over $X$ and let $x\in X(K)$. Then the following hold:
\begin{enumerate}
\item[$(i)$]  the map $\rho_{\mathcal T,R}$ is an injection,
\item[$(ii)$] the map $\Psi_x$ is a bijection,
\item[$(iii)$] when $K$ is a number field the map $\Phi_0$ is a bijection.
\end{enumerate}
\end{thm}
\begin{proof} The map $\rho_{\mathcal T,R}$ is injective by part $(a)$ of Theorem 8.5 of \cite{CTSSD2} on page 86. Claim $(ii)$ is true by Theorem 8.8 of \cite{CTSSD2} on page 89 while claim $(iii)$ holds by Theorem 8.10 of \cite{CTSSD2} on page 91.
\end{proof}
Let $X$ be a Ch\^atelet surface over a number field  $K$. Then for every $x\in X(K)$ let $R_0(x)\subseteq X(K)/R$ denote the set of those $R$-equivalences classes $s$ such that for every $y\in s$ the points $x,y\in X_v(K_v)$ are $R$-equivalent for every $v\in|K|$. 
\begin{cor}\label{orderofr0} For every $X$, $K$ and $x$ as above the set $R_0(x)$ has the  same cardinality as $\Sh^1(K,S)$.
\end{cor}
\begin{proof} Note that $\Psi_x$ maps $R_0(x)$ into $\Sh^1A_0(X)$. Therefore it will be enough to show that the induced map
\begin{equation}\label{9.7.1}
\Psi_x|_{R_0(x)}:R_0(x)\longrightarrow\Sh^1A_0(X)
\end{equation}
is a bijection by part $(iii)$ of Theorem \ref{CTSSW}. It is injective by part $(ii)$ of Theorem \ref{CTSSW}. Let $s\in\Sh^1A_0(X)$ be arbitrary; by part $(ii)$ of Theorem \ref{CTSSW} there is a $y\in X(K)$ such that $[y]-[x]$ is $s$. Because $\Psi_x$ is a bijection over $K_v$ for every $v\in|K|$ by part $(ii)$ of 
Theorem \ref{CTSSW} we get that the points $x,y\in X_v(K_v)$ are $R$-equivalent for every $v\in|K|$. Therefore the map in (\ref{9.7.1}) is surjective, too.
\end{proof}
\begin{thm}\label{chatelet} Let $K$ be a number field and let $X$ be a generalised Ch\^atelet surface over $K$. Then HSP holds for $X$.
\end{thm}
\begin{rem} According to Remark 8.10.2 of \cite{CTSSD1} on page 91 there is a Ch\^atelet surface $X$ defined over a number field $K$ on which $R$-equivalence is strictly finer than Brauer equivalence. Hence by Theorem \ref{chatelet} we get that $H$-equivalence is finer than \'etale--Brauer equivalence over number fields. Moreover there are two rational points $x,y\in X(K)$ such that 
$x,y\in X_v(K_v)$ are $H$-equivalent for every $v\in|K|$, but $x$ and $y$ are not $H$-equivalent over $K$. The theorem above is also interesting because it covers a whole class of varieties $X$ for which HSP holds, but which are not homogeneous spaces, moreover every $R$-equivalence class of $X(K)$ is Zariski dense (see part $(b)$ of Theorem 8.5 of \cite{CTSSD2} on page 86) and, by Theorem 8.13 of \cite{CTSSD2} on page 95, the set $X(K)/R$ could be arbitrarily large for $X$ defined over $\mathbb Q$.
\end{rem}
\begin{proof}[Proof of Theorem \ref{chatelet}] Because Ch\^atelet surfaces are geometrically rational it will be both necessary and sufficient to show that the map:
$$\iota_{X/K,R}:X(K)/R\longrightarrow
\textrm{\rm Sel}(X/K)$$
is a bijection by Theorem \ref{bloch}. Let $c\in H^1(X,S)$ be the cohomology class corresponding to the universal torsor $\pi:\mathcal T\rightarrow X$. Let $x,y\in X(K)$ be in two different $R$-equivalence classes. By part $(i)$ of Theorem \ref{CTSSW} the pull-back classes $x^*(c ),y^*(c )\in H^1(K,S)$ are different. Therefore the classes
$$\delta(x^*(c ))=x^*(\delta( c  )),\ \ 
\delta(y^*(c ))=y^*(\delta( c ))\in H^2(K,C(S)\otimes\widehat{\mathbb Z}(1))$$
are also different by part $(i)$ of Lemma \ref{ordersofshas}. Therefore $x$ and $y$ are not $H$-equivalent by Lemma \ref{basicinvariant}. The injectivity of $\iota_{X/K,R}$ follows.

Now we prove that it is surjective, too. Let $s\in\textrm{\rm Sel}(X/K)$. By Proposition \ref{approximate} there is an $x\in X(K)$ such that $\iota_{X_v/K_v}(i_v(x))=r_{v*}(s)$ for every $v\in|K|$. Note that
$\iota_{X/K,R}$ maps $R_0(x)$ into $\textrm{\rm Sel}_0(s)$. Therefore it will be enough to show that the induced map
\begin{equation}\label{9.7.2}
\iota_{X/K,R}|_{R_0(x)}:R_0(x)\longrightarrow\textrm{\rm Sel}_0(s)
\end{equation}
is a bijection. By the above this map is injective, so it will be enough to show that the cardinality of $\textrm{\rm Sel}_0(s)$ is at most $\Sh^2(K,C(S)
\otimes\widehat{\mathbb Z}(1))$ by part $(ii)$ of Lemma \ref{ordersofshas} and Corollary \ref{orderofr0}. By the definition of $S$ we have that
$$C(S)\otimes\widehat{\mathbb Z}(1)=
\textrm{Hom}(\textrm{Pic}(\overline X),\mathbb Z)\otimes
\widehat{\mathbb Z}(1)=
\textrm{Hom}(\textrm{Pic}(\overline X)\otimes\widehat{\mathbb Z},\widehat{\mathbb Z}(1)).$$
We already noted in the proof of Lemma \ref{sha2finite} that
$$\pi_2(\overline X)=\textrm{Hom}(\textrm{Pic}(\overline X)\otimes\widehat{\mathbb Z},\widehat{\mathbb Z}(1))
\textrm{, so }
\pi_2(\overline X)=C(S)\otimes\widehat{\mathbb Z}(1).$$
Therefore it will be enough to show that the cardinality of $\textrm{\rm Sel}_0(s)$ is at most $\Sh^2(K,\pi_2(\overline X))$. But this is the content of Lemma \ref{sha_0bound}.
\end{proof}

\end{document}